\documentclass[]{article}

\PassOptionsToPackage{%
	maxalphanames	= 4,
	style			= alphabetic%
}{biblatex}

\usepackage[doi=false,isbn=false,url=false,
	hyperref=auto,
	sorting=nyt,
	maxnames=10,
	maxcitenames=4,
	backend=biber,
	texencoding=auto,
	giveninits=true,
	block=space,
]{biblatex}

\DefineBibliographyStrings{english}{%
	andothers = {et al}
}

\DeclareFieldFormat{eprint:arxiv}{%
	\href{https://arxiv.org/abs/#1}{\emph{arXiv:\allowbreak #1}}}

\DeclareFieldFormat{eprint:mrnumber}{%
	\href{http://www.ams.org/mathscinet-getitem?mr=MR#1}{MR\allowbreak #1}}

\DeclareFieldFormat{eprint:jstor}{%
	\href{http://www.jstor.org/stable/#1}{JSTOR\allowbreak #1}}

\DeclareFieldFormat{eprint:customeprint}{%
	\em Available at \upshape \ttfamily \href{https://www.#1}{#1}}

\DeclareFieldFormat{eprint:online}{%
	Manuscript available \href{https://#1}{online}}

\DeclareFieldFormat{eprint:inprep}{%
	\emph{In preparation}}

\DeclareFieldFormat{eprint:manual}{%
	\emph{#1}}

\DeclareFieldFormat{eprint:onarxiv}{%
	\emph{#1}}

\DeclareFieldFormat{eprint:toappear}{%
	\mbox{\emph{To appear}}}

\DeclareFieldFormat{eprint:accepted}{%
	\emph{Accepted at {#1}; to appear}}

\DeclareSourcemap{
	\maps[datatype=bibtex]{
		\map{
			\step[fieldsource=mrnumber,		fieldtarget=eprint, final]
			\step[fieldset=eprinttype,		fieldvalue=mrnumber]
		}
		\map{
			\step[fieldsource=arxiv,		fieldtarget=eprint, final]
			\step[fieldset=eprinttype,		fieldvalue=arxiv]
		}
		\map{
			\step[fieldsource=jstor,		fieldtarget=eprint, final]
			\step[fieldset=eprinttype,		fieldvalue=jstor]
		}
		\map{
			\step[fieldsource=customeprint,	fieldtarget=eprint, final]
			\step[fieldset=eprinttype,		fieldvalue=customeprint]
		}
		\map{
			\step[fieldsource=online,		fieldtarget=eprint, final]
			\step[fieldset=eprinttype,		fieldvalue=online]
		}
		\map{
			\step[fieldsource=inprep,		fieldtarget=eprint, final]
			\step[fieldset=eprinttype,		fieldvalue=inprep]
		}
		\map{
			\step[fieldsource=manual,		fieldtarget=eprint, final]
			\step[fieldset=eprinttype,		fieldvalue=manual]
		}
		\map{
			\step[fieldsource=onarxiv,		fieldtarget=eprint, final]
			\step[fieldset=eprinttype,		fieldvalue=onarxiv]
		}
		\map{
			\step[fieldsource=toappear,		fieldtarget=eprint, final]
			\step[fieldset=eprinttype,		fieldvalue=toappear]
		}
		\map{
			\step[fieldsource=accepted,		fieldtarget=eprint, final]
			\step[fieldset=eprinttype,		fieldvalue=accepted]
		}
	}
}

\DeclareFieldFormat{doi}{%
	\href{https://doi.org/#1}{DOI}%
}

\DeclareFieldFormat{url}{%
	\href{https://#1}{\texttt{#1}}%
}

\DeclareFieldFormat{comments}{%
	\emph{#1}%
} 
 
\DeclareFieldFormat
	[article,inbook,incollection,inproceedings,patent,thesis,unpublished]
	{title}{{#1\isdot}}
\DeclareFieldFormat
	[article,book,inbook,incollection,inproceedings,patent,thesis,unpublished, online]
	{date}{\mkbibparens{#1}}
\DeclareFieldFormat
	[article,book,inbook,incollection,inproceedings,patent,thesis,unpublished]
	{volume}{\mkbibbold{#1}}
\DeclareFieldFormat
	{pages}{\mkbibparens{#1}}

\DeclareFieldFormat{edition}{%
	\ifinteger{#1}
	{\mkbibordedition{#1}~\bibstring{edition}\addcomma}
	{#1~\bibstring{edition}\addcomma}}

\renewbibmacro*{volume+number+eid}{%
	\printfield{volume}%
\setunit*{\adddot}%
	\printfield{number}%
\setunit{\space}%
	\printfield{eid}%
}

\DeclareBibliographyDriver{article}{%
	\usebibmacro{bibindex}%
	\usebibmacro{begentry}%
	\usebibmacro{author}%
\setunit{\addspace}%
	\usebibmacro{date}%
\newunit\newblock
	\usebibmacro{title}%
\newunit\newblock
	\printfield{journaltitle}\addperiod%
\setunit*{\addspace}%
	\usebibmacro{volume+number+eid}%
\setunit*{\addspace}\newblock
	\printfield{pages}
\setunit*{\addspace}%\newblock
	\printfield{comments}%
	\usebibmacro{eprint}%
\setunit*{\addspace}%
	\printfield{doi}%
	\usebibmacro{finentry}%
}

\DeclareBibliographyDriver{online}{%
	\usebibmacro{bibindex}%
	\usebibmacro{begentry}%
	\usebibmacro{author}%
\setunit{\addspace}%
	\usebibmacro{date}%
\newunit\newblock
	\usebibmacro{title}%
\newunit\newblock
	\usebibmacro{eprint}%
	\usebibmacro{finentry}%
}

\DeclareBibliographyDriver{book}{%
	\usebibmacro{bibindex}%
	\usebibmacro{begentry}%
	\usebibmacro{author}%
\setunit{\addspace}\newblock
	\usebibmacro{date}%
\newunit\newblock
	\usebibmacro{title}%
\setunit*{\addspace}\newblock
	\printfield{edition}%
\newunit\newblock
	\printlist{publisher}%\addcomma
\setunit*{\addspace}%
	\printfield{volume}%
\setunit*{\addspace}\newblock%
	\printfield{comments}%
	\usebibmacro{eprint}%
\setunit*{\addspace}
	\printfield{doi}
	\usebibmacro{finentry}%
}

\DeclareBibliographyDriver{thesis}{%
	\usebibmacro{bibindex}%
	\usebibmacro{begentry}%
	\usebibmacro{author}%
\setunit{\addspace}\newblock
	\usebibmacro{date}%
\newunit\newblock
	\usebibmacro{title}.%
\setunit*{\addspace}\newblock
	\space Thesis, \printlist{institution}
	\usebibmacro{eprint}%
	\printfield{comments}%
\setunit*{\addspace}
	\printfield{doi}%
	\usebibmacro{finentry}%
}

\DeclareBibliographyDriver{inproceedings}{%
	\usebibmacro{bibindex}%
	\usebibmacro{begentry}%
	\usebibmacro{author}%
\setunit{\addspace}%
	\usebibmacro{date}%
\newunit\newblock
	\usebibmacro{title}%
\newunit\newblock
	\printfield{booktitle}\addcomma%
\setunit*{\addspace}%
	\printfield{series}\addcomma%
\setunit*{\addspace}\newblock
	\usebibmacro{editor}\addcomma
\setunit*{\addspace}\newblock
	\printlist{publisher}
\setunit*{\addspace}%
	\printfield{volume}%
\setunit*{\addspace}%
	\printfield{pages}
\setunit*{\addspace}%\newblock
	\usebibmacro{eprint}%
	\printfield{comments}%
\setunit*{\addnnspace}
	\printfield{doi}
	\usebibmacro{finentry}%
}

\DeclareBibliographyDriver{inbook}{%
	\usebibmacro{bibindex}%
	\usebibmacro{begentry}%
	\usebibmacro{author}%
\setunit{\addspace}%
	\usebibmacro{date}%
\newunit\newblock
	\usebibmacro{title}%
\newunit\newblock
	\printfield{booktitle}\addcomma%
\setunit*{\addspace}%
	\printfield{series}\addcomma%
\setunit*{\addspace}\newblock
	\usebibmacro{editor}\addcomma
\setunit*{\addspace}\newblock
	\printlist{publisher}%\addcomma%
\setunit*{\addspace}%
	\printfield{volume}%
\setunit*{\addspace}%
	\printfield{pages}
\setunit*{\addspace}%\newblock
	\usebibmacro{eprint}%
	\printfield{comments}%
\setunit*{\addspace}
	\printfield{doi}
	\usebibmacro{finentry}%
}

\DeclareBibliographyDriver{incollection}{%
	\usebibmacro{bibindex}%
	\usebibmacro{begentry}%
	\usebibmacro{author}%
\setunit{\addspace}%
	\usebibmacro{date}%
\newunit\newblock
	\usebibmacro{title}%
\newunit\newblock
	\printfield{booktitle}\addcomma%
\setunit*{\addspace}%
	\printfield{series}\addcomma%
\setunit*{\addspace}\newblock
	\printlist{publisher}%\addcomma%
\setunit*{\addspace}%
	\printfield{volume}%
\setunit*{\addspace}%
	\printfield{pages}
\setunit*{\addspace}%\newblock
	\usebibmacro{eprint}%
	\printfield{comments}%
\setunit*{\addspace}
	\printfield{doi}
	\usebibmacro{finentry}%
}

\AtEveryBibitem{%
	\clearfield{day}%
	\clearfield{month}%
} 
 
\DeclareNameAlias{sortname}{given-family}
\DeclareNameAlias{default}{given-family}

\usepackage{geometry}
\usepackage[hidelinks]{hyperref}
\hypersetup{
	colorlinks,
	linkcolor={blue},
	citecolor={ForestGreen},
	urlcolor={blue}
}

\usepackage{graphicx}
\graphicspath{{FMMC_Figures/}}
\usepackage[labelfont = {bf, sf}, labelsep = period, singlelinecheck = false]{caption}
\usepackage{titlesec}
\titleformat{\section}
	{\sffamily \Large \bfseries \boldmath}{\thesection}{1em}{}
\titleformat{\subsection}
	{\sffamily \large \bfseries \boldmath}{\thesubsection}{1em}{}
\titleformat{\subsubsection}
	{\sffamily \normalsize \bfseries \boldmath}{\thesubsubsection}{1em}{}
\titleformat{\paragraph}[runin]
	{\normalsize \itshape}{\theparagraph}{1em}{}
\titlespacing*{\paragraph}
	{0pt}{\medskipamount}{*2.5}
\titleformat{\subparagraph}[runin]
	{\normalsize \itshape}{\thesubparagraph}{1em}{}
\usepackage[utf8]{inputenc}

\usepackage[UKenglish]{babel}
\usepackage{amsmath, amsthm, amssymb, mathrsfs, bm, mathtools}
\usepackage{wasysym}
\allowdisplaybreaks

\usepackage{titlesec}
\usepackage{xifthen}
\usepackage{xspace}

\usepackage[normalem]{ulem}

\usepackage{footnotebackref}
\usepackage{enumitem}

\setlist[itemize,enumerate]{%
	itemsep	= 0pt,%
	topsep	= \smallskipamount%
}

\newcommand{\romannumbering}{%
	\renewcommand{\labelenumi}{\upshape(\roman{enumi})}%
	\renewcommand{\theenumi}{\upshape(\roman{enumi})}}

\usepackage{scrextend}

\setlist[description]{%
	topsep = \smallskipamount,		% space before start / after end of list
	itemsep = \smallskipamount,		% space between items
	font = {\mdseries\itshape},		% set the label font
	leftmargin = \parindent
}

\newlist{ranlist}{enumerate}{3}
\setlist[ranlist,1]{
	label=(\roman*) }
\setlist[ranlist,2]{
	label=(\textit{\alph*}),
	ref=(\roman{ranlisti}.\textit{\alph*})	}
\setlist[ranlist,3]{
	label=\arabic*.,
	ref=(\roman{ranlisti}.\textit{\alph{ranlistii}}.\arabic*)	}
\usepackage[dvipsnames,hyperref]{xcolor}

\newcommand{\green} [1]{{\color{ForestGreen} #1}}

\def\IfAmpersandUseAlign#1#2&#3\EndIfAmpersandUseAlign%
{%
	\if\relax\detokenize{#3}\relax
	\begin{equation*}%
	#1%
	\end{equation*}%
	\else
	\begin{align*}%
	#1%
	\end{align*}%
	\fi
}
\def\[#1\]%
{%
	\IfAmpersandUseAlign{#1}#1&\EndIfAmpersandUseAlign
}

\usepackage{eqparbox}

\newcommand{\tv}[1]{
	\mathchoice
	{\bigl\| #1 \bigr\|_{\TV}}%
	{\| #1 \|_{\TV}}%
	{\| #1 \|_{\TV}}%
	{\| #1 \|_{\TV}}
}

\newcommand{\one}  [1]{\bm1\{ #1 \}}

\let\originalexp\exp
\let\exp\relax
\DeclareRobustCommand{\exp} [1]{\originalexp(#1)}

\newcommand{\expb}[2][]{
	\ifthenelse{\equal{}{#1}}
	{\originalexp\bigl( #2 \bigr)}
	{\originalexp_{#1}\bigl( #2 \bigr)}
}

\newcommand{\logn}[1][]{
	\ifthenelse{\equal{}{#1}}
	{\log n}
	{(\log n)^{#1}}
}
\newcommand{\logk}[1][]{
	\ifthenelse{\equal{}{#1}}
	{\log k}
	{(\log k)^{#1}}
}

\newcommand{\Quad}[1]{
	\mathchoice
	{\quad\text{#1}\quad}%	\displaystyle
	{\text{ #1 }}%			\textsyle
	{\text{ #1 }}%			\scriptstyle
	{\text{ #1 }}%			\scriptscriptstyle
}

\newcommand{\Qand}{\Quad{and}}
\newcommand{\Qfor}{\Quad{for}}
\newcommand{\Qforall}{\Quad{for all}}

\newcommand{\Qwhere}{\Quad{where}}
\newcommand{\Qwith}{\Quad{with}}

\newcommand{\fnrestrict}[2]{\left. {#1} \right|_{#2}\!}

\newcommand{\abs}  [1]{| #1 |}

\newcommand{\rbr} [1]{ ( #1 ) }
\newcommand{\rbb} [1]{\bigl( #1 \bigr)}

\newcommand{\bra} [1]{ \{ #1 \} }
\newcommand{\brb} [1]{\bigl\{ #1 \bigr\}}

\newcommand{\brbb}[1]{\biggl\{ #1 \biggr\}}

\newcommand{\floor}  [1]{\lfloor #1 \rfloor}

\newcommand{\midb}{\bigm|}

\newcommand{\pr}[2][]{
	\mathchoice
	{\ifthenelse{\isempty{#1}}
		{\mathbb{P}\bigl(#2\bigr)}
		{\mathbb{P}_{#1}\bigl(#2\bigr)}}%	\displaystyle
	{\ifthenelse{\isempty{#1}}
		{\mathbb{P}(#2)}
		{\mathbb{P}_{#1}(#2)}}%	\textsyle
	{\ifthenelse{\isempty{#1}}
		{\mathbb{P}(#2)}
		{\mathbb{P}_{#1}(#2)}}%	\scriptstyle
	{\ifthenelse{\isempty{#1}}
		{\mathbb{P}(#2)}
		{\mathbb{P}_{#1}(#2)}}%	\scriptscriptstyle
}
\newcommand{\prt}[2][]{
	\ifthenelse{\equal{}{#1}}
	{\mathbb{P}( #2 )}
	{\mathbb{P}_{#1}( #2 )}
}
\newcommand{\prb}[2][]{
	\ifthenelse{\equal{}{#1}}
	{\mathbb{P}\bigl( #2 \bigr)}
	{\mathbb{P}_{#1}\bigl( #2 \bigr)}
}
\newcommand{\prB}[2][]{
	\ifthenelse{\equal{}{#1}}
	{\mathbb{P}\Bigl( #2 \Bigr)}
	{\mathbb{P}_{#1}\Bigl( #2 \Bigr)}
}
\newcommand{\prbb}[2][]{
	\ifthenelse{\equal{}{#1}}
	{\mathbb{P}\biggl( #2 \biggr)}
	{\mathbb{P}_{#1}\biggl( #2 \biggr)}
}
\newcommand{\prBB}[2][]{
	\ifthenelse{\equal{}{#1}}
	{\mathbb{P}\Biggl( #2 \Biggr)}
	{\mathbb{P}_{#1}\Biggl( #2 \Biggr)}
}
\newcommand{\prs}[2][]{
	\ifthenelse{\equal{}{#1}}
	{\mathbb{P}\left( #2 \right)}
	{\mathbb{P}_{#1}\left( #2 \right)}
}

\newcommand{\ex}[2][]{
	\mathchoice
	{\ifthenelse{\isempty{#1}}
		{\mathbb{E}\bigl(#2\bigr)}
		{\mathbb{E}_{#1}\bigl(#2\bigr)}}%	\displaystyle
	{\ifthenelse{\isempty{#1}}
		{\mathbb{E}(#2)}
		{\mathbb{E}_{#1}(#2)}}%	\textsyle
	{\ifthenelse{\isempty{#1}}
		{\mathbb{E}(#2)}
		{\mathbb{E}_{#1}(#2)}}%	\scriptstyle
	{\ifthenelse{\isempty{#1}}
		{\mathbb{E}(#2)}
		{\mathbb{E}_{#1}(#2)}}%	\scriptscriptstyle
}
\newcommand{\ext}[2][]{
	\ifthenelse{\equal{}{#1}}
	{\mathbb{E}( #2 )}
	{\mathbb{E}_{#1}( #2 )}
}
\newcommand{\exb}[2][]{
	\ifthenelse{\equal{}{#1}}
	{\mathbb{E}\bigl( #2 \bigr)}
	{\mathbb{E}_{#1}\bigr( #2 \bigr)}
}
\newcommand{\exB}[2][]{
	\ifthenelse{\equal{}{#1}}
	{\mathbb{E}\Bigl( #2 \Bigr)}
	{\mathbb{E}_{#1}\Bigl( #2 \Bigr)}
}
\newcommand{\exbb}[2][]{
	\ifthenelse{\equal{}{#1}}
	{\mathbb{E}\biggl( #2 \biggr)}
	{\mathbb{E}_{#1}\biggl( #2 \biggr)}
}
\newcommand{\exBB}[2][]{
	\ifthenelse{\equal{}{#1}}
	{\mathbb{E}\Biggl( #2 \Biggr)}
	{\mathbb{E}_{#1}\Biggl( #2 \Biggr)}
}

\newcommand{\Varb}[2][]{
	\ifthenelse{\equal{}{#1}}
	{\mathbb{V}\textnormal{ar} \bigl(#2\bigr)}
	{\mathbb{V}\textnormal{ar}_{#1} \bigl(#2\bigr)}
}
\newcommand{\VAR}[2][]{
	\ifthenelse{\equal{}{#1}}
	{\textnormal{Var}(#2)}
	{\textnormal{Var}_{#1}(#2)}
}
\usepackage{calc}
\newlength{\halfplusheight}
\newcommand{\maxt}[1]{
	\mathchoice
	{\textstyle \max_{#1} \displaystyle}
	{\max_{#1}}
	{\max_{#1}}
	{\max_{#1}}
}

\newcommand{\MAX}[1]{%
	\mathop{\raisebox{\halfplusheight}{\(\displaystyle\max_{#1}\)}}\:%
}

\newcommand{\mint}[1]{
	\mathchoice
	{\textstyle \min_{#1} \displaystyle}
	{\min_{#1}}
	{\min_{#1}}
	{\min_{#1}}
}

\newcommand{\MIN}[1]{%
	\mathop{\raisebox{\halfplusheight}{\(\displaystyle\min_{#1}\)}}\:%
}

\newcommand{\sumt}[2][]{
	\mathchoice
	{\ifthenelse{\isempty{#1}}
		{\textstyle \sum_{#2}      \displaystyle}
		{\textstyle \sum_{#2}^{#1} \displaystyle}}%	\displaystyle
	{\ifthenelse{\isempty{#1}}
		{\sum_{#2}}
		{\sum_{#2}^{#1}}}%			\textsyle
	{\ifthenelse{\isempty{#1}}
		{\sum_{#2}}
		{\sum_{#2}^{#1}}}%			\scriptstyle
	{\ifthenelse{\isempty{#1}}
		{\sum_{#2}}
		{\sum_{#2}^{#1}}}%			\scriptscriptstyle
}
\newcommand{\sumd}[2][]{
	\ifthenelse{\isempty{#1}}
		{\sum_{#2}}
		{\sum_{#2}^{#1}}
}

\newcommand{\intt}[2][]{
	\mathchoice
	{\ifthenelse{\isempty{#1}}
		{\textstyle \int_{#2}      \displaystyle}
		{\textstyle \int_{#2}^{#1} \displaystyle}}
	{\ifthenelse{\isempty{#1}}
		{\int_{#2}}
		{\int_{#2}^{#1}}}
	{\ifthenelse{\isempty{#1}}
		{\int_{#2}}
		{\int_{#2}^{#1}}}
	{\ifthenelse{\isempty{#1}}
		{\int_{#2}}
		{\int_{#2}^{#1}}}
}

\newcommand{\prodt}[2][]{
	\mathchoice
	{\ifthenelse{\isempty{#1}}
		{\textstyle \prod_{#2}      \displaystyle}
		{\textstyle \prod_{#2}^{#1} \displaystyle}}
	{\ifthenelse{\isempty{#1}}
		{\prod_{#2}}
		{\prod_{#2}^{#1}}}
	{\ifthenelse{\isempty{#1}}
		{\prod_{#2}}
		{\prod_{#2}^{#1}}}
	{\ifthenelse{\isempty{#1}}
		{\prod_{#2}}
		{\prod_{#2}^{#1}}}	
}
\newcommand{\prodd}[2][]{
	\ifthenelse{\isempty{#1}}
		{\prod_{#2}}
		{\prod_{#2}^{#1}}
}
\newcommand{\asinf}[1]{\text{as \(#1\to\infty\)}\xspace}

\newcommand{\dist}{\textnormal{dist}}

\newcommand{\TV}{%
	\ifmmode%
		\mathsf{TV}%
	\else%
		\textsf{TV}\xspace%
	\fi%
}

\DeclareMathOperator{\diam}{diam}

\newcommand{\st}{{\textnormal{ st }}}

\newcommand{\dmax}{d_\mathsf{max}}
\newcommand{\dmin}{d_\mathsf{min}}
\newcommand{\pimin}{\pi_\mathsf{min}}

\newcommand{\mbn}{\mathbb{N}}

\newcommand{\mbr}{\mathbb{R}}

\newcommand{\mbz}{\mathbb{Z}}

\newcommand{\cq}{\coloneqq}

\newenvironment{Proof}[1][\proofname]{%
	\proof[\upshape\bfseries\sffamily\boldmath#1]
}{\endproof}

\newcommand{\eps}{\varepsilon}

\newcommand{\blfootnote}[1]{%
	\begingroup
	\renewcommand\thefootnote{}\footnote{\sffamily#1}%
	\addtocounter{footnote}{-1}%
	\endgroup
}

\makeatletter
\let\save@mathaccent\mathaccent
\newcommand*\if@single[3]{%
	\setbox0\hbox{${\mathaccent"0362{#1}}^H$}%
	\setbox2\hbox{${\mathaccent"0362{\kern0pt#1}}^H$}%
	\ifdim\ht0=\ht2 #3\else #2\fi
}
\newcommand*\rel@kern[1]{\kern#1\dimexpr\macc@kerna}
\newcommand*\widebar[1]{\@ifnextchar^{{\wide@bar{#1}{0}}}{\wide@bar{#1}{1}}}
\newcommand*\wide@bar[2]{\if@single{#1}{\wide@bar@{#1}{#2}{1}}{\wide@bar@{#1}{#2}{2}}}
\newcommand*\wide@bar@[3]{%
	\begingroup
	\def\mathaccent##1##2{%
		\let\mathaccent\save@mathaccent
		\if#32 \let\macc@nucleus\first@char \fi
		\setbox\z@\hbox{$\macc@style{\macc@nucleus}_{}$}%
		\setbox\tw@\hbox{$\macc@style{\macc@nucleus}{}_{}$}%
		\dimen@\wd\tw@
		\advance\dimen@-\wd\z@
		\divide\dimen@ 3
		\@tempdima\wd\tw@
		\advance\@tempdima-\scriptspace
		\divide\@tempdima 10
		\advance\dimen@-\@tempdima
		\ifdim\dimen@>\z@ \dimen@0pt\fi
		\rel@kern{0.6}\kern-\dimen@
		\if#31
		\overline{\rel@kern{-0.6}\kern\dimen@\macc@nucleus\rel@kern{0.4}\kern\dimen@}%
		\advance\dimen@0.4\dimexpr\macc@kerna
		\let\final@kern#2%
		\ifdim\dimen@<\z@ \let\final@kern1\fi
		\if\final@kern1 \kern-\dimen@\fi
		\else
		\overline{\rel@kern{-0.6}\kern\dimen@#1}%
		\fi
	}%
	\macc@depth\@ne
	\let\math@bgroup\@empty \let\math@egroup\macc@set@skewchar
	\mathsurround\z@ \frozen@everymath{\mathgroup\macc@group\relax}%
	\macc@set@skewchar\relax
	\let\mathaccentV\macc@nested@a
	\if#31
	\macc@nested@a\relax111{#1}%
	\else
	\def\gobble@till@marker##1\endmarker{}%
	\futurelet\first@char\gobble@till@marker#1\endmarker
	\ifcat\noexpand\first@char A\else
	\def\first@char{}%
	\fi
	\macc@nested@a\relax111{\first@char}%
	\fi
	\endgroup
}
\makeatother
\usepackage{xparse}
\ExplSyntaxOn
\NewDocumentCommand{\mref}{m}{\quinn_mref:n {#1}}
\seq_new:N \l_quinn_mref_seq
\cs_new:Npn \quinn_mref:n #1
{
	\seq_set_split:Nnn \l_quinn_mref_seq { , } { #1 }
	\seq_pop_right:NN \l_quinn_mref_seq \l_tmpa_tl
	( % print the left bracket
	\seq_map_inline:Nn \l_quinn_mref_seq
	{ \ref{##1},\nobreakspace } % print the first references
	\exp_args:NV \ref \l_tmpa_tl % print the last or only one
	) % print the right bracket
}
\ExplSyntaxOff

\numberwithin{equation}{subsection}

\makeatletter
\newenvironment{subtheorem}[1]{%
	\def\subtheoremcounter{#1}%
	\refstepcounter{#1}%
	\protected@edef\theparentnumber{\csname the#1\endcsname}%
	\setcounter{parentnumber}{\value{#1}}%
	\setcounter{#1}{0}%
	\expandafter\def\csname the#1\endcsname{\theparentnumber\alph{#1}}%
	\expandafter\def\csname theH#1\endcsname{thm.\theparentnumber\alph{#1}}%
	\unskip\ignorespaces
}{%
	\setcounter{\subtheoremcounter}{\value{parentnumber}}%
	\ignorespacesafterend
}
\makeatother
\newcounter{parentnumber}

\makeatletter
\newenvironment{subtheorem-num}[1]{%
	\def\subtheoremcounter{#1}%
	\refstepcounter{#1}%
	\protected@edef\theparentnumber{\csname the#1\endcsname}%
	\setcounter{parentnumber}{\value{#1}}%
	\setcounter{#1}{0}%
	\expandafter\def\csname the#1\endcsname{\theparentnumber.\arabic{#1}}%
	\expandafter\def\csname theH#1\endcsname{thm.\theparentnumber.\arabic{#1}}%
	\unskip\ignorespaces
}{%
	\setcounter{\subtheoremcounter}{\value{parentnumber}}%
	\ignorespacesafterend
}
\makeatother
\newcommand{\qedtriangle}{\renewcommand{\qedsymbol}{\ensuremath{\triangle}}}

\usepackage[capitalise]{cleveref}

\crefname{figure}{Figure}{Figures}

\newcommand{\nextresult}{%
	\setcounter{introthm}{\value{introthm}}
	\setcounter{introcor}{\value{introthm}}
	\setcounter{introconj}{\value{introthm}}
	\setcounter{introdefn}{\value{introthm}}
	\setcounter{intrormkT}{\value{introthm}}
}

\newtheoremstyle{sfsl}
{1\baselineskip}		% Space above
{1\baselineskip}		% Space below
{\slshape}				% Theorem body font
{}						% Indent amount
{\bfseries\sffamily}	% Theorem head font
{.}						% Punctuation after theorem head
{0.5em}					% Space after theorem head
{\thmname{#1}\thmnumber{ #2}\thmnote{ {\mdseries(#3)}}}
\newtheoremstyle{sfup}
{1\baselineskip}		% Space above
{1\baselineskip}		% Space below
{\upshape}				% Theorem body font
{}						% Indent amount
{\bfseries\sffamily}	% Theorem head font
{.}						% Punctuation after theorem head
{0.5em}					% Space after theorem head
{\thmname{#1}\thmnumber{ #2}\thmnote{ {\mdseries(#3)}}}
\theoremstyle{sfsl}

\newtheorem*{thm*}{Theorem}
\newtheorem{thm} {Theorem}[section]
\crefname{thm}{Theorem}{Theorems}

\newtheorem*{introthm*}{Theorem}
\newtheorem{introthm}{Theorem}

\crefname{introthm}{Theorem}{Theorems}

\newtheorem*{cor*}{Corollary}
\newtheorem{cor} [thm]{Corollary}
\crefname{cor}{Corollary}{Corollaries}

\newtheorem*{introcor*}{Corollary}

\crefname{introcor}{Corollary}{Corollaries}

\newtheorem*{introconj*}{Conjecture}

\crefname{introconj}{Conjecture}{Conjectures}

\newtheorem*{introques*}{Question}

\crefname{introques}{Question}{Questions}

\newtheorem*{lem*}    {Lemma}
\newtheorem{lem} [thm]{Lemma}
\crefname{lem}{Lemma}{Lemmas}

\newtheorem*{introlem*}{Lemma}

\crefname{introlem}{Lemma}{Lemmas}

\newtheorem*{prop*}    {Proposition}
\newtheorem{prop} [thm]{Proposition}
\crefname{prop}{Proposition}{Propositions}

\newtheorem*{clm*}    {Claim}

\crefname{clm}{Claim}{Claims}

\newtheorem*{defn*}    {Definition}
\newtheorem{defn} [thm]{Definition}
\crefname{defn}{Definition}{Definitions}

\newtheorem*{introdefn*}{Definition}
\newtheorem{introdefn}{Definition}

\crefname{introdefn}{Definition}{Definitions}

\newtheorem*{nota*}{Notation}

\newtheorem{innercustomgenericsl}{\customgenericnamesl}
\providecommand{\customgenericnamesl}{}
\newcommand{\newcustomtheoremsl}[2]{%
	\newenvironment{#1}[1]
	{%
		\renewcommand\customgenericnamesl{#2}%
		\renewcommand\theinnercustomgenericsl{##1}%
		\innercustomgenericsl
	}
	{\endinnercustomgenericsl}
}

\newcustomtheoremsl{customthm}{Theorem}
\crefname{customthm}{Theorem}{Theorems}
\newcustomtheoremsl{customprop}{Proposition}
\newcustomtheoremsl{customlemma}{Lemma}
\newcustomtheoremsl{customrmk}{Remark}
\newcustomtheoremsl{customconj}{Conjecture}
\newcustomtheoremsl{customdefn}{Definition}
\newcustomtheoremsl{customquestion}{Question}
\newcustomtheoremsl{customopenproblem}{Open Problem}
\newcustomtheoremsl{customhyp}{Hypothesis}

\newtheorem*{conj*}   {Conjecture}

\crefname{conj}{Conjecture}{Conjectures}

\newtheorem*{rmk*}{Remark}

\theoremstyle{sfup}

\newtheorem{innercustomgenericup}{\customgenericnameup}
\providecommand{\customgenericnameup}{}
\newcommand{\newcustomtheoremup}[2]{%
	\newenvironment{#1}[1]
	{%
		\renewcommand\customgenericnameup{#2}%
		\renewcommand\theinnercustomgenericup{##1}%
		\innercustomgenericup
	}
	{\endinnercustomgenericup}
}

\crefname{exm} {Example}{Examples}
\crefname{exmT}{Example}{Examples}

\newtheorem{rmkT}[thm]{Remark}
	\newenvironment{rmkt}
	{\pushQED{\qed}\renewcommand{\qedsymbol}{\ensuremath{\triangle}}\rmkT}
	{\popQED\endrmkT}
\crefname{rmk} {Remark}{Remarks}
\crefname{rmkT}{Remark}{Remarks}

\newtheorem*{rmkTT}{Remark}
\newenvironment{rmkt*}
	{\pushQED{\qed}\renewcommand{\qedsymbol}{\ensuremath{\triangle}}\rmkTT}
	{\popQED\endrmkTT}

\newcustomtheoremup{customrmkT}{Remark}

\crefname{rmks} {Remarks}{Remarks}
\crefname{rmksT}{Remarks}{Remarks}

\newtheorem*{rmks*} {Remarks}
\newtheorem*{rmksTT}{Remarks}
\newenvironment{rmkst*}
	{\pushQED{\qed}\renewcommand{\qedsymbol}{\ensuremath{\triangle}}\rmksTT}
	{\popQED\endrmksTT}

\crefname{intrormk} {Remark}{Remarks}
\crefname{intrormkT}{Remark}{Remarks}

\newtheorem*{intrormk*} {Remark}
\newtheorem*{intrormkTT}{Remark}
\newenvironment{intrormkt*}
	{\pushQED{\qed}\renewcommand{\qedsymbol}{\ensuremath{\triangle}}\intrormkTT}
	{\popQED\endintrormkTT}

\newtheorem*{exm*} {Example}
\newtheorem*{exmTT}{Lemma}
	\newenvironment{exmt*}
	{\pushQED{\qed}\renewcommand{\qedsymbol}{\ensuremath{\triangle}}\exmTT}
	{\popQED\endexmTT}

\newtheorem*{note*} {Note}
\newtheorem*{noteTT}{Note}
	\newenvironment{notet*}
	{\pushQED{\qed}\renewcommand{\qedsymbol}{\ensuremath{\triangle}}\noteTT}
	{\popQED\endnoteTT}

\crefformat{section}{\S#2#1#3}
\crefformat{subsection}{\S#2#1#3}
\crefformat{subsubsection}{\S#2#1#3}

\crefformat{equation}{(#2#1#3)}
\crefmultiformat{equation}{(#2#1#3}{, #2#1#3)}{, #2#1#3}{, #2#1#3)}
\crefrangeformat{equation}{(#3#1#4--#5#2#6)}

\usepackage{pgfplots,tikz-3dplot}
\usepackage{tikz,pgf}
\usepackage{tikz-3dplot} 

\usetikzlibrary{shapes.callouts,decorations.pathmorphing,arrows,snakes,shapes,graphs,graphs.standard,calc,positioning}

\usepackage{multicol,multirow}
\usepackage{xcolor}

\newenvironment{center-small}
	{\par\centering\smallskip}
	{\par\smallskip}

\usepackage{chngcntr}

\newcommand{\FMMC}{\textsf{FMMC}\xspace}

\newcommand{\RW}{\textsf{RW}\xspace}
\newcommand{\RWs}{\textsf{RW}s\xspace}
\newcommand{\SDP}{\textsf{SDP}\xspace}

\newcommand{\U}{\textsf{U}\xspace}

\newcommand{\unif}[1][]{%
	\mathcal U%
	\ifthenelse{\equal{#1}{}}{}{_{#1}}%
}

\newcommand{\wrt}{w.r.t.\xspace}

\newcommand{\LL}{\ensuremath{L}}
\newcommand{\PP}{\ensuremath{P}}
\newcommand{\QQ}{\ensuremath{Q}}
\newcommand{\qq}{\ensuremath{q}}

\newcommand{\dumbbell}{D_\star}
\newcommand{\binarytree}{\mathbb T}
\newcommand{\stargraph}{G_\star}
\newcommand{\matching}{\mathcal M}
\newcommand{\source}{\Sigma}

\newcommand{\optSG}[1][]{%	optimal spectral gap
	\gamma_\star%
		\ifthenelse{\equal{}{#1}}%
		{}%
		{\rbr{#1}}%
}

\NewDocumentCommand{\MIX}{oo}{%
	\tau%
	\IfValueT{#1}{_{#1}}%
	\IfValueT{#2}{(#2)}%
}

\newcommand{\MM}[1][]{%		set of matchings
	\ifthenelse{\equal{}{#1}}%
	{\mathcal M}%
	{\mathcal M(#1)}%
}

\newcommand{\MW}[2][]{%		matching weight
	\ifthenelse{\equal{}{#1}}%
	{\upsilon(#2)}%
	{\upsilon_{#1}(#2)}%
}
\newcommand{\MWb}[2][]{%	matching weight \big
	\ifthenelse{\equal{}{#1}}%
	{\upsilon\bigl(#2\bigr)}%
	{\upsilon_{#1}\bigl(#2\bigr)}%
}

\newcommand{\PD}[1]{%	set of probability measures/distributions
	\mathcal D(#1)%
}

\newcommand{\TM}[1]{%	set of transition matrices
	\mathcal M(#1)%
}

\NewDocumentCommand{\EC}{soo}{%	(optimal) edge conductance
	\Phi%
	\IfBooleanT{#1}{^\star}%	add star
	\IfValueT{#3}{_{#3}}%		add probability measure
	\IfValueT{#2}{(#2)}%		add graph/set
}

\NewDocumentCommand{\ECRW}{soo}{%	(optimal) edge conductance for RW
	\Phi%
	\IfBooleanT{#1}{^\star}%	add star
	\IfValueT{#2}{_{#2}}%		add probability measure
	\IfValueT{#3}{(#3)}%		add graph/set
}

\NewDocumentCommand{\aECRW}{soo}{%	(optimal) adjusted edge conductance for RW
	\widetilde \Phi%
	\IfBooleanT{#1}{^\star}%	add star
	\IfValueT{#2}{_{#2}}%		add probability measure
	\IfValueT{#3}{(#3)}%		add graph
}

\NewDocumentCommand{\MC}{soo}{%	(optimal) vertex conductance
	\Upsilon%
	\IfBooleanT{#1}{^\star}%	add star
	\IfValueT{#3}{_{#3}}%		add probability measure
	\IfValueT{#2}{(#2)}%		add graph
}

\makeatletter
\newcommand{\SG}{%
	\@ifstar%
	\SGopt%[#1][#2]%
	\SGgen%[#1]%
}
\NewDocumentCommand{\SGopt}{oo}{%
	\gamma^\star%
	\IfValueT{#2}{_{#2}}%		add probability measure
	\IfValueT{#1}{(#1)}%		add graph
}
\NewDocumentCommand{\SGgen}{o}{%
	\gamma%
	\IfValueT{#1}{_{#1}}%		add transition matrix
}
\makeatother

\makeatletter
\newcommand{\REL}{%
	\@ifstar%
	\RELopt%[#1][#2]%
	\RELgen%[#1]%
}
\NewDocumentCommand{\RELopt}{oo}{%
	\rho^\star%
	\IfValueT{#2}{_{#2}}%		add probability measure
	\IfValueT{#1}{(#1)}%		add graph
}
\NewDocumentCommand{\RELgen}{o}{%
	\rho%
	\IfValueT{#1}{_{#1}}%		add transition matrix
}
\makeatother

\newcommand{\bestgap}[1]{%	spectral gap of the fastest mixing chain
	\gamma^\star(#1)
}
\newcommand{\onedimgap}[1]{%	one-dimensional Roch's characterisation
	\gamma^{(1)}(#1)
}

\NewDocumentCommand{\VC}{soo}{%	(optimal) vertex conductance
	\Psi%
	\IfBooleanT{#1}{^\star}%	add star
	\IfValueT{#3}{_{#3}}%		add probability measure
	\IfValueT{#2}{(#2)}%		add graph
}

\DeclareMathOperator{\vol}{vol}		

\newcommand{\DIAM}%	diameter
	{\ensuremath{\Delta}}
\newcommand{\uu}%	weight
	{\ensuremath{u}}
\newcommand{\ww}%	weight
	{\ensuremath{w}}

\newcommand{\anc}{\mathsf{anc}}
\newcommand{\BFS}{\textsf{BFS}\xspace}
\DeclareMathOperator{\depth}{depth}
\newcommand{\prnt}{\mathsf{prt}}

\title{\sffamily%
	Geometric Bounds on the Fastest Mixing Markov Chain
}
\author{\sffamily%
	Sam Olesker-Taylor\qquad Luca Zanetti
}
\date{}

\begin{document}

\maketitle

\renewcommand{\abstractname}{\sffamily Abstract}

\begin{abstract}
In the Fastest Mixing Markov Chain problem, we are given a graph $G = (V, E)$ and desire the discrete-time Markov chain with smallest mixing time $\tau$ subject to having equilibrium distribution uniform on $V$ and non-zero transition probabilities only across edges of the graph.

It is well-known that the mixing time $\tau_\textsf{RW}$ of the lazy random walk on $G$ is characterised by the edge conductance $\Phi$ of $G$ via Cheeger's inequality: $\Phi^{-1} \lesssim \tau_\textsf{RW} \lesssim \Phi^{-2} \log |V|$. Analogously, we characterise the fastest mixing time $\tau^\star$ via a Cheeger-type inequality but for a different geometric quantity, namely the vertex conductance $\Psi$ of $G$: $\Psi^{-1} \lesssim \tau^\star \lesssim \Psi^{-2} (\log |V|)^2$.

This characterisation forbids fast mixing for graphs with small vertex conductance. To bypass this fundamental barrier, we consider Markov chains on $G$ with equilibrium distribution which need not be uniform, but rather only $\varepsilon$-close to uniform in total variation. We show that it is always possible to construct such a chain with mixing time $\tau \lesssim \varepsilon^{-1} (\operatorname{diam} G)^2 \log |V|$.

Finally, we discuss analogous questions for continuous-time and time-inhomogeneous chains.
\end{abstract}

\small
\begin{quote}
\begin{description}
	\item [Keywords:]
	mixing time,
	random walks,
	conductance,
	fastest mixing Markov chain
	
	\item [MSC 2020 subject classifications:]
	05C81;
	60J10,
	60J20,
	60J27
\end{description}
\end{quote}
\normalsize

%05C12: Distance in graphs
%05C48: Expander graphs
%05C80:	Random graphs
%05C81:	Random walks on graphs

%20C15: Ordinary representations and characters
%20C30: Representations of finite symmetric groups

%42A61: Probilistic methods in Fourier analysis

%43A30: Fourier and Fourier-Stieltjes transforms on non-Abelian groups and on semigroups, etc
%43A65: Representations of groups, semigroups, etc
%43A75: Analysis on specific compact groups
%43A90: Spherical functions (in Abstract harmonic analysis)

%60B15: Probability measures on groups or semigroups, Frouier transforms, factorization
%60C05: Combinatorial probability
%60G50: Sums of independent random variables; random walks
%60J10: Markov chains (discrete-time Markov processes on discrete state spaces)
%60J20: Applications of Markov chains and discrete-time Markov processes on general state spaces
%60J90: Coalescent processes
%60J27:	Continuous-time Markov processes on discrete state spaces
%60K35: Interacting random processes; statistical mechanics type models; percolation theory
%60K37:	Processes in random environments

\blfootnote{{}%
\\
	\eqmakebox[fm][l]{Sam Olesker-Taylor,}\quad%
	\href{mailto:oleskertaylor.sam@gmail.com}{oleskertaylor.sam@gmail.com}%
\\
	\eqmakebox[fm][l]{Luca Zanetti,}\quad%
	\href{mailto:lz2040@bath.ac.uk}{lz2040@bath.ac.uk}%
\smallskip
\\
	Department of Mathematical Sciences, University of Bath, UK%
\smallskip
\\
%\hspace*{1em}%
	SOT was supported by EPSRC grant EP/N004566/1
}

%\newpage
\vfill
\setcounter{tocdepth}{5}
\setcounter{tocdepth}{1}
\sffamily
\tableofcontents
\normalfont
\vspace*{\bigskipamount}

\romannumbering

%\section*{Some \LaTeX\ St

\newpage
\section{Introduction}
\label{sec:intro}

\subsection{Fastest Mixing Markov Chain Set-Up and Motivation}
\label{sec:intro:set-up}

Sampling objects from a finite set is a basic primitive which has a myriad of applications.
Sampling directly from such a set, however, may be computationally too expensive or even impossible, for example, if the objects are nodes of a distributed network.
A common approach in these scenarios is to design a random walk (\textit{\RW}), or, more generally, a Markov chain
	with state space corresponding to the set from which we wish to sample
and
	appropriate equilibrium distribution.
Furthermore, to ensure our sampling procedure is computationally efficient, 
	we desire our Markov chain to converge to equilibrium in a small number of steps,
	i.e., have fast mixing time.

This has wide-ranging applications:
	from shuffling cards
	\cite{BD:riffle-shuffle,DS:random-trans},
	to approximating statistical physics models
	\cite{JS:poly-ising,DHJM:poly-ising}
and
	analysing load-balancing protocols in distributed computing
\cite{RSW:load-balancing,TS:rw-load-balancing}.
Furthermore, approximately sampling from the uniform distribution of a set can be used to estimate the size of the set itself
\cite{JS:mcmc}.
This has been applied to
	approximating the permanent of a matrix
	\cite{JS:permanent,JSV:permanent:jour}
and
	counting the number of independent sets
	\cite{DJMV:beyond-permanent},
	 perfect matchings
	\cite{DM:pm-switch-chain}
and
	forests
	\cite{A:num-forests}
	in graphs.

Fundamental to  these applications is a fast mixing time.
Understanding in which instances fast mixing is achievable and what the intrinsic obstacles to fast mixing are is the focus of this paper.
More precisely, we consider the following scenario.
\begin{itemize}[topsep = \medskipamount]
	\item 
	We are given a finite, undirected graph $G = (V, E)$:
		the vertex set represents the underlying state space,
	while
		the edge set $E$ defines the transitions allowed.
	
	\item 
	Our goal is to study the \emph{fastest} mixing Markov chain satisfying these constraints.
\end{itemize}
We assume throughout that graphs are finite, undirected and connected.
%We do not repeat this.

This problem was originally introduced by \textcite{BDX:fmmc-graph} as the \textit{Fastest Mixing Markov Chain} (\textit{\FMMC}) problem.
Specifically, by considering only reversible chains and optimising the \emph{spectral gap} as a proxy for the mixing time, they recast the problem of finding the fastest mixing Markov chain on a graph as a convex optimisation problem.
Analogously to \textcite{BDX:fmmc-graph}, we dedicate most of our attention to reversible, time-homogenous chains in discrete-time.
We do, however, dedicate one section to questions in the continuous-time setting, first studied in \cite{SBXD:fmmc-cts}, and one short final section to time-inhomogeneous chains. Compared with discrete-time chains, continuous-time and time-inhomogeneous chains are considerably more powerful, but perhaps less natural from an application viewpoint.
%We do, however, dedicate one section to questions in the continuous-time setting, first studied in \cite{SBXD:fmmc-cts}.
%Compared with discrete-time, continuous-time affords considerable more power to Markov chains on graphs.
%We also dedicate the final short section to time-inhomogeneous chains, which are even more powerful, \green{albeit less natural from an application viewpoint}.

% optimisation problem was introduce by \textcite{BDX:fmmc-graph};
%they restricted to \textit{reversible} Markov chains and termed it the 
%\textit{Fastest Mixing Markov Chain} (\textit{\FMMC}) problem.
%The majority of our analysis is in the discrete-time set-up, with one part on continuous-time.

\smallskip

To be explicit and precise,
a transition matrix $\PP$ is \textit{reversible} with respect to (\textit{\wrt}) $\pi$ if
\(
	\pi(u) \PP(u,v) = \pi(v) \PP(v,u)
\)
for all
\(
	u,v \in V.
\)
Results on the spectral gap below---both our own and those referenced---are always in the reversible set-up.
%Reversible Markov chains can be viewed as random walks (\textit{\RW}s) on weighted graphs.
We also restrict to \textit{lazy} chains,
	ie chains with $\PP(v,v) \ge \tfrac12$ for all $v \in V$.
%A chain can be made lazy by only taking a step half the time:
	%replace $\PP$ with $\tfrac12 (I + \PP)$.
This is without loss of generality, since we are interested in maximising the spectral gap and this restriction costs a factor of at most $\tfrac12$ in the optimal~spectral~gap.
%Since we are interested in maximising the spectral gap, ; this restriction costs a factor at most $\tfrac12$.
%Our results are all regarding the \emph{order} of the optimum, so this constant factor is unimportant.
%Also, some tools apply only to lazy chains.

\smallskip

There are a variety of choices for how ``convergence to equilibrium'' is measured.
It is typically measured in the total variation (\textit{\TV}), or equivalently $\ell_1$, distance.
Other popular measures, particularly in the statistics literature, include $\ell_2$, or $\chi^2$, distance and relative entropy, or Kullback--Leibler divergence.
We also recall that $\ell_2$ is equivalent to $\ell_\infty$, or uniform, distance for reversible~chains.

Nevertheless, no matter which one of these measure we choose, the long-term convergence to equilibrium of a lazy, reversible Markov chain is governed by its spectral gap $\SG[\PP]$.
More precisely,
given a transition matrix $\PP$ of a lazy, reversible Markov chain,
	let $d_\PP(t, x)$ denote the distance between $\PP^t(x,\cdot)$ and its equilibrium distribution according to any of the aforementioned measures
and
	let $d_\PP(t) \cq \max_{x \in V} d_P(t, x)$.
Then,
\[
%	\LIM{\toinf t}
	d_\PP(t)^{1/t}
\to
	\SG[\PP]
\quad
	\asinf t.
\]
See \cite[Theorems~12.4 and~12.5]{LPW:markov-mixing} for details.
%This asymptotic relations holds for any of the aforementioned distances.
The spectral gap thus determines the asymptotic convergence to equilibrium without having to select a specific measure.

We now define formally the class of reversible Markov chains on a graph and then the spectral gap and the relaxation and mixing times.
%We now define precisely the spectral gap and the relaxation and mixing times.

\begin{introdefn*}[Markov Chains on a Graph]
Let $G = (V, E)$ be a graph and $\pi$ a probability measure on $V$.
A transition matrix $\PP \in [0,1]^{V \times V}$ is \textit{on $G$} if
	$\PP(u,v) > 0$ then either $\bra{u,v} \in E$ or $u = v$.
Let $\TM{G, \pi}$ denote the set of lazy transition matrices on $G$ which are reversible \wrt $\pi$.
\end{introdefn*}

\begin{introdefn*}[Spectral Gap, Relaxation Time and Mixing Time]
Let $G = (V, E)$ be a graph and $\pi$ a probability measure on $V$.
%Denote the set of lazy transition matrices on $G$ which are reversible \wrt $\pi$ by $\TM{G, \pi}$.
%
Let $\PP \in \TM{G, \pi}$.
%Let $\PP$ be the transition matrix of a lazy, reversible Markov chain.
%
The \textit{spectral gap} is $\SG[\PP] \cq 1 - \lambda_\PP$, where $\lambda_\PP$ is the largest non-unitary eigenvalue of $\PP$.
The \textit{relaxation time} is $1/\SG[\PP]$.
The \textit{(uniform) mixing time} is
%defined as
\(
	\MIX[\PP][\xi]
\cq
	\inf\brb{ t \ge 0 \midb d^\infty_\PP(t) \le \xi }
\Qfor
	\xi \in [0, 1],
\)
where $d^\infty$ is the $\ell_\infty$-distance.
Write $\MIX[\PP] \cq \MIX[\PP][\tfrac14]$.
%We use the standard abbreviation $\MIX[\PP] \cq \MIX[\PP][\tfrac14]$.
	%
\end{introdefn*}

There is a standard relation between the relaxation and mixing times:
\[
	\SG[\PP]^{-1}
\lesssim
	\MIX[\PP]
\lesssim
	\SG[\PP]^{-1} \cdot \log \pimin^{-1};
%\Quad{recalling that}
%	\MIX[\PP]
%=
%	\MIX[\PP][\tfrac14].
\]
see \cite[Theorems~12.4 and~12.5]{LPW:markov-mixing} for details.
The ``$\lesssim$'' symbol hides an implicit universal constant; we use the symbols ``$\gtrsim$'' and ``$\asymp$'' similarly.
Typically, $\log \pimin^{-1} \asymp \log \abs V$.
So, the relaxation time is a proxy for the mixing time, as well as characterising
long-term convergence to~equilibrium.

%\begin{introdefn*}[`Fast' Mixing]
%	%
%We say that $\PP$ is \textit{fast mixing} if $\SG[\PP]^{-1}$ is polylogarithmic in $\abs V$.
%	%
%\end{introdefn*}

%We use the following terminology.
%\begin{itemize}[noitemsep]
%	\item 
%	Mixing is \textit{fast} if
%		$\SG[\PP]^{-1}$ is polylogarithmic in $\abs V$.
%	
%	\item 
%	Mixing is \textit{slow} if
%		$\SG[\PP]^{-1}$ is polynomial in $\abs V$, i.e.\ grows like some positive power~of~$\abs V$.
%\end{itemize}

We are now finally ready to formally introduce the \FMMC problem.

%Let $G = (V, E)$ be a connected, undirected graph.
%Let $\PD{V}$ denote the set of positive probability distributions on $V$ and let $\pi \in V$
%and let $\pi$ be a probability distribution on $V$.
%We include all self-loops in $E$ when we speak of a graph $G = (V, E)$,
%	ie $\bra{v} = \bra{v,v} \in E$ for all $v \in V$,
%unless stated otherwise. \blue{I think it would be easier to not add self-loops directly, just add that $\PP(u,u)$ can be non-zero}
%Let $\TM{G, \pi}$ denote the set of all transition matrices $\PP$ on $G$ which are lazy and reversible \wrt $\pi$:
%\begin{itemize}[noitemsep]
%	\item 
%	\textit{on $G$} means
%	\(
%		\PP \in [0,1]^{V \times V},
%	\
%		\sumt{v \in V} \PP(u,v) = 1 \text{ for all } u \in V
%	\Qand
%		\PP(u,v) = 0 \text{ if } \bra{u,v} \notin E;
%	\)
%	
%	\item 
%	\textit{lazy} and \textit{reversible} means
%	\(
%		\mint{v} \PP(v,v) \ge \tfrac12
%	\Qand
%		\pi(u) \PP(u,v) = \pi(v) \PP(v,u) \text{ for all } u,v \in V.
%	\)
%\end{itemize}
%Write $\PD{V}$ for the set of positive probability distributions on the set $V$. \blue{Should we say why we ask for lazy?}
%\green{as I note above, by the time we explain \emph{why} we study lazy and reversible chains, we may as well define them properly at the same time. I'm not suggesting we copy+paste this format. Rather, we just say informally what ``on a graph'' means, then give the formal definition. ``lazy''/''reversible'' is probably easiest to just define formally}
%%, i.e.\ those giving non-zero weight to all elements of $V$.
%
%Given $G$ and $\pi$,
%the \FMMC problem aims to find $\PP \in \TM{G, \pi}$ with largest spectral gap.

\begin{defn*}[Fastest Mixing Markov Chain]
	Let $G = (V, E)$ be a graph and let $\pi$ be a probability measure on $V$.
	The \textit{optimal spectral gap} is defined as
	\[
		\SG*[G][\pi]
	\cq
		\max\brb{ \SG[\PP] \mid \PP \in \TM{G, \pi} }.
	\]
%	The space $\TM{G, \pi}$ is compact, thus it is legitimate to write $\max$ rather than $\sup$.
	The \textit{optimal relaxation time} is
	\(
		1 / \SG*[G][\pi].
	\)
	We write $\SG*[G]$, omitting the $\pi$, when $\pi = \unif[V]$ is uniform.
	
	A transition matrix $\PP$ is \textit{fast mixing} if $1/\SG[\PP]$ is polylogarithmic in $\abs V$, asymptotically.
	Analogously,
	a graph $G$ \textit{admits a fast mixing chain} if $1/\SG*[G]$ is polylogarithmic in $\abs V$,~asymptotically.
\end{defn*}

%The \FMMC problem was originally introduced by \textcite{BDX:fmmc-graph}, which was the first in a series of articles \cite{BDX:fmmc-graph,BDSX:fmmc-path,SBXD:fmmc-cts,BDPX:symm-anal,BDPX:fmmc-symm} by the those authors along with Parrilo and Sun.
%We give an overview of the history of the question in \cref{sec:intro:previous-work}. \blue{We already said this above. Consider cutting it or rewriting it}.

Previous work has been mainly focussed on
	finding useful formulations of the problem
or
	on solving particular cases;
see \cref{sec:intro:previous-work} for further details.
The primary aim of our work, instead, is~twofold:
\begin{enumerate}[noitemsep]
	\item 
%	first,
	to control the optimal spectral gap in terms of geometric barriers
	in the graph;
	
	\item 
%	second,
	to find ways to overcome these geometric barriers by slightly relaxing the \FMMC problem.
\end{enumerate}
Finally, we centre our attention on the case where $\pi = \unif[V]$ is the uniform distribution on $V$.
This case was also the main focus of the original series of papers studying the \FMMC problem.

\subsection{Main Results}
\label{sec:intro:main}

\begingroup

\titleformat{\subsubsection}
	{\sffamily \large \bfseries \boldmath}{\thesubsubsection}{1em}{}

\renewcommand{\thesubsubsection}{\Alph{subsubsection}}

This article includes multiple avenues of study, all on the theme of finding fast mixing Markov chains.
We introduce these and the main theorems that we prove in the following subsections.

\subsubsection{Characterisation of Fast Mixing on Graphs}
\label{sec:intro:main:vc}

\begin{subtheorem-num}{introdefn}
	\label{def:intro:main:vc}

We are looking for some natural statistic of the graph $G$ which \emph{characterises} fast mixing:
	we desire necessary and sufficient conditions for $1/\SG*[G]$ to be `small',
	namely polylogarithmic~in~$|V|$.
	
How well-connected a graph is should, intuitively, influence how fast a chain on the graph can mix.
Thus, we would like to understand what kind of connectivity measure best characterises fast mixing.
A natural candidate is the \textit{edge conductance} $\EC*$ of a graph, which is defined as follows.

\begin{introdefn}[Edge Conductance]
\label{def:intro:main:vc:edge}
	%
%Let $G = (V, E)$ be a graph.
The \textit{edge conductance} $\EC*[G]$ of a graph $G = (V, E)$ is
%defined by
\[
	\EC*[G]
\cq
	\MIN{S \subseteq V : 0 < \vol(S) \le \vol(V)/2}
	\EC[S]
\Qwhere
	\EC[S]
\cq
	%\ww(\partial_\mathrm{edge} S) / \ww(S),
	\abs{ E(S,S^c) }/ \vol(S)
\Qfor
	S \subseteq V,
\]
where
$E(S,S^c)$ is the \textit{edge boundary} of $S \subseteq V$ and $\vol(S)$ is the \textit{volume} of $S \subseteq V$:
\[
	E(S,S^c)
\cq
	\brb{ \bra{x,y} \in E \mid x \in S, \: y \notin S }
\Qand
	\vol(S)
\cq
	\sumt{x \in S}
	\deg(x).
\]
%Define the \textit{uniform edge conductance} $\EC*[G][\textup{\U}]$ similarly but with $\vol(S)$ replaced by $\abs S$.
%
%The \RW and uniform conductances satisfy
%\(
%	\dmin
%\le
%	\EC*[G][\textup{\U}] / \EC*[G][\textup{\RW}]
%\le
%	\dmax.
%\)
	%
\end{introdefn}

It is well-known the edge conductance $\EC*$ characterises the spectral gap of the lazy random walk (abbreviated \textit{\RW}) $P^\RW$ on $G$ via the discrete Cheeger inequality,
discovered in \cite{JS:approx-counting,LS:cheeger}:
\[
	\EC*[G]^2
\lesssim
	\SG[P^\RW]
\lesssim
	\EC*[G].
\]

The lazy \RW on a graph, however, does not have uniform equilibrium distribution, unless the graph is regular.
For this reason, it is natural to consider
the \textit{uniform}, or \textit{maximum degree}, \RW $P^\U$ which is defined by adding the appropriate number of self-loops to each vertex so that the becomes graph regular.
%This is equivalent to adding the appropriate vertex-dependent laziness so that the simple \RW on this modified graph has uniform equilibrium distribution. 
A simple calculation with the Dirichlet characterisation gives
\[
	\SG[P^\U]
\le
	\SG*[G]
\le
	2 \dmax \SG[P^\U];
\]
see \cite[\S 7.2]{BDX:fmmc-graph} for details.
Applying this along with the discrete Cheeger inequality gives
\[
	\EC*[G]^2 \cdot \dmin/\dmax
\lesssim
	\SG[P^\U]
\le
	\SG*[G]
\lesssim
	\dmax \SG[P^\U]
\lesssim
	\dmax \EC*[G] \cdot \dmax/\dmin.
\]
Fast mixing for low-degree graphs is thus characterised by the edge conductance $\EC*[G]$:
% via the uniform \RW:
\begin{enumerate}[noitemsep]
	\item 
		$G$ admits a fast mixing chain
	if and only if
		$P^\U$ is fast mixing;
	
	\item 
		$P^\U$ is fast mixing
	if and only if
		$1/\EC*[G]$ is polylogarithmic in $\abs V$.
\end{enumerate}

Such a simple characterisation does not hold if $\dmax$~is~large.
This may be slightly counter-intuitive at first:
	adding edges can only increase the optimum $\SG*$;
	but the lower bound above gets worse as $\dmax$ increases.
A striking example is given by taking two cliques on $n$ vertices and connecting them by a perfect matching; see \cref{fig:exm:matching} in \cref{sec:intro:exm:matching}.
It is a regular graph with
%edge conductance
$\EC* \asymp 1/n$,
but, as we will see later, it has
%optimal spectral gap
$\SG* \asymp 1$.
Informally, $\SG* \asymp 1$ because we can replace the two cliques with two bounded degree expander graphs without overly damaging its connectivity properties.
This shows that edge conductance is not the correct conductance measure for the \FMMC problem.

This prompts us to consider an alternative notion of connectivity: the \textit{vertex conductance} $\VC*$.
It measures how well connected a set is by comparing the number of \emph{vertices} in the boundary with its size.
Contrastingly,
edge conductance compares the number of edges in the boundary with the total number of edges inside the set.

\begin{introdefn}[Vertex Conductance]
\label{def:intro:main:vc:vtx}
The \textit{vertex conductance} $\VC*[G]$ of a finite graph $G = (V, E)$~is
%defined by
\[
	\VC*[G]
\cq
	\MIN{S \subseteq V : 0 < \abs S \le \abs V / 2}
	\VC[S]
\Qwhere
	\VC[S]
\cq
	\abs{ \partial S } / \abs{ S }
\Qfor
	S \subseteq V,
\]
where $\partial S$ is the \textit{vertex boundary} of $S \subseteq V$:
\[
	\partial S
\cq
	\brb{ y \notin S \mid \exists \, x \in S \st \bra{x,y} \in E }.
\]
\end{introdefn}

\end{subtheorem-num}

The example above in which two equisized cliques are connected by a perfect matching has vertex conductance $\VC* \asymp 1$. This agrees with our claimed optimal spectral gap $\SG* \asymp 1$.

Vertex conductance has been used to provide upper bounds
	on the time to spread a rumour in a graph
and
	on the hitting times of \RWs,
by \textcite{G:vertex-expansion,CRRST:electrical}, respectively, amongst others.
%	Giakkoupis, \green{Sauerwald and Sauerwald}
%	\cite{G:vertex-expansion,GSS:rumour-dynamic:conf}
%and
%	\textcite{CRRST:electrical},
%respectively, amongst others.
\textcite[Proposition~2]{R:fmmc-dual} showed that vertex conductance represents a fundamental barrier to fast mixing:
%	he showed that
	\(
		\SG*[G]
	\lesssim
		\VC*[G].
	\)
%	see \cite[Proposition~2]{R:fmmc-dual}.
This can be seen directly via a simple calculation comparing
the edge conductance
of any reweighing of $G$,
for which the \RW on this weighted graph has uniform equilibrium distribution,
with the vertex conductance.
%of the graph.

The edge and vertex conductances are comparable for low-degree graphs:
%precisely,
\[
	\EC*[G] \le \VC*[G] \le \dmax \EC*[G].
\]
Thus, the fact that the edge conductance $\EC*[G]$ characterises fast mixing for low-degree graphs means that the same holds for the vertex conductance $\VC*[G]$:
%precisely,
\[
	\VC*[G]^2 / \dmax^2
\lesssim
	\SG[P^\U]
\le
	\SG*[G]
\lesssim
	\VC*[G].
\]
We remove this $\dmax$ factor, at the cost of a $\log \abs V$ factor, thus showing that vertex conductance characterises the existence of a fast mixing chain \emph{for any graph}.
%Contrastingly, the clique perfect matching of \cref{fig:exm:matching} gave a counter-example to the claim for edge conductance when $\dmax$ is large.%
The graph of \cref{fig:exm:matching} in \cref{sec:intro:exm:matching} shows this does not hold for the edge conductance.  

%Moreover, a lower bound on
%%the optimal spectral gap
%$\SG*[G]$
%can be obtained in terms of the vertex conductance
%via Cheeger's inequality and the relations $\EC*[G] \le \VC*[G] \le \dmax \EC*[G]$:
%%by applying the discrete Cheeger inequality to the uniform \RW $P^\U$,
%\[
%	\VC*[G]^2 / \dmax^2
%\lesssim
%	\SG[P^\U]
%\le
%	\SG*[G]
%\lesssim
%	\VC*[G].
%\]
%These relations jointly characterise fast mixing if $\dmax$ is small, but not if it is large.

%\blue{I think the sentence leading to the theorem [in the old phrasing] makes it feel a bit underwhelming. We should try to emphasise a bit more the result.}
%\purple{I agree. We didn't used to mention edge conductance, so it was fairly important to explain where the vertex conductance lower bound came from. But since your suggestion of adding edge conductance earlier, the previous paragraph became a little reductant. I updated the edge conductance part too, in line with this new phrasing above. I really want to emphasise that low degree graphs are trivial---otherwise people might think, ``well, just removing some factor on the degree... low degree graphs are more natural... etc''}

%We establish bounds on
%%the optimal spectral gap
%$\SG*$
%which depend only on
%%the vertex conductance
%$\VC*$,
%up to a factor of $\log \abs V$:
%We remove any dependence on $\dmax$ in the lower bound,
%at the cost of a factor of $\log \abs V$.

\begin{introthm}[Characterisation of Fast Mixing]
\label{res:intro:main:vc}
	Let $G = (V, E)$ be a finite graph.
	Then
%	the optimal spectral gap
	$\SG*[G]$ satisfies%
	\[
		\VC*[G]^2 / \log \abs V
	\lesssim
		\SG*[G]
	\lesssim
		\VC*[G].
	\]
	Thus,
	vertex conductance characterises fast mixing for any graph.
\end{introthm}

%We give two remarks on this theorem and its proof.

%\begin{subtheorem-num}{intrormkT}
%	\label{rmk:intro:main:vc}

%\begin{intrormkt}[Tightness of Lower Bound]
	%
The quadratic dependence on the vertex conductance in the lower bound is needed for graphs such as the cycle.
	This has
		optimal spectral gap
		$\SG* \asymp 1/n^2$
	and
		vertex conductance
		$\VC* \asymp 1/n$;
	see \cite{BDSX:fmmc-path}.
We are not aware of a graph for which the $\log \abs V$ factor
%, or at least a factor $\log \dmax$,
is needed, but we have reasons to believe that such a factor,
or at least a factor $\log \dmax$,
is necessary. We elaborate.

\textcite{LRV:vertex} have essentially shown the following:
	under the so-called Small Set Expansion Conjecture of \textcite{RS:sse},
	for any $\eps > 0$,
	there is no polynomial-time algorithm that can distinguish between $\VC*[G] \le \eps$ and  $\VC*[G] \gtrsim \sqrt{\eps \log \dmax}$
	for any graph $G = (V,E)$.
Since the optimal spectral gap $\SG*[G]$ can be computed in polynomial time,
getting rid of the logarithmic factor altogether in \cref{res:intro:main:vc} would violate the Small Set Expansion Conjecture.
We leave open the problem of reducing the factor $\log \abs V$ to $\log \dmax$.

\smallskip

%\begin{intrormkt}[Matching Conductance]
	%
One of the most interesting aspects of the proof of \cref{res:intro:main:vc},
	given in \cref{sec:vc},
is that it does not directly relate the vertex conductance to the spectral gap.
Rather, it relates
	a variational characterisation of the optimal spectral gap,
	due to \textcite[Proposition~1]{R:fmmc-dual},
to
	a new connectivity measure for graphs which we introduce.
We term it \textit{matching conductance} and denote it $\MC*$.
It is defined similarly to vertex conductance, but it replaces
	the size of the vertex boundary of a set $S$
	in the numerator
with
	the size of a maximum matching between $S$ and $S^c$ in $E$.
A formal definition is given in \cref{def:vc:prelim:vc}.
It can be viewed as a measure of \emph{fault tolerance} of a graph:
	a graph has small matching conductance
if and only if
	we can remove a few vertices of the graph and split the graph into two large, disconnected subsets.

It turns out the matching conductance of a graph is essentially equivalent to its vertex conductance:
	\(
		\MC*[G]
	\asymp
		\VC*[G],
	\)
	uniformly over all graphs $G$;
	see \cref{res:vc:prelim:vertex-matching}.
A specific set $S$ of vertices, however, can have
	matching conductance $\MC[S]$
much smaller than its
	vertex conductance $\VC[S]$.
This fact makes using matching, rather than vertex, conductance essential in our proof.
	%
%\end{intrormkt}

%\end{subtheorem-num}

\smallskip

Finally, it is natural to ask for a statement analogous to \cref{res:intro:main:vc}, but for general $\pi$, rather than specific to $\pi = \unif[V]$.
This is not immediate, for slightly technical reasons.
We elaborate~in~\cref{sec:conc}.

\subsubsection{`Almost Mixing'}
\label{sec:intro:main:am}

\nextresult

We introduced the \FMMC problem to formalise our desire to construct a fast mixing Markov chain.
\cref{res:intro:main:vc}, however, implies there are certain graphs,
	namely those with small vertex conductance,
for which this desire cannot be attained.
It is then natural to ask if we can slightly relax the constraints we imposed to overcome this fundamental obstacle.

We answer this question affirmatively:
	we show that if the Markov chain is not required to have equilibrium distribution \emph{exactly} uniform, but only sufficiently close to uniform,
	then all graphs with small diameter admit a fast-mixing Markov chain.
Before formalising this claim, we gain some intuition by considering the following simple example, known as the \textit{dumbbell graph}.
\begin{itemize}[noitemsep]
	\item 
	Take two complete graphs $H_\pm = (V_\pm, E_\pm)$, each on $n$ vertices.
	Choose $v_\pm \in V_\pm$, respectively.
	
	\item 
	Form $\dumbbell = (V, E)$ by connecting both $v_\pm$ to a single `external' vertex $v_\star \notin V_+ \cup V_-$:
	\begin{center-small}
		$V \cq V_+ \cup V_- \cup \bra{v_\star}$
	and
		$E \cq E_+ \cup E_- \cup \bra{ \bra{v_+, v_\star}, \bra{v_-, v_\star} }$.
	\end{center-small}
\end{itemize}
Since $\dumbbell$ has vertex conductance equal to $1/n$,
\cref{res:intro:main:vc} implies that no chain with uniform equilibrium distribution can have relaxation time of smaller order than $n$.

In light of the above,
we propose the following \RW, described by a weighting on the edges~of~$\dumbbell$:
\begin{itemize}[noitemsep]
	\item 
	give all edges which do not include any of $\bra{v_+, v_-, v_\star}$ unit weight;
	
%	\item 
%	Set $\ww(e) \cq 1$ for all edges $e \in E$ which do not include any of $\bra{v_+, v_-, v_\star}$.
	
	\item 
	give the remaining edges weight $\eps n$.
	
%	\item 
%	Set $\ww\rbr{\bra{v_\pm, v}} \cq \eps n$ for all $v \in H_\pm$, respectively, and $\ww\rbr{\bra{v_+, v_\star}} \cq \eps n \eqqcolon \ww\rbr{\bra{v_-, v_\star}}$.
\end{itemize}
The \RW takes steps with distribution proportional to the edge weights.
It is straightforward to check that the equilibrium distribution induced is at most $\eps$ far from uniform in \TV.

The fundamental barrier to fast mixing in $\dumbbell$ is that any chain with uniform equilibrium gets stuck in one side of the graph for a time at least order $n$ in expectation.
Up-weighting the edges through the bottleneck means that the new \RW transitions between the two sides with expected time order $1/\eps$.
This leads to a relaxation time order $1/\eps$.
This all comes at a cost of having invariant distribution $\eps$ far from uniform.
Further details
%and generalisations
%including generalisation beyond having $H_\pm$ complete graphs,
are given in \S\ref{sec:intro:exm}.

We are able to generalise this construction to general graphs and general equilibrium distributions.
The fast-mixing Markov chain we design is a \RW on a carefully weighted breadth-first search (\textit{\BFS}) spanning tree, supplemented with self-loops.
We establish an upper bound of $12 (\diam G)^2 / \eps$ on the relaxation time when we allow the \RW to have equilibrium distribution $\eps$-far from uniform.

We now define precisely our set-up and formally state our result.

\begin{introdefn}[Almost-Mixed Distributions]
\label{def:intro:main:am:inv-dist}
	Let $V$ be a finite set.
	Let $\PD{V}$ denote the set of positive probability distributions on $V$.
	For $\eps \in [0, 1]$ and $\pi \in \PD{V}$,
	define
	\[
		\PD{\pi, \eps}
	\cq
		\brb{ \pi' \in \PD{V} \midb \mint{x \in V} \pi'(x) / \pi(x) \ge 1 - \eps }.
	\]
	In particular, if $\pi' \in \PD{\pi, \eps}$, then $\tv{ \pi' - \pi } \le \eps$.
\end{introdefn}

We actually establish a stronger result than the one described above.
The above description says that \emph{there exists some} reversible chain which is fast mixing:
	there exist $\pi' \in \PD{\pi, \eps}$ and $\PP \in \TM{G, \pi'}$ such that
	$\SG[\PP] \gtrsim \eps / (\diam G)^2$.
We prove that \emph{any} reversible chain can be perturbed into a fast mixing chain:
	for all $\pi \in \PD{V}$ and all $P \in \TM{G, \pi}$,
	there exist $\pi' \in \PD{\pi, \eps}$ and $\QQ \in \TM{G, \pi'}$ such that
	$\SG[\QQ] \gtrsim \eps / (\diam G)^2$ and $\QQ(e) \ge (1 - \eps) \PP(e)$ for all $e \in E$.

\begin{introthm}[Almost Mixing]
\label{res:intro:main:am}
	Let $G = (V, E)$ be a finite, connected graph and $\pi \in \PD{V}$.
	Let $\eps \in (0, 1)$ and $\PP \in \TM{G, \pi}$.
	There exist $\pi' \in \PD{\pi, \eps}$ and $\QQ \in \TM{G, \pi'}$ such that
	\[
		\SG[\QQ]
	\ge
		\tfrac1{12} \eps (\diam G)^{-2}
	\Qand
		\QQ(e) \ge (1 - \eps) \PP(e)
	\Qforall
		e \in E.
	\]
	A consequence of this spectral gap estimate is that
	\[
		\MIX[\QQ]
	\le
		24 \eps^{-1} (\diam G)^2 \log\rbr{ 16 \pimin^{-1} }.
	\]
\end{introthm}

The matrix $\QQ$ is obtained as a perturbation of $\PP$.
Moreover, this perturbation is actually independent of $\PP$:
	we construct a weighted \BFS tree, as describe above, and `superimpose' it with the weights corresponding to $\PP$.
A more refined statement, making this independence of the perturbation explicit, is given in \cref{res:am:res:main}.

\smallskip

This diameter bound is a substantial improvement over the vertex conductance lower bound on the optimal spectral gap from \cref{res:intro:main:vc}.
It comes at the cost of having invariant distribution only `almost' uniform---hence the name ``almost mixing''.
We show in the next section that passing to the continuous-time setting allows this diameter-squared bound to be maintained while having \emph{exactly} uniform invariant distribution.
We use fundamentally the same chain:
	it is a \RW on the same weighted \BFS tree,
	where the weights now represent the rate at which an edge is crossed.

\subsubsection{Continuous-Time Markov Chains}
\label{sec:intro:main:cts}

The discussion and results above all concern \emph{discrete-time} Markov chains.
It is also natural to study the question of the fastest mixing Markov chain in \emph{continuous-time}.
We restrict to the case where the target distribution $\pi = \unif[V]$.
The question of the \FMMC in continuous-time was originally raised by \textcite{SBXD:fmmc-cts} and has been studied subsequently by \textcite{S:fmmc-thesis,MT:mixing-times}.
We review their work~in~\cref{sec:intro:previous-work}.

A continuous-time Markov chain on a graph $G = (V, E)$ with uniform equilibrium distribution can be represented by the \RW on a weighted graph $(G, \qq)$, where $\qq : E \to \mbr_+$, as follows.

\begin{subtheorem-num}{introdefn}
	\label{def:intro:main:cts}

\begin{introdefn}[\RW on Weighted Graph]
Let $G = (V, E)$ be a graph and $\qq : E \to \mbr_+$ a collection of non-negative weights.
The \RW on $(G, \qq)$ jumps from $x$ to $y$ at rate $\qq(\bra{x,y})$ for $x,y \in V$ with $\bra{x,y} \in E$.
The \textit{Laplacian} $\LL^\qq \in \mbr^{V \times V}$ of the weighting $\qq$ is defined by
\[
	\LL^\qq_{x,y}
\cq
	\one{\bra{x,y} \in E}
\cdot
	\qq\rbr{\bra{x,y}}
-	\one{x = y}
 	\sumt{z \in V : \bra{x,z} \in E}
	\qq\rbr{\bra{x,z}}
\Qfor
	x,y \in V.
\]
The \textit{spectral gap}, which we denote $\SG[\qq]$, is given by the second smallest eigenvalue of $\LL^\qq$.
\end{introdefn}

The spectral gap $\SG[\qq]$ is intrinsically related to the mixing time $\MIX[\qq]$, as in discrete-time:%
\[
	\SG[\qq]^{-1}
\lesssim
	\MIX[\qq]
\lesssim
	\SG[\qq]^{-1} \cdot \log \abs V
\Qwhere
	\MIX[\qq]
\cq
	\MIX[\qq][\tfrac14];
\]
see \cite[Lemma~4.23]{AF:book}.
Again, \emph{fast mixing} means relaxation time $\SG[\qq]^{-1}$ polylogarithmic in $\abs V$.
The above relations thus imply that a polylogarithmic relaxation time characterises fast mixing.

It is immediate to see that if all the rates are multiplied by some factor, then the spectral gap changes by that factor too:
\(
	\SG[c \qq]
=
	c \SG[\qq]
\)
for any $c > 0$.
We must therefore impose some~normalisation.

\begin{introdefn}[Normalisation]
The rate at which the walk leaves the vertex $x$ is given by
\[
	\qq(x)
\cq
	\sumt{y \in V}
	\one{\bra{x,y} \in V}
	\qq\rbr{\bra{x,y}}
\Qfor
	x \in V.
\]
We call
	$\maxt{x \in V} \qq(x)$
	the \textit{maximal} leave-rate
and
	$\abs V^{-1} \sumt{x \in V} \qq(x)$
	the \textit{average} leave-rate.
\end{introdefn}

\end{subtheorem-num}

A natural normalisation is to require a \emph{maximal} leave-rate of $1$.
It can be seen, however, that this reduces to the discrete-time case via exponential-$1$ waiting times.
We impose instead an \emph{average} leave-rate of $1$, or, equivalently, $\qq(E) \le \tfrac12 \abs E$.
This allows a few vertices to have abnormally large leave-rate, but rarely enough that the average is not significantly affected. This will allow the \RW to exit small `bottlenecks' quickly, where the discrete-time walk would remain stuck for significant time.
This average leave-rate normalisation was considered in \cite{SBXD:fmmc-cts,S:fmmc-thesis,MT:mixing-times}.
%A slightly more general normalisation is considered in \cite{SBXD:fmmc-cts}, but the principle is the same.
\textcite[\S 7.1]{MT:mixing-times} describe this normalisation as ``rather powerful [compared with discrete-time]'' due to the fact that the maximal leave-rate may be very large.

The main result of this section states that, for any graph, it is possible to construct a weighting with average leave-rate of $1$ such that its spectral gap depends only on the diameter of the graph.
% can be summarised as follows:
%	given a graph $G = (V, E)$,
%	we construct
%	a weighting
%		$\ww : E \to \mbr_+$
%	with
%		$\ww(E) \le \tfrac12$
%	such that
%		$\SG[\ww] \lesssim (\diam G)^2$.
	
\begin{introthm}[Continuous-Time]
\label{res:intro:main:cts}
	Let $G = (V, E)$ be a graph.
	There exists a weighting $\ww : E \to \mbr_+$ with average leave-rate at most $1$ such that the \RW on $(G, \ww)$
	satisfies
%	has spectral gap $\SG[\ww]$ and mixing time $\MIX[\ww] = \MIX[\ww][\tfrac14]$ satisfying
	\[
		\SG[\ww]
	\ge
		\tfrac1{16} (\diam G)^{-2}
	\Qand
		\MIX[\ww]
	\le
		8 (\diam G)^2 \log\rbr{ 16 n }.
	\]
\end{introthm}

An upper bound on $1/\SG*[G]$ of order $(\diam G)^2$ is required for graphs with diffusive behaviour, such as the cycle or the path.
A lower bound of order $\diam G$, however, \emph{is not} necessary, in general. This is in stark contrast to the discrete-time case.
Indeed, in continuous-time, a few edges can be up-weighted significantly with little affect on the average.
So if the `typical' distance is much less than the maximal, a relaxation time of smaller order than the diameter may well be achievable.
%More details are given in the next paragraph.
%See the next paragraph for more details.

\smallskip

A highly related theorem is given in \citeauthor{S:fmmc-thesis}'s PhD thesis \cite[\S 3.3]{S:fmmc-thesis};
%, supervised by Tetali.
see also \cite[\S 7.1]{MT:mixing-times}.
They bound
the optimal spectral gap
$\SG*[G]$ in terms of the \textit{spread constant} $c(G)$, introduced in \cite{ABS:spread-constant}, which is the maximal variance of a function that is Lipschitz on the edges of $G$:
\[
	2 c(G)
\le
	1/\SG*[G]
\lesssim
	c(G) \log \abs V.
\]
The spread constant $c(G)$ can be upper bounded by $\tfrac14 (\diam G)^2$, but there are examples for which this is far from tight.
Still, if a very general, easy to calculate, bound is desired, then we do not know of a better bound than $c(G) \lesssim (\diam G)^2$, which reduces to approximately our bound.
The spread constant $c(G)$ can also be lower bounded by a type of `typical' distance; see \cite[Corollary~7.2]{MT:mixing-times}.

In contrast with our result, however, theirs is non-constructive, relying on the famous, but non-constructive, Johnson--Lindenstrauss lemma \cite{JL}.
%Our result is differs from theirs in a key way:
%	it is constructive, whilst theirs is not.
%Theirs uses the famous, but non-constructive, Johnson--Lindenstrauss lemma \cite{JL} to reduce a high-dimensional optimisation problem to a $1$-dimensional one.
\textcite[Remark~7.3]{MT:mixing-times} comment on the difficulty of explicitly constructing such a process:
	``It might be challenging and in general impractical a task to actually find such a process explicitly.''
Our construction is explicit and can actually be constructed in time linear in the size of the graph.

\smallskip

\textcite[Remark~7.3]{MT:mixing-times} also comment on the existence of such a fast mixing Markov chain in continuous-time:
%such as theirs or \cref{res:intro:main:cts}:
	``The key [to the existence of such a chain] ... might be that we were implicitly providing the continuous-time chain with more power ... by not requiring the rates in each row to sum to $1$, but only the [average rate to be $1$].''
This significant additional power allows bottlenecks to be traversed quickly while maintaining an \emph{average} leave-rate of $1$.
Indeed, the weighting $\ww$ that we construct has $\max_{x \in V} \ww(x) \asymp n / \diam G$, which may be far larger than $1$.

This really emphasises the strength of our `almost mixing' result, \cref{res:intro:main:am}:
	the chain there is in discrete-time%
	---or, equivalently, has $\max_{x \in V} \qq(x) \le 1$---%
	but still attains a spectral gap only order $\eps$ smaller than that attained in the continuous-time case of \cref{res:intro:main:cts}.
Of course, the cost is that the equilibrium distribution $\pi'$ only satisfies $\min_{x \in V} \pi'(x) / \pi(x) \ge 1 - \eps$, not $\pi' = \pi$.

\smallskip

We expect that our continuous-time analysis can be adjusted to handle general equilibrium distributions $\pi$ with relatively little changed.
We have not checked the details, however.
We focussed on the uniform case because it is, arguably, the most important and the cleanest to present.

\subsubsection{Time-Inhomogeneous Markov Chains}
\label{sec:intro:main:tinhom}

Our attention has been so far restricted to \emph{time-homogeneous} Markov chains, in which the transition probabilities do not change over time and are described by a single transition matrix $P$.
A \textit{time-inhomogeneous} Markov chain, instead,
%$(X_t)_{t \in \mbn}$
is described by a sequence $(P_t)_{t \in \mbn}$ of transition matrices and an initial law $\mu_0$:
	the time-$t$ law
	\(
		\mu_t
	\cq
		\pr{X_t \in \cdot}
	\)
	is given by
	\(
		\mu_t
	=
		\mu_0 P_1 P_2 \cdots P_t
	\Qfor
		t \in \mbn.
	\)
A \textit{time-homogeneous} chain has $P_t = P$ for all $t \in \mbn$, for some $P$. 
%Such a chain is somewhat less natural from an implementation point of view. \blue{I would cut this sentence: it breaks the flow.}
Nevertheless, we close our section of results by showing that they can lead to improvements over time-homogeneous chains.
%Namely, we can sample from any distribution in at most $2 \diam G$ steps, \green{using only transitions permitted by $G$}.

\begin{introthm}%[Time-Inhomogeneous]
\label{res:intro:main:tinhom}
	Let $G = (V, E)$ be a connected graph and let $\pi \in \PD{V}$.
	There exists a time-inhomogeneous Markov chain on $G$ that perfectly mixes to $\pi$ after $2 \diam G$ steps:
	\(
		\mu_{2 \diam G}
	=
		\pi.
	\)
\end{introthm}

It is easy to see that $\diam G$ is a lower bound on the fastest `perfectly mixing' chain.
If one only requires
\(
	\tv{ \pr{X_t \in \cdot} - \pi }
<
	\tfrac12,
\)
then $\tfrac12 \diam G$ is a lower bound.
Thus the bound of $2 \diam G$
above
%from \cref{res:intro:main:tinhom}
is tight up to a factor of at most $4$.

\endgroup

%\titleformat{\subsubsection}
%	{\sffamily \normalsize \bfseries \boldmath}{\thesubsubsection}{1em}{}
%
%\renewcommand{\thesubsubsection}{\thesubsection.\arabic{subsubsection}}

\subsection{Notable Examples}
\label{sec:intro:exm}

We discuss briefly a few examples which are of particular interest.
We always consider the uniform distribution, i.e.\ $\pi = \unif[V]$, unless specified to the the contrary.

\subsubsection*{Dumbbell Graph}
\label{sec:intro:exm:dumbbell}

Let $\dumbbell$ be the dumbbell graph with bells $H_\pm$ of size $n$.
The bells $H_\pm$ need not be cliques $K_n$; they can be arbitrary connected graphs on $n$ vertices.
See \cref{fig:exm:dumbbell} for an illustration when $H_\pm = K_n$.
%This is the graph $\dumbbell$ described in the motivation for `almost mixing'.

\tikzset{%
	mynode3/.style={%
		circle, minimum size=4pt,
		inner sep=0pt, outer sep=0pt, draw=darkgray, fill=white
	}
}

\begin{figure}
%	\centering
\begin{minipage}{\textwidth}
\begin{minipage}[b]{0.475\textwidth}
	\centering

\begin{tikzpicture}

%\draw[gray,fill=lightgray, opacity=0.2] (0,0) circle (1.3cm);
%\draw[gray,fill=lightgray, opacity=0.2] (3,0) circle (1.3cm);

%rotate = 321.5
\node[minimum size=3cm, regular polygon, regular polygon sides=7, rotate=-90] (epta) at (0,0) {};
    \foreach \x in {1,2,...,7}{%
        \node[mynode3] at (epta.corner \x) (e\x) {};
    }
   
    \foreach \x in {1,2,...,7}{%
        \foreach \y in {1,2,...,7}{%
            \ifthenelse{\x>\y}{}{\draw[black] (e\x) -- (e\y);}
         }       
    }
    
\node[mynode3] at (2,0) (c) {};

\node[minimum size=3cm, regular polygon, regular polygon sides=7, rotate=90] (epta2) at (4,0) {};
    \foreach \x in {1,2,...,7}{%
        \node[mynode3] at (epta2.corner \x) (f\x) {};
    }

    \foreach \x in {1,2,...,7}{%
        \foreach \y in {1,2,...,7}{%
            \ifthenelse{\x>\y}{}{\draw[black] (f\x) -- (f\y);}       
    }}

   	\draw[black] (c) -- (e1);
	\draw[black] (c) -- (f1);

%\node[align = center] at (1.5,-1.7) {$1-\lambda(\overline{P})=\Omega(1)$};
% \node[align = center, below = of image] {Odd $t$};

\end{tikzpicture}%
\captionof{figure}{%
	Dumbbell graph $\dumbbell$ with $n = 7$:
	two cliques connected to a single external vertex%
}

\label{fig:exm:dumbbell}
\end{minipage}
%	\quad
\hfill
%	\quad
\begin{minipage}[b]{0.475\textwidth}
	\centering

\begin{tikzpicture}

%\draw[gray,fill=lightgray, opacity=0.2] (0,0) circle (1.3cm);
%\draw[gray,fill=lightgray, opacity=0.2] (3,0) circle (1.3cm);

\node[minimum size=3cm, regular polygon, regular polygon sides=7, rotate=0] (octa) at (0,0) {};
    \foreach \x in {1,2,...,7}{%
        \node[mynode3] at (octa.corner \x) (e\x) {};
    }
    \node[mynode3] at (octa.center) (e0) {};
   
    \foreach \x in {1,2,...,7}{%
       \draw[black] (e0) -- (e\x);
    }

%\node[minimum size=3cm, regular polygon, regular polygon sides=7, rotate=90] (octa) at (0,0) {};
%    \foreach \x in {1,2,...,7}{%
%        \node[mynode3] at (octa.corner \x) (e\x) {};
%    }
%    \node[mynode3] at (octa.center) (e0) {};
%   
%    \foreach \x in {1,2,...,7}{%
%       \draw[black] (e0) -- (e\x);
%    }

%\node[align = center] at (1.5,-1.7) {$1-\lambda(\overline{P})=\Omega(1)$};
% \node[align = center, below = of image] {Odd $t$};

\end{tikzpicture}%
\captionof{figure}{%
	Star graph $\stargraph$ with $n = 7$
	a central vertex connected to leaves%
}

\label{fig:exm:star}
\end{minipage}
\end{minipage}

\bigskip

\begin{minipage}{\textwidth}
\begin{minipage}[b]{0.475\textwidth}
	\centering

\begin{tikzpicture}

%\draw[gray,fill=lightgray, opacity=0.2] (0,0) circle (1.3cm);
%\draw[gray,fill=lightgray, opacity=0.2] (3,0) circle (1.3cm);

\node[minimum size=3cm, regular polygon, regular polygon sides=7, rotate=-90] (epta) at (0,0) {};
    \foreach \x in {1,2,...,7}{%
        \node[mynode3] at (epta.corner \x) (e\x) {};
    }
   
    \foreach \x in {1,2,...,7}{%
        \foreach \y in {1,2,...,7}{%
            \ifthenelse{\x>\y}{}{\draw[black] (e\x) -- (e\y);}       
    }} 

\node[minimum size=3cm, regular polygon, regular polygon sides=7, rotate=90] (epta2) at (4,0) {};
    \foreach \x in {1,2,...,7}{%
        \node[mynode3] at (epta2.corner \x) (f\x) {};
    }

    \foreach \x in {1,2,...,7}{%
        \foreach \y in {1,2,...,7}{%
            \ifthenelse{\x>\y}{}{\draw[black] (f\x) -- (f\y);}       
    }}

%     \foreach \x in {1,2,...,7}{
%     	\draw[black] (e\x) -- (f\x);
%     } 

	\draw[black] (e1) -- (f1);
	\draw[black] (e7) -- (f2);
	\draw[black] (e6) -- (f3);
	\draw[black] (e5) -- (f4);
	\draw[black] (e4) -- (f5);
	\draw[black] (e3) -- (f6);
	\draw[black] (e2) -- (f7);

%\node[align = center] at (1.5,-1.7) {$1-\lambda(\overline{P})=\Omega(1)$};
% \node[align = center, below = of image] {Odd $t$};
 
\end{tikzpicture}%
\captionof{figure}{%
	Matching graph $\matching$ with $n = 7$:
	two cliques connected via a matching%
}

\label{fig:exm:matching} 
\end{minipage}
%	\quad
\hfill
%	\quad
\begin{minipage}[b]{0.475\textwidth}
	\centering

\begin{tikzpicture}
% \node[align = center] at (1.5,1.65) {Average transition probabilities $\overline{P}$};

%\draw[gray,fill=lightgray, opacity=0.2] (0,0) circle (1.3cm);
%\draw[gray,fill=lightgray, opacity=0.2] (3,0) circle (1.3cm);

\node[minimum size=3cm, regular polygon, regular polygon sides=7, rotate=-90] (epta) at (0,0) {};
    \foreach \x in {1,2,...,7}{%
        \node[mynode3] at (epta.corner \x) (e\x) {};
    }
   
    \foreach \x in {1,2,...,7}{%
        \foreach \y in {1,2,...,7}{%
            \ifthenelse{\x>\y}{}{\draw[black] (e\x) -- (e\y);}       
    }} 

\node[minimum size=3cm, regular polygon, regular polygon sides=7, rotate=90] (epta2) at (4,0) {};
    \foreach \x in {1,2,...,7}{%
        \node[mynode3] at (epta2.corner \x) (f\x) {};
    }

    \foreach \x in {1,2,...,7}{%
        \foreach \y in {1,2,...,7}{%
            \ifthenelse{\x>\y}{}{\draw[black] (f\x) -- (f\y);}       
    }}

   	\draw[black] (f1) -- (e1);
	\draw[black] (f1) -- (e2);
	\draw[black] (f1) -- (e7);

%\node[align = center] at (1.5,-1.7) {$1-\lambda(\overline{P})=\Omega(1)$};
%\node[align = center, below = of image] {Odd $t$};

\end{tikzpicture}%
\captionof{figure}{%
	Source graph $\source$ with $n = 7$ and $k = 3$:
	two cliques connected via a `source'%
}

\label{fig:exm:source}
\end{minipage}
\end{minipage}

\end{figure}

\paragraph*{Conductance Measures.}
It is straightforward to see that the set with the worst vertex conductance is given by one side of the dumbbell graph: $S = H_-$ or $S = H_+$.
This shows that
\[
	\VC*[\dumbbell]
=
	\VC[H_\pm]
=
	1/n.
\]
This implies that the optimal relaxation time $1/\SG*[\dumbbell]$ satisfies
\(
	n \lesssim 1/\SG*[\dumbbell] \lesssim n^2 \log n.
\)

It is easy to find a chain attaining the correct order of $n^{-1}$ when $H_\pm = K_n$, i.e.\ each bell is a complete graph on $n$ vertices.
Define a weighting as follows.
Each vertex gets the same total~weight.
\begin{itemize}[noitemsep]
	\item 
	Place unit weights on all edges which do not include the centre $v_\star$.
%	This gives weight $n-1$ to each vertex of $V \setminus \bra{v_\star}$.
	
	\item 
	Give the edges $\bra{v_-, v_\star}$ and $\bra{v_+, v_\star}$ weight $n-1$.
%	The vertices $\bra{v_-, v_+, v_\star}$ now have weight $2(n-1)$.
	
	\item 
	Add self-loops of weight $2(n-1)$ to each of the vertices in $V \setminus \bra{v_-, v_+, v_\star}$.
%	These now also have weight $2(n-1)$.
\end{itemize}
%In particular, each vertex has the same weight.

The probability of stepping to $v_\star$ from either of $v_\pm$ is $\tfrac12$.
This gives an order-$n$ hitting time of one clique from the other.
This implies that our chain has relaxation time order $n$.
Contrast this with the suboptimal order $n^2$ hitting and relaxation time for the uniform \RW.

\paragraph*{Almost Mixing.}
If the graphs $H_\pm$ have polylogarithmic diameter
then \cref{res:intro:main:am} provides a chain with polylogarithmic relaxation time.
This is a substantial improvement from linear.
This is true regardless of the particular structure of $H_\pm$:
	it just needs $\log \diam H_\pm \lesssim \log \log n$.
If $\diam H_\pm \asymp 1$, then we obtain relaxation time order $1$, known also as an expander.

We now explain very roughly how to construct this chain for this dumbbell example $\dumbbell$.
The general idea is to up-weight edges towards the central vertex $v_\star$, which is the bottleneck. We do this is such a way to make the distance to $v_\star$ behave somewhat like an unbiased \RW on $\mbz$. This way it should take time order $(\diam H_\pm)^2$ to move from $H_\pm$ to $H_\mp$.

It is natural to try to achieve this bias by rooting a spanning tree $T$ at $v_\star$ and then up-weighting the vertices towards the root. This leads to a worst-case hitting time for the root $v_\star$ of order $(\diam T)^2$.
We choose $T$ to be a breadth-first search (\textit{\BFS}) tree since this has $\diam T \le 2 \diam G$.

We give a more detailed overview in \cref{sec:am:outline}.
We specifically chose the bottleneck vertex $v_\star$ to be the root of $T$ above.
It turns out that actually any choice of root suffices.
The reader may find this surprising at first; we did.
%But the paths from one side to the other all need to go through $v_\star$, so this path automatically has its weights increased, regardless of the choice of root.
More generally, suppose that $o \in V$ is any vertex and a \BFS is rooted at $o$; we up-weight the edges towards $o$. Paths from $v \ne o$ to $o$ naturally go through bottlenecks. This automatically up-weights edges in bottlenecks.

\subsubsection*{Binary Tree}
\label{sec:intro:exm:binarytree}

Let $\binarytree = (V, E)$ be the complete binary tree on $n = 2^N - 1$ vertices with depth $N \approx \log_2 n$.
%See \cref{fig:exm:binarytree} for an illustration.

\paragraph*{Conductance Measures.}

It is straightforward to see that the set with the worst vertex conductance is given by one side of the tree:
	the root,
	a child
and
	all its descendants
This gives
\[
	\VC*[\binarytree]
=
	1/n.
\]
This implies that the optimal relaxation time $1/\SG*[\binarytree]$ satisfies
\(
	n \lesssim 1/\SG*[\binarytree] \lesssim n^2 \log n.
\)

$\binarytree$ has bounded degree, so the maximum degree chain attains the correct order relaxation time.
This chain is just the simple \RW, but with extra laziness at the leaves to make the invariant distribution uniform.
The correct relaxation time is order $n$; see \cite{S:spectral} for details.

\paragraph*{Almost Mixing.}

It is very natural to root the \BFS tree at the root $o$ of $\binarytree$.
The up-weighting will help pull the walk up the tree towards the root, allowing it to spread across the width of the tree more easily.
The weight given to the edge from $x \ne o$ to its parent is given by the number of vertices in the subtree rooted at $x$. Precisely, if $x$ is at distance $d \ge 1$ from the root, then the weight is $2^{N-d}-1$.
The up-weighting means that the distance from the root behaves roughly as an unbiased \RW.
The hitting time of the root is then order $N^2 \asymp (\log n)^2$.
Once the \RW hits the root, which branch it takes after is uniformly distributed.
Thus, once it hits the leaves again, it is uniform on the leaves and so approximately mixed.
The total time for this is order $N^2 \asymp (\log n)^2$.

Our method does not know, however, that there is anything special about the root of the binary tree. \emph{Any} vertex can be picked as the root of the \BFS tree.
Viewed from this vertex, $\binarytree$ is like a complete binary tree of some depth, but with some of the branches pruned.
The depth is at most $2N$.
The same ideas give relaxation and mixing time order $N^2 \asymp (\log n)^2$.

\subsubsection*{Star Graph}
\label{sec:intro:exm:star}

Let $\stargraph = (V, E)$ be the star graph with centre $v_\star$ and $n$ leaves.
See \cref{fig:exm:star} for an illustration.

\paragraph*{Conductance Measures.}

There is a simple dichotomy for the vertex boundary of a set $S \ne \emptyset$ in the star graph $\stargraph$:
	if $v_\star \in S$,
		then $\partial S = S^c$;
	if $v_\star \notin S$,
		then $\partial S = \bra{v_\star}$.
Thus any $S \subseteq V$ with $\abs S = \floor{\tfrac12 n}$ and $v_\star \notin S$ satisfies $\VC[S] = 1/\abs S \asymp 1/n$.
It is straightforward to see that this gives the correct order:
	place unit weights on all the edges
and
	weight-$(n-1)$ self-loops on all the non-central vertices;
	it is easy to see that this chain has mixing time order $n$.

Another measure of vertex conductance replaces the $\partial S$ with the symmetric union $\partial_\mathrm{sym} S \cq \partial S \cup \partial S^c$; denote the vertex conductance with $\partial_\mathrm{sym}$ by $\Psi_\mathrm{sym}$.
Again, there is a simple dichotomy for $S \ne \emptyset$:
	if $v_\star \in S$,
		then $\partial_\mathrm{sym} S = S^c \cup \bra{v_\star}$;
	if $v_\star \notin S$,
		then $\partial_\mathrm{sym} S = S \cup \bra{v_\star}$.
Thus $\Psi_\mathrm{sym}^\star \asymp 1$.

The difference between the two measures is that if the boundary of $S$ is small, then that of $S^c$ is large.
The use of $\partial S$, as opposed to $\partial_\mathrm{sym} S$, is thus important in \cref{res:intro:main:vc}.

\smallskip

It is well-known that the spectral gap $\SG$
%, often denoted $\lambda_2$,
is characterised by a variational form.
The relationship between
spectral gap $\SG$
%$\lambda_2$
and the edge conductance $\EC*$ of the simple \RW
is given by the well-known Cheeger inequalities.
A related variational form was introduced by \textcite{BHT:lambda-inf}, which they denote $\lambda_\infty$.
They establish various Cheeger-type relationships between $\lambda_\infty$ and the vertex conductance $\Psi_\mathrm{sym}^\star$.
In particular, they show, for any graph, that
\[
%	(\VC*)^2
	(\Psi_\mathrm{sym}^\star)^2
\lesssim
	\lambda_\infty
\lesssim
	\Psi_\mathrm{sym}^\star;
\]
see \cite[Theorems~1 and~2]{BHT:lambda-inf}.
It is immediate from the definitions that
%\(
%	\SG*
%\lesssim
%	\lambda_\infty
%\lesssim
%	\dmax \SG*.
%\)%
\(
	\lambda_\infty / \dmax
\lesssim
	\SG*
\lesssim
	\lambda_\infty.
\)

In light of \cite{BHT:lambda-inf} and our \cref{res:intro:main:vc},
it is natural to wonder whether $\lambda_\infty$ can be directly related to the optimal spectral gap $\SG*$, without a $\dmax$ factor.
The example of the star graph $\stargraph$ shows that this is not possible:
$\Psi_\mathrm{sym}^\star \asymp 1$ and thus $\lambda_\infty \asymp 1$, but $\VC* \asymp 1/n$ and so $\SG* \lesssim 1/n$.
This shows that $\lambda_\infty$ is really not the correct parameter for the \FMMC problem.
%In fact, we showed that $\SG* \asymp 1/n$ via the exhibition of a simple chain with spectral gap order $1/n$.

%We do not need to use the upper bound $\SG* \lesssim \lambda_\infty \lesssim \Psi_\mathrm{sym}^\star$; rather, we can use $\SG* \lesssim \VC*$ from before.
%Precisely,
%we can combine \cite[Theorems~1 and~2]{BHT:lambda-inf} with 
%\(
%	\SG*
%\lesssim
%	\VC*
%\)
%to obtain
%\[
%	(\Psi_\mathrm{sym}^\star)^2 / \dmax
%\lesssim
%	\lambda_\infty / \dmax
%\lesssim
%	\SG*
%\lesssim
%	\VC*.
%\]
%All these bounds are now tight for the star graph.
%However, the star graph is somewhat unusual, albeit not unique, in having $\VC*$ so much smaller than $\Psi_\mathrm{sym}^\star$.
%These quantities are the same order for many graphs, in which case the above display is no improvement on
%%the almost trivial bound
%\(
%	(\VC*)^2 / \dmax
%\lesssim
%	\SG*
%\lesssim
%	\VC*.
%\)

%\blue{In light of the above result and our \cref{res:intro:main:vc},
%it is then natural to wonder if $\lambda_\infty$ could be directly related to the optimal spectral gap $\SG$. Our discussion on the star graph, however, shows it cannot be the case. Indeed, since $\Psi_\mathrm{sym}^\star \asymp 1$, also $\lambda_\infty \asymp 1$, while $\SG^\star \lesssim \VC^\star \lesssim 1/n$.}

\paragraph*{Almost Mixing.}

Obtaining an `almost mixing' chain with order $1/\eps$ mixing time is simple:
	place $\eps$ weights on all the edges
and
	weight-$1$ self-loops to all non-central vertices.
The total weight of the central vertex is $\eps (n-1)$.
The remainder of the weight is spread uniformly.
Thus the distribution $\pi$ induced by this weighting is in $\PD{\stargraph, \eps}$.
It is easy to see that the mixing time is order $1/\eps$.

\subsubsection*{Complete Graphs Connected via a Matching}
\label{sec:intro:exm:matching}

Let $\matching = (V, E)$ be two $n$-cliques $H_\pm = (V_\pm, E_\pm)$, connected via a matching:
%, rather than a single edge:
	enumerate $V_\pm$ as $\bra{v_\pm^1, ..., v_\pm^n}$
and
	connect $v_+^i$ with $v_-^i$ for all $i \in \bra{1, ..., n}$.
See \cref{fig:exm:matching} for an illustration.

%\purple{The analysis in this section is extremely easy. I wondered whether it's worth extending it to general expanders in the `bells'. I've written the details, but I'm not sure how much we should include---if any! We can always just state that $\VC*[H_\pm] \asymp 1 \implies \VC*[\matching] \asymp 1$ and say that we omit the details}

\paragraph*{Conductance Measures.}

If $v_+^i \in S$ but $v_-^i \notin S$, then $v_-^i \in \partial S$.
From this, $\VC*[\matching] \asymp 1$ follows~easily.
It is not too hard to show that the optimal spectral gap $\SG*[\matching]$ is of the same order:
	just replace the two cliques by two $3$-regular expanders and leave the perfect matching in place.
The lazy simple \RW on this edge-induced subgraph of $\matching$ has order-$1$ spectral gap.

\paragraph*{Almost Mixing.}

The optimal spectral gap is order $1$, so there is no need for `almost mixing'.

%The optimal spectral gap is order $1$ when $H_\pm$ are expanders, so there is no need for `almost mixing'.
%More generally, there are a variety of different cases which can be constructed.
%Certainly as a worst case,
%\(
%	\diam \matching
%\le
%	\diam H_+ + \diam H_- + 1;
%\)
%this inequality is tight, up to constants, for the dumbbell graph, but may not be remotely tight for the matching graph.
%%but this may not be remotely tight.
%%, depending on the graph.

\subsubsection*{Complete Graphs Connected via a `Source'}

Let $\source = (V, E)$ be two $n$-cliques $H_\pm$ as `bells', connected via a `source':
	choose $v_0 \in H_-$ and $v_1, ..., v_k \in H_+$;
	connect $v_0$ with each of $\bra{ v_i }_{i=1}^k$.
See \cref{fig:exm:source} for an illustration.

\paragraph*{Conductance Measures.}

One may think at first that this `source' of $k$ edges, rather than just a single edge, give rise to faster mixing; indeed, $\VC[H_-] = k/n$.
However, removing the source from the set gives $\VC[H_- \setminus v_0] = 1/(n-1)$.
So in fact the vertex conductance of the source graph $\source$ is almost the same as that of the dumbbell graph $\dumbbell$.

The edge conductance of the graph does improve with $k$:
	$\EC[\source] \asymp k/n^2$.
But this is always at most $1/n$.
So the improvement from $k$ is not enough to outweigh the fact that the uniform \RW has spectral gap far from the optimal---unless $k \asymp n$.

An optimal spectral gap can be achieved by choosing an arbitrary $3$-regular expander as a subgraph of each of the cliques and connecting these via a single edge.
The uniform \RW on this sparse subgraph then has spectral gap order $1/n$.

\paragraph*{Almost Mixing.}

We can use exactly the same construction as in the dumbbell graph $\dumbbell$, picking an arbitrary edge amongst the $k$ connecting edges.

%We choose an arbitrary vertex $o \in V$ and root a \BFS tree at $o$, as always.
%The two sources are natural bottlenecks in the graph. Thus many of the paths from vertices $v \ne o$ to $o$ pass through these. This up-weights the edges in the bottleneck.
%This allows the achievement of a chain with relaxation time order
%\(
%	(\diam \source)^2
%\asymp
%	1,
%\)
%which is order $1$ if $H_\pm$ are complete graphs.\footnote{%
%	Can we compare $\diam \source$ with $\diam H_\pm$ in any nice way for general $H_\pm$?}

\subsection{Review of Previous Work}
\label{sec:intro:previous-work}

\newcommand{\diac}[1][]{\cite[#1]{SBXD:fmmc-cts}\xspace}
\newcommand{\diat}[1][]{\textcite[#1]{SBXD:fmmc-cts}\xspace}
\newcommand{\diaa}{\citeauthor{SBXD:fmmc-cts}\xspace}
\newcommand{\mtc}[1][]{\cite[#1]{MT:mixing-times}\xspace}
\newcommand{\mtt}[1][]{\textcite[#1]{MT:mixing-times}\xspace}
\newcommand{\mta}{\citeauthor{MT:mixing-times}\xspace}
\newcommand{\samc}[1][]{\cite[#1]{S:fmmc-thesis}\xspace}
\newcommand{\samt}[1][]{\textcite[#1]{S:fmmc-thesis}\xspace}
\newcommand{\sama}{\citeauthor{S:fmmc-thesis}\xspace}

We now review previous related work.
The \FMMC question was originally introduced by \textcite{BDX:fmmc-graph}, which was the first in a series of articles
\cite{BDX:fmmc-graph,BDSX:fmmc-path,SBXD:fmmc-cts,BDPX:symm-anal,BDPX:fmmc-symm}
%\cite{BDPX:symm-anal,BDPX:fmmc-symm,BDSX:fmmc-path,BDX:fmmc-graph,SBXD:fmmc-cts}
by those authors along with Parrilo and Sun.
It has subsequently been studied by
\cite{R:fmmc-dual,S:fmmc-thesis,MT:mixing-times,AS:fmmc-cliques,JJ:fmmc-partite,FK:fmmc-comparison,CA:fmmc-srw-dist}.
%\cite{AS:fmmc-cliques,CA:fmmc-srw-dist,FK:fmmc-comparison,JJ:fmmc-partite,MT:mixing-times,R:fmmc-dual,S:fmmc-thesis}.
We roughly collect these by theme.

\subsubsection*{Finding Useful Formulations}

\paragraph*{\textcite{BDX:fmmc-graph}.}

This original work introduces the \FMMC question and then primarily studies equivalent formulations.
In our view, the most important contribution of that paper,
beyond the introduction of the very interesting \FMMC question,
is their formulation of the \FMMC optimisation problem as a semi-definite program (\textit{\SDP}).
This allows the computation of an optimal solution in polynomial time \green{via standard convex optimisation techniques}.
%They also give a sub-gradient method which can be used to solve larger problems.
The \SDP leads naturally to a dual formulation, which found use in subsequent work
%\cite{R:fmmc-dual,S:fmmc-thesis,MT:mixing-times}.
\cite{MT:mixing-times,R:fmmc-dual,S:fmmc-thesis}.

\paragraph*{\textcite{R:fmmc-dual}.}

\citeauthor{R:fmmc-dual} takes the dual formulation of \cite{BDSX:fmmc-path} much further, writing the optimal spectral gap $\SG*[G][\pi]$ as a minimisation of the variance of a certain constrained graph embedding.
To quote him,
	``Informally, to obtain [a lower bound on the optimal spectral gap] we seek to embed the graph into $\mbr^{\abs V}$ so as to `spread' 	the nodes as much as possible under constraints over the distances separating nodes connected by edges.''
He re-derives the upper bound $\SG*[G][\pi] \lesssim \VC[G][\pi]$ using this formulation. This shows vertex conductance is a fundamental barrier to fast mixing. Our result shows that vertex conductance is essentially \emph{the} fundamental barrier to fast mixing.
%This shows that a particular order, namely $\VC[G][\pi]$, can be achieved by some chain.
%Contrastingly, our results are \emph{lower} bounds for $\SG*[G]$, i.e.\ showing that no chain can do better than a particular order, namely $\VC[G]^2 / \log n$.

\paragraph*{\textcite{SBXD:fmmc-cts}.}

The paper \cite{SBXD:fmmc-cts} is of a similar flavour to \cite{BDX:fmmc-graph} but in the continuous-time set-up. We discuss it in detail in \cref{sec:intro:previous-work} below.

\subsubsection*{Special Cases and Particular Examples}

\paragraph*{\textcite{BDSX:fmmc-path}.}

The special case of the path with uniform distribution is studied in the short note \cite{BDSX:fmmc-path}, as a follow-on from \cite{BDX:fmmc-graph}.
They show that the `uniform chain', i.e., the unbiased \RW with $\tfrac12$-holding at the ends, has the largest spectral gap.

\paragraph*{\textcite{BDPX:symm-anal,BDPX:fmmc-symm}.}

The \FMMC problem on graphs with rich symmetry properties is studied in \cite{BDPX:fmmc-symm}.
They are able to solve various cases analytically:
	edge-transitive graphs,
		such as the cycle;
	Cartesian products of graphs,
		such as the two-dimensional torus and the hypercube;
	distance-transitive graphs,
		such as Petersen, Hamming and Johnson graphs.
They then use  algebraic methods to study \FMMC on orbit graphs. This uses powerful representation theory arguments developed in \cite{BDPX:symm-anal}.

\paragraph*{\textcite{CA:fmmc-srw-dist}.}

Many similar scenarios, such as edge-transitive graphs, are studied in \cite{CA:fmmc-srw-dist}. The focus is on two \SDP methods. They study the degree-biased and uniform equilibria.

\paragraph*{\textcite{JJ:fmmc-partite}.}

Symmetric $K$-partite graphs and connections to sensor networks are considered in \cite{JJ:fmmc-partite}.
They compare numerically with a Metropolis--Hastings algorithm.

\paragraph*{\textcite{AS:fmmc-cliques}.}

Graphs which are overlapping unions of two cliques are studied in \cite{AS:fmmc-cliques}.
Here are are two cliques, say of sizes $r+s$ and $r+t$, respectively, and there are $s$ overlapping vertices.
The \FMMC problem is solved analytically for such graphs.

\paragraph*{\textcite{FK:fmmc-comparison}.}

A rather different approach is taken in \cite{FK:fmmc-comparison}. Their paper is focussed on comparison inequalities and majorisation of measures.
They use these to analyse the \FMMC problem.
This use of majorisation allows them to study a distance, such as \TV, separation or $\ell_2$, rather than the spectral gap, which is only a proxy for the mixing time.

\subsubsection*{Continuous-Time Set-Up}

\paragraph*{\textcite{SBXD:fmmc-cts}.}

The study of the \FMMC question in continuous-time was initiated in \diac.
The structure and goals of this paper are similar to \cite{BDX:fmmc-graph}.
The primary contribution is a convex \SDP formulation as well as some dual formulations.

Recall that a normalisation on the weights was required. Indeed, doubling all the weights doubles the spectral gap. We imposed an ``average leave-rate of $1$''.
A slightly more general `weighted average' is considered in \diac.
A number of physical interpretations of this~normalisation~are~given.

\paragraph*{\textcite{S:fmmc-thesis,MT:mixing-times}.}

The \FMMC problem is considered by \samt[\S 3.3]. It is referenced and discussed by \mtt[\S 7].
We discussed their work in detail immediately after \cref{res:intro:main:cts}.
We add a small caveat to that discussion.

\mtt[\S 7.1] claim to impose a scaling of $\qq(V) \le 1$ on their edge weightings $\qq : E \to \mbr_+$; contrast this with our imposition of $\qq(V) \le n$.
Their scaling immediately implies that the relaxation time is at least order $n$; this contradicts their theorem.
There are a couple of other points where there seem to be issues with the scalings, in particular in application of results from \diac.
It may be possible to rectify these issues,
but we have not checked~carefully.

\section{Vertex Conductance and the Optimal Spectral Gap}
\label{sec:vc}

%\green{I suggest that we add an small comment at the start explaining that we work with matching conductance throughout and finally show that $\MC \asymp \VC$}
This section is devoted to a proof of \cref{res:intro:main:vc}. In \cref{sec:vc:matching-def} we define the matching conductance of a graph, which plays a central role in the proof of \cref{res:intro:main:vc}. We also show in \cref{res:vc:prelim:vertex-matching} that matching and vertex conductance of a graph differ by at most a universal constant. In \cref{sec:vc:prelim} contains the necessary notation and preliminaries needed in the proof of \cref{res:intro:main:vc}. In \cref{sec:vc:fm} we relate the optimal spectral gap of a graph to its matching conductance. This relation is formalised in \cref{res:vc:fm:vc-cheeger}. Notice that  \cref{res:vc:fm:vc-cheeger} together with \cref{res:vc:prelim:vertex-matching} directly imply \cref{res:intro:main:vc}.

%\purple{I changed the case of the title, I hope that's ok? That way it is consistent with all the other ones I wrote. I prefer to use title case for titles and sentence case for sentences. I know it sounds a bit silly saying that, but I realise that not everyone feels the same way! Some just find sentence case easier to read so use it always}

\subsection{The Matching Conductance of a Graph}
\label{sec:vc:matching-def}

A \emph{matching} is a set of edges that do not share an endpoint. Given a set of (undirected) edges $E$ together with a weight function $w \colon E \to \mathbb{R}_{\ge 0}$, a \emph{maximum matching} for $E$ is a matching with maximum total weight (if $E$ is the edge set of an unweighted 
 graph, we assume $w$ is equal to one on $E$). We denote with $\nu(E)$ the weight of a maximum matching for $E$; that is, 
\[
\nu(E) \triangleq \max_{\text{matching} \; F \subseteq E} \sum_{e \in F} w(e).
\]

We can now define the \emph{matching conductance} of a graph.

\begin{defn}
\label{def:vc:prelim:vc}
	Let $G = (V, E)$ be a graph and $\emptyset \ne S \subset V$. The \textit{matching conductance} of $S$ is defined as
	\[
		\MC[S] \triangleq \frac{\nu(E(S,S^c))}{|S|}.
	\]
	The matching conductance of $G$ is defined as
	\[
		\MC*[G] \triangleq \min_{S \colon 0 < |S| \le |V|/2} \MC[S].
	\]
\end{defn}

The next proposition relates matching and vertex conductance.

\begin{prop}
\label{res:vc:prelim:vertex-matching}
	Let $G = (V, E)$ be a graph. Then, it holds that 
	\[
		\MC*[G] \le \VC*[G] \le 4 \MC*[G].
	\] 
\end{prop}
\begin{Proof}
The inequality $\MC*[G] \le \VC*[G]$ is obvious: for any $S \subseteq V$, $\nu(E(S,S^c))$ must be smaller than the size of the vertex boundary of $S$. Therefore, $\MC[S] \le \VC[S]$ for any $S \subseteq V$, which yields the inequality.

The proof of $ \VC*[G] \le 4 \MC*[G]$ is slightly more involved. In particular it's not true that $\MC[S] \le \VC[S]$ for any $S \subseteq V$. \cref{fig:exm:source} provides an example of a graph with a set with small matching conductance, but large vertex conductance. Nevertheless, the \emph{worst} vertex conductance of a set in a graph is related to the matching conductance of the graph. To prove this, consider $|S| \le |V|/2$ with $\MC[S] = \nu(E(S,S^c))/|S|$. We can assume $\MC[S] \le 1/4$, otherwise $\MC*[G] > 1/4 \ge \VC*[G]/4$. Let $M$ be a maximum matching for $E(S,S^c)$, that is $|M| = \nu(E(S,S^c))$, and $V(M) \subset V$ be the set of vertices adjacent to edges in $M$. Now consider the set $T = S \setminus V(M)$. We claim $T$ has small vertex conductance. To this end, consider  $\partial T$. It holds that $|\partial T| \le |V(M)|$. Indeed, any $u \in S^c \setminus V(M)$ cannot be in $\partial T$, otherwise there would exist an edge between a vertex $S \setminus V(M)$ and a vertex in $S^c \setminus V(M)$, and this would contradict the maximality of $M$. Therefore, since $|V(M)| = 2 |M|$, we have that
\[
\VC*[G] \le \frac{|\partial T|}{|T|} \le \frac{2\nu(E(S,S^c))}{|S| - 2|M|} \le \frac{2\nu(E(S,S^c))}{|S|/2} = 4 \MC[S] = 4 \MC*[G],
\]
where in the last inequality we have used the fact that $|M| = \MC[S] \cdot |S| \le |S|/4$.
\end{Proof}

\subsection{Definitions and Preliminaries}
\label{sec:vc:prelim}

Given a set of vertices $V$, together with a set of edges $E$ on $V$, a \emph{fractional matching} is a function $f \colon E \to [0,1]$ such that, for any $v \in V$, $\sum_{e \ni v} f(e) \le 1$. Moreover, the fractional matching number of $E$, denoted by $\nu^*(E)$, is the maximum total weight of a fractional matching for $E$:
\[
\nu^*(E) = \max_{f \colon E \to [0,1]} \sum_{e \in E} w(e) f(e),
\]
where the maximisation is over valid fractional matchings.

Notice that $\nu^*(E)$ is the solution of a linear program which is a convex relaxation for $\nu(E)$. As such, $\nu(E) \le \nu^*(E)$.
A useful characterisation of $\nu^*(E)$ is the following.

\begin{prop}
\label{res:vc:prelim:dualfm}
	The fractional matching number of $E$, $\nu^*(E)$, is equal to the minimum of the following linear program.
	\begin{alignat*}{2}
		&\min_{g \colon V \to \mbr_{\ge 0}} 	&\qquad 	& \sum_{u \in V} g(u)\\
		&\text{subject to} 				&		& g(u) + g(v) \ge w(u,v) \quad \forall \, \{u,v\} \in E.
	\end{alignat*}
\end{prop}
\begin{Proof}
This simply follows {from} linear programming duality.
\end{Proof}

With a slight abuse of notation we will also use $\nu(G)$ and $\nu^*(G)$ to denote the maximum (fractional) matching weight on the edge set of $G$.

Up {until} now we have considered only matchings in undirected graphs. For technical reasons, however, we will also need to consider matchings in directed graphs. Given a set $\overrightarrow{E} \subseteq V \times V$ of directed edges, a \emph{directed matching} $\overrightarrow{M} \subseteq E$ is a set of edges such that, if $(u,v),(w,z) \in \overrightarrow{M}$ and $(u,v) \ne (w,z)$, then $u \ne w$ and $v \ne z$. Alternatively, a directed matching can be seen as a subgraph where each vertex has indegree and outdegree at most one, whereas an undirected matching is a subgraph where each vertex has degree at most one. Analogously to the undirected case, we denote with $\nu(\overrightarrow{E})$ the weight of a maximum matching in $\overrightarrow{E}$. 

\begin{defn}
\label{def:vc:prelim:disc}
Given an undirected (weighted) graph $G=(V,E,w)$, an \emph{orientation} $\overrightarrow{G} = (V, \overrightarrow{E}, w)$ is a directed graph constructed by replacing each undirected edge $\{u,v\} \in E$ with a directed edge $(u,v)$ (with arbitrary orientation) having weight $w(u,v)$.
\end{defn}

The next lemma relates the maximum matching weight in an undirected graph with the maximum matching weight of its orientation.

\begin{prop}
\label{res:vc:fm:dirmatch}
	For any graph $G=(V,E,w)$ and an orientation $\overrightarrow{G}$, it holds that $\nu(\overrightarrow{G}) \le 4 \nu(G)$.
\end{prop}

\begin{Proof}
Let $M \subseteq E$ be a matching returned by the greedy algorithm for finding a maximal matching on $G$, which works as follows: let $e_1, \dots, e_m$ be an ordering of the edges of $G$ such that $w(e_1) \ge \cdots \ge w(e_m)$.
Then, greedy incrementally construct $M$ by adding $e_i$ to it, for $i=1,\dots,m$, as long as this operation maintain the property that $M$ is a matching. Denote with $e_{i_1},\dots, e_{i_{|M|}}$ the edges of $M$ ordered nondecreasingly according to their weight.

Let  $\overrightarrow{M^*}$ be a maximum matching in $\overrightarrow{G}$. We upper bound its total weight as follow. 
Let $\overrightarrow{M}_0 = \overrightarrow{M^*}$ and, for $j=1,\dots,|M|$, let $\overrightarrow{M}_j$ be the directed graph obtained from $\overrightarrow{M}_{j-1}$ by removing all edges incident to one of the endpoints of $e_{i_j}$. Since $M$ is maximal by construction,  $\overrightarrow{M}_{|M|}$ is empty. Moreover, at each iteration $j$ we remove at most four edges, since there are at most four edges in $\overrightarrow{M^*}$ that share an endpoint with $e_{i_j}$. These edges all have weight less than or equal to $w(e_{i_j})$. This is because the matching $\{e_1,\dots,e_{i_{j-1}}\}$ can be augmented by adding any one of these edges (or rather, their undirected equivalent) without breaking the property of it being a matching. But then, their weight must be less than or equal to $w(e_{i_j})$, since otherwise greedy would have chosen one of those instead of $e_{i_j}$.
Therefore, we have proved that $4\nu(G) \ge \nu(\overrightarrow{M^*})$, from which the proposition follows.
%
% Consider now an edge $e = \{u,v\} \in M$. There are at most four edges in $\overrightarrow{M^*}$ that share an endpoint with $e$. It is clear from the execution of greedy that all these edges have weight less than or equal to $w(u,v)$. The proposition follows.
\end{Proof}

\subsection{Matching Conductance and the Fastest Mixing Problem}
\label{sec:vc:fm}
The following result is due to \textcite{R:fmmc-dual} and gives a variational characterisation of $\bestgap{G}$. It follows from the fact that $\bestgap{G}$ can be expressed as the solution to a semidefinite program for which strong duality holds.
\begin{prop}
\label{res:vc:fm:bestgap}
	Let $G = (V, E)$ be a graph of $n$ vertices. Then, $\bestgap{G}$ is equal to the minimum of the following optimisation problem.
	\begin{alignat*}{2}
		&\min_{\substack{f \colon V \to \mbr^n \\ g \colon V \to \mbr_{\ge 0}}} 	
						&\qquad 	& 	\frac{\sum_{u \in V} g(u)}{\sum_{u \in V} \|f(u)\|^2} 		\\
		&\text{subject to} 	&		& 	\sum_{u \in V} f(u) = 0						\\
		&				&		& 	g(u) + g(v) \ge \|f(u) - f(v)\|^2 \quad \forall \, \{u,v\} \in E.
	\end{alignat*}
\end{prop}
\begin{rmkt}
\label{rmk:vc:fm:variational}
The variational characterisation actually given by \textcite{R:fmmc-dual} doesn't include a non-negativity constraint for the function $g$. The function $g$, however, needs to be non-negative whenever, as in our case, Markov chains on $G$ are allowed non-negative holding probabilities. More precisely, for any $u \in V$ and $P$ transition matrix of a Markov chain on $G$, if we allow $P(u,u) > 0$, then we need to require $g(u) \ge 0$.
\end{rmkt}

The variational characterisation of $\bestgap{G}$ given by \cref{res:vc:fm:bestgap} requires minimising over $n$-dimensional embeddings of the vertices in the graph. It is often more convenient to work with one-dimensional embeddings. For this reason, we introduce the following parameter.

\begin{defn}
\label{def:vc:fm:onedimgap}
	Let $G = (V, E)$ be a graph of $n$ vertices. We denote with $\onedimgap{G}$ the minimum of the following optimisation problem.
	\begin{alignat*}{2}
		&\min_{\substack{f \colon V \to \mbr \\ g \colon V \to \mbr_{\ge 0}}} 	
						&\qquad 	& 	\frac{\sum_{u \in V} g(u)}{\sum_{u \in V} f(u)^2} 		\\
		&\text{subject to} 	&		& 	\sum_{u \in V} f(u) = 0						\\
		&				&		& 	g(u) + g(v) \ge (f(u) - f(v))^2 \quad \forall \, \{u,v\} \in E.
	\end{alignat*}
\end{defn}

The following proposition shows that $\onedimgap{G}$ is a $O(\log{n})$-approximation of $\bestgap{G}$. 

\begin{prop}
\label{res:vc:fm:JL}
	Let $G = (V, E)$ be a graph. 
	It holds that
	\[
		\bestgap{G} \le \onedimgap{G} \lesssim \log{n} \cdot \bestgap{G}.
	\]
\end{prop}

%The proof of this proposition uses a standard trick (see, e.g., \textcite{MT:mixing-times}): first, we apply the Johnson--Lindenstrauss lemma \cite{JL} to show that considering only $O(\log{n})$-dimensional embeddings suffices to obtain a constant approximation for $\bestgap{G}$. Then, we transform such $O(\log{n})$-dimensional embedding into a one-dimensional embedding, but in doing so we will lose a $O(\log{n})$ factor.

%\purple{the above paragraph overflows into the margin. I added the word ``only'' so as to stop this. Also, perhaps use an enumerated list, as follows? It really separates the two steps.}

%\color{ForestGreen}
The proof of this proposition uses a standard trick (see, e.g., \textcite{MT:mixing-times}):
\begin{enumerate}
	\item 
	we apply the Johnson--Lindenstrauss lemma \cite{JL} to show that considering only $O(\log{n})$-dimensional embeddings suffices to obtain a constant approximation for $\bestgap{G}$;
	
	\item 
	we transform such $O(\log{n})$-dimensional embedding into a one-dimensional embedding,
	but in doing so we will lose a $O(\log{n})$ factor.
\end{enumerate}
%We lose an $O(\log n)$ factor in this final transformation.
%\color{black}

\begin{Proof}
The relation $\bestgap{G} \le \onedimgap{G}$ follows trivially since computing $ \onedimgap{G}$ can be seen as minimising over the same set of $n$-dimensional embeddings as for $\bestgap{G}$, with the additional constraint that only the first coordinate can be non-zero.

To prove the upper bound, let $f \colon V \to \mbr^n, g \colon V \to \mbr_{\ge 0}$ be the minimiser achieving $\bestgap{G}$ in \cref{res:vc:fm:bestgap}. Then, the Johnson--Lindenstrauss lemma ensures there exists an embedding $\widetilde{f} \colon V \to \mbr^d$ such that $d = O(\log{n})$ and, for any $u,v \in V$,
\[
\frac{1}{2} \|f(u) - f(v)\|^2 \le \|\widetilde{f}(u) - \widetilde{f}(v)\|^2 \le \|f(u) - f(v)\|^2,
\]
and
\[
\frac{1}{2} \|f(u)\|^2 \le \|\widetilde{f}(u)\|^2 \le \|f(u)\|^2.
\]

Now let $i \in \{1,\dots,d\}$ maximising  $\sum_{u,v \in V} \left(\widetilde{f}(u)_j - \widetilde{f}(v)_j\right)^2$ over $j \in \{1,\dots,d\}$, and
define $h \colon V \to \mbr$ as $h(u) = \widetilde{f}(u)_i - \frac{1}{n} \sum_{v \in V} \widetilde{f}(v)_i$ for any $u \in V$.

By construction, 
\[
\sum_{u \in V} h(u) = \sum_{u \in V} \widetilde{f}(u)_i - \sum_{u \in V} \frac{1}{n} \sum_{v \in V} \widetilde{f}(v)_i = 0.
\] 
Moreover, for any $u,v \in V$,
\[
(h(u)-h(v))^2 = (\widetilde{f}(u)_i - \widetilde{f}(v)_i)^2 \le  \|\widetilde{f}(u) - \widetilde{f}(v)\|^2 \le \|f(u) - f(v)\|^2.
\]
Therefore, $(h,g)$ is a feasible solution to the optimisation problem of \cref{def:vc:fm:onedimgap}.

Finally,
\[
\sum_{u \in V} h(u)^2 &= \frac{1}{2n} \sum_{u,v \in V} (h(u)-h(v))^2 \\
	&= \frac{1}{2n} \sum_{u,v \in V} (\widetilde{f}(u)_i - \widetilde{f}(v)_i)^2 \\
	&\ge \frac{1}{2nd}  \sum_{u,v \in V} \left\|\widetilde{f}(u) - \widetilde{f}(v)\right\|^2 \\
	&\ge \frac{1}{4nd}  \sum_{u,v \in V} \left\|{f}(u) - {f}(v)\right\|^2 \\
	&= \frac{1}{2d} \|f(u)\|^2,
\]
where we used the fact that both $h$ and $f$ are centred at zero. Therefore,
\[
\onedimgap{G} \le \frac{\sum_{u \in V} g(u)}{\sum_{u \in V} h(u)^2} 
		\le 2d \cdot \frac{\sum_{u \in V} g(u)}{\sum_{u \in V} h(u)^2} \le 2d \bestgap{G} \lesssim \log{n} \cdot \bestgap{G}.
\qedhere
\]
\end{Proof}

%\purple{you were missing the \texttt{\textbackslash qedhere} command at the end of the display, so the $\qedsymbol$ was on the wrong line. I've added this.}
% Incidentally, the $\qedsymbol$ will be on exactly the same line once we remove the equation numbering

Suppose we fix $f \colon V  \to \mbr$ that minimises the optimisation problem above. Then, by \cref{res:vc:prelim:dualfm}, $\onedimgap{G}$ can be seen as the fractional matching value of a graph $G_f$ which is constructed from $G$ by reweighing each edge $\{u,v\}$ by $(f(u) - f(v))^2$. Together with \cref{res:vc:fm:JL}, this hints towards a connection between the matching conductance of $G$ and $\bestgap{G}$. This connection is formalised in \cref{res:vc:fm:vc-cheeger}, which is the main result of this section.

\begin{thm}
\label{res:vc:fm:vc-cheeger}
	Let $G = (V, E)$ be a graph. 
	It holds that
	\[
		\MC*[G]^2 \lesssim \onedimgap{G} \lesssim \MC*[G].
	\]
	Moreover, this implies that
	\[
		\frac{\MC*[G]^2}{\log{n}} \lesssim \bestgap{G} \lesssim \MC*[G].
	\]
\end{thm}

The proof of \cref{res:vc:fm:vc-cheeger} follows the standard template of the proof of the discrete Cheeger inequality. To upper bound $\onedimgap{G}$ it suffices to construct test functions $f,g$ from a set $S$ minimising the matching conductance of $G$. The other direction is more complicated and, similarly to the case of the ``hard direction'' of the discrete Cheeger inequality, it requires using the function $f$ that minimises $\onedimgap{G}$ to construct \emph{sweep sets} and analyse the matching conductance of such sets. Analysing the matching conductance of these sweep sets, however, is not as straightforward as analysing their edge conductance as in the proof of the standard discrete Cheeger inequality.

We split the proof of \cref{res:vc:fm:vc-cheeger} in several lemmata. The first one, \cref{res:vc:fm:coarea}, relates the maximum matching of cuts in the graph to the maximum matching of a weighted directed graph appositely constructed.

\begin{lem}
\label{res:vc:fm:coarea}
	Let $G = (V, E)$ be an unweighted undirected graph and let $f \colon V \to \mbr_{\ge 0}$. 
	Let $\overrightarrow{G}_f = (V, \overrightarrow{E}_f,w_f)$ be a directed weighted graph constructed as follows: 
	\begin{enumerate}
	\item for any $u,v \in V$, $(u,v) \in \overrightarrow{E}_f$ if and only if $\{u,v\} \in E$ and $f(u) < f(v)$;
	\item for any $(u,v) \in \overrightarrow{E}_f$, $w_f(u,v) = f(v)^2 - f(u)^2$. 
	\end{enumerate}
	For any $t > 0$, define $S_t = \{u \in V \colon f(u)^2 > t\}$. Then, it holds that,
	\[
		\int_0^{\infty}\nu(E(S_t, S_t^c))  \, dt \le 2 \nu(\overrightarrow{G}_f).
	\]
\end{lem}
\begin{Proof}
For any $t \in [0,\infty)$, let $M_t \subseteq E(S_t, S_t^c)$ be a matching achieving value $\nu(E(S_t, S_t^c))$. Notice that, for any $t$, there might be several distinct maximum matching for  $E(S_t, S_t^c)$: we just pick one of them arbitrarily. We have that
\[
 \int_0^{\infty} \nu(E(S_t, S_t^c)) \, dt = \int_0^{\infty} \sum_{e \in E} \one{e \in M_t} \, dt.
\]

Notice that, for any edge $\{u,v\} \in E$ with $f(u) < f(v)$, 
\[\{ t \in [0, \infty) \colon \{u,v\} \in M_t \} \subseteq  \left(f(u)^2,f(v)^2\right].\]

Let $\overrightarrow{M}$ of $\overrightarrow{G}_f$ output by an execution of greedy on $\overrightarrow{G}_f$, which works as follows. We first order the edges $\overrightarrow{E}_f = \{e_1,\dots,e_m\}$ of $\overrightarrow{G}_f$ such that $w_f(e_1) \ge w_f(e_2) \ge \cdots \ge  w_f(e_m)$. For $i=1,\dots,m$, we incrementally construct $\overrightarrow{M}$ by including $e_i$ for $i=1,\dots,m$ as long as adding this edge doesn't break the property of $\overrightarrow{M}$ being a directed matching. This is the same algorithm as the one described in the proof of \cref{res:vc:fm:dirmatch}, with the difference that we are now constructing a directed instead of an undirected matching. %Notice that $\overrightarrow{M}$ is maximal: no edge can be added to it without breaking the property of being a directed matching.

Consider now $\{u,v\} \in E$ such that $f(u) < f(v)$ and $\{u,v\} \in M_t$ for some $t \ge 0$. Then, there must exist an edge $(u',v') \in \overrightarrow{M}$ such that 
$u = u'$ or $v = v'$ (since greedy outputs a maximal matching) and  
\[
w_f(u',v') = f(v')^2 - f(u')^2  \ge f(v)^2 - f(u)^2.
\]
The inequality above holds because otherwise greedy would have picked $(u,v)$ instead of $(u',v')$. This implies that $[f(u),f(v)) \subseteq [f(u'),f(v'))$ and $t \in [f(u'),f(v'))$. 

For any $t \ge 0$, let $h_t \colon E \to E$ be the function that maps any  $\{u,v\} \in E$ such that $f(u) < f(v)$ and $\{u,v\} \in M_t$ to an edge $(u',v')$ as above. Notice that for any edge $(u',v') \in \overrightarrow{M}$ and any $t\ge 0$, there can be at most two edges in $M_t$ that share an endpoint with $(u',v')$. Hence, $|h_t^{-1}(u',v')| \le 2$.

Therefore, we have that
\[
\sum_{\{u,v\} \in E}  \int_0^{\infty}  \one{\{u,v\} \in M_t} \, dt 
		&\le 2  \sum_{(u',v') \in \overrightarrow{M}}  \int_0^{\infty}  \one{h_t^{-1}(u',v') \cap M_t \ne \emptyset} dt \\
 		&\le 2 \sum_{(u',v') \in \overrightarrow{M}}  \left(f(v')^2 - f(u')^2\right) \\
		&\le 2 \nu(\overrightarrow{G}_f).
\]

The lemma follows by observing that
\[
\int_0^{\infty} \nu(E(S_t, S_t^c)) =  \int_0^{\infty} \sum_{e \in E}   \one{e \in M_t} \, dt = \sum_{e \in E}  \int_0^{\infty}  \one{e \in M_t} \, dt \le 2 \nu(\overrightarrow{G}_f),
\]
where we can swap the signs of integration and summation since the matchings $\{M_t : t\ge 0\}$ can be chosen so that we need to consider only at most $n-1$ different matchings (since there are at most $n-1$ different sets $S_t$), which implies that the integral can actually~be~computed~as~a~finite~sum.
\end{Proof}

%\purple{ditto the qedhere command. Also, ``sum'' was orphaned, all alone, so I added unbreaking spaces so that it's on the previous line. I did similarly with a ``that'' below}

%Let $t(e)$ be the contribution of $e \in E$ to the integral above, i.e., 
%\[
%\sum_{e \in E} t(e) = \int_0^{\infty} \nu(E(S_t, S_t^c)) \, dt.
%\]
%First notice that $t(\{u,v\}) \le w_f(u,v)$, since $\{u,v\}$ belongs to $E(S_t, S_t^c)$ only for $t \in \left(f(u)^2,f(v)^2\right]$, where we assume $f(u) < f(v)$ w.l.o.g.
%
%Now consider a matching $M$ of $\overrightarrow{G}_f$ output by an execution of greedy on $\overrightarrow{G}_f$ (see the proof of \cref{res:vc:fm:dirmatch} for a description of the algorithm). Let $(u,v) \in M$. Then, we have that
%\[
%\sum_{u' \colon e=(u',v) \in \overrightarrow{E} } t(e) \le w_f(u,v).
%\]
%This is because, since $(u,v)$ was chosen by greedy, 
%\[ \label{eq:vc:fm:dirmatch}
%w_f(u',v) = \left|f(u')^2 - f(v)^2\right| < \left|f(u)^2 - f(v)^2 \right| = w_f(u,v), 
%\]
%for any $u'$ such that  $(u',v) \in \overrightarrow{E}_f$, which implies $f(v) > f(u') > f(u)$. This in turn implies that $v$ is the arriving endpoint of edges in a matching for $E(S_t, S_t^c)$ only for $t \in \left(f(u)^2,f(v)^2\right)$. Moreover, for any fixed $t$, at most one edge incident to $v$ can belong to a matching for $E(S_t, S_t^c)$. Altogether, this proves the inequality above. Analogously,
%\[
%\sum_{v' \colon e=(u,v') \in \overrightarrow{E} } t(e) \le w_f(u,v).
%\]
%Therefore,
%\[
%\sum_{e \in E} t(e) \le 2 \sum_{(u,v) \in M} w_f(u,v) \le 2 \nu(\overrightarrow{G}_f).
%\]

The next lemma shows how to construct a set of small matching conductance given a ``good'' non-negative function $f \colon V \to \mbr_{\ge 0}$.

\begin{lem}
\label{res:vc:fm:maincheeger}
	Let $G = (V, E)$ be a graph and $f \colon V \to \mbr_{\ge 0}$. Let $\lambda$ be the minimum of the following optimisation problem.
	\begin{alignat*}{2}
		&\min_{g \colon V \to \mbr_{\ge 0}} 	
						&\qquad 	& 	\frac{\sum_{u \in V} g(u)}{\sum_{u \in V} f(u)^2} 		\\
		&\text{subject to} 	&		& 	g(u) + g(v) \ge (f(u) - f(v))^2 \quad \forall \, \{u,v\} \in E.
	\end{alignat*}
	Then, there exists a set $S \subseteq \{u \in V \colon f(u) > 0 \}$ such that $\MC[S] \le 8\sqrt{2\lambda}$.
\end{lem}
\begin{Proof}
Let $S_t = \{u \in V \colon f(u)^2 > t\}$ for $t \ge 0$. We have that
\[
\min_t \Upsilon(S_t) \le \frac{\int_0^{\infty} \nu(E(S_t, S_t^c)) \, dt }{\int_0^{\infty} |S_t| \, dt}.
\]

First notice the denominator is equal to
\[
\int_0^{\infty} |S_t| \, dt = \sum_{u\in V} f(u)^2.
\]

We now upper-bound the numerator. Let $\overrightarrow{G}_f = (V, \overrightarrow{E}_f,w_f)$ be the directed weighted graph defined in \cref{res:vc:fm:coarea}: 
\begin{enumerate}
\item for any $u,v \in V$, $(u,v) \in \overrightarrow{E}_f$ if and only if $\{u,v\} \in E$ and $f(u) < f(v)$;
\item for any $(u,v) \in \overrightarrow{E}_f$, $w_f(u,v) = f(v)^2 - f(u)^2$. 
\end{enumerate}
Notice that this is an orientation of a graph $G_f = (V, E_f,w_f)$ where each directed edge $(u,v) \in \overrightarrow{E}_f$ is replaced by $\{u,v\} \in E$. Therefore, by \cref{res:vc:fm:coarea} and \cref{res:vc:fm:dirmatch}, we have that
\[
\int_0^{\infty} \nu(E(S_t, S_t^c)) \, dt \le 2 \nu(\overrightarrow{G}_f) \le 8 \nu(G_f).
\]

We now want to relate $\nu(G_f)$ to $\lambda$. Let $M$ be a maximum matching in $G_f$. By applying Cauchy--Schwartz and the triangle inequality,
\[
\nu(G_f) &= \sum_{\{u,v\} \in M} w_f(u,v) \\
	&=  \sum_{\{u,v\} \in M} \left|f(u)^2 - f(v)^2\right| \\
	&= \sum_{\{u,v\} \in M} \left|f(u) - f(v)\right| \cdot \left|f(u) + f(v)\right| \\
	&\le \sqrt{\sum_{\{u,v\} \in M} \left(f(u) - f(v)\right)^2} 
		\sqrt{\sum_{\{u,v\} \in M} \left(f(u) + f(v)\right)^2} \\
	&\le \sqrt{\sum_{\{u,v\} \in M} (f(u) - f(v))^2} 
		\sqrt{\sum_{\{u,v\} \in M} 2\left(f(u)^2 + f(v)^2\right)} \\
	&\le \sqrt{\sum_{\{u,v\} \in M} (f(u) - f(v))^2}  \sqrt{2 \sum_{u \in V} f(u)^2},
\]
where the first inequality follows from Cauchy--Schwartz, while the second from the inequality $(a+b)^2 \le 2a^2 + 2b^2$ for any $a,b \in \mbr$. Notice that $\sum_{\{u,v\} \in M} (f(u) - f(v))^2$ can be interpreted as the weight of a matching in an undirected graph obtained from $G$ by reweighing each edge $\{u,v\}$ with weight $(f(u) - f(v))^2$. Therefore, we can apply \cref{res:vc:prelim:dualfm} and the definition of $\lambda$ to show~that
\[
\sum_{\{u,v\} \in M} (f(u) - f(v))^2 &\le \min_{g \colon V \to \mbr_{\ge 0}} 
	\left\{\sum_{u \in V} g(u) \, \middle| \, g(u) + g(v) \ge (f(u) - f(v))^2 \; \forall \, \{u,v\} \in E \right\} \\
	&\le \lambda \sum_{u \in V} f(u)^2.
\]

Putting all together, we obtain
\[
\min_t \MC[S_t] \le \frac{\int_0^{\infty} \nu(S_t) \, dt }{\int_0^{\infty} |S_t| \, dt}
	\le \frac{8\nu(G_f)}{\sum_{u\in V} f(u)^2} 
	\le \frac{8 \sqrt{2\lambda} \, \sum_{u \in V} f(u)^2}{\sum_{u\in V} f(u)^2} 
	\le 8 \sqrt{2 \lambda}.
\qedhere
\]
\end{Proof}

%\purple{qedhere}

We are now finally ready to prove \cref{res:vc:fm:vc-cheeger}.
\begin{Proof}[Proof of \cref{res:vc:fm:vc-cheeger}]
We start by proving the ``easy side'' of the Cheeger-type inequality, i.e., $\onedimgap{G} \lesssim \MC*[G]$.
Let $S \subset V$ such that $|S| \le |V|/2$ and $\MC[S] = \MC*[G]$. Define $f \colon V \to \mbr$ as 
\[
f(u) = \begin{cases}
	\frac{1}{\sqrt{2|S|}} & \text{ if } u \in S \\
	- \frac{1}{\sqrt{2|V \setminus S|}} & \text{ if } u \not\in S.
\end{cases}
\]
Let $M$ be a maximum matching for $E(S,S^c)$, i.e., $|M| = \nu(E(S,S^c))$. Define $g \colon V \to \mbr$ as
\[
g(u) = \begin{cases}
	1/|S| & \text{ if } u \in  M \\
	0 & \text{ otherwise}.
\end{cases}
\]
By construction, $f,g$ satisfy the constraints of the optimisation problem of \cref{def:vc:fm:onedimgap}. Moreover, 
$\sum_{u \in V} g(u) = 2 \nu(E(S,S^c)) / |S|$, while $\sum_{u \in V} f(u)^2 = 1$. Therefore, 
\[
\onedimgap{G} \le \frac{\sum_{u \in V} g(u)}{\sum_{u \in V} f(u)^2} = 2 \MC[S] = 2 \MC*[G]
\].

We now turn the attention to the ``harder side'', i.e., $\MC*[G]^2 \lesssim \onedimgap{G}$. 
Let $f,g \colon V \to \mbr$ minimise $\onedimgap{G}$. We cannot directly apply \cref{res:vc:fm:maincheeger} with $f$ since $f$ is not non-negative. For this reason, we define two non-negative functions $h^-,h^+ \colon V \to \mbr_{\ge 0}$ as follows. Let $c$ be the median of $f$, i.e., order the vertices in $V$ such that $f(u_1) \le f(u_2) \le \cdots \le f(u_n)$ and set $c \triangleq f(u_{\lceil n/2 \rceil})$. For any $u \in V$ define $h^-(u) \triangleq \max\{0, - (f(u)-c)\}$ and $h^+(u) \triangleq \max\{0, f(u)-c\}$. If $\sum_{u \in V} h^-(u)^2 \ge \sum_{u \in V} h^+(u)^2$, we define $h \triangleq h^-$, otherwise $h \triangleq h^+$. We now apply \cref{res:vc:fm:maincheeger} with $h$.

First notice that 
\[
\sum_{u \in V} h(u)^2 \ge \frac{1}{2} \sum_{u \in V} (f(u)-c)^2 \ge \frac{1}{2} \sum_{u \in V} f(u)^2,
\]
since $\sum_{u \in V} f(u) = 0$ because $f$ is a feasible solution to the optimisation problem of \cref{def:vc:fm:onedimgap}. Moreover, for any $\{u,v\} \in E$,
\[
(h(u) - h(v))^2 \le (f(u)-f(v))^2 \le g(u) + g(v).
\]
We can then apply \cref{res:vc:fm:maincheeger} with 
\[
\lambda \triangleq \frac{\sum_{u \in V} g(u)}{\sum_{u \in V} h(u)^2} \ge \frac{2\sum_{u \in V} g(u)}{ \sum_{u \in V} f(u)^2} = 2\onedimgap{G}
\]
Therefore, there exists $S \subseteq \{u \in V \colon h(u) > 0 \}$ such that $\MC[S] \le 8\sqrt{2\lambda}$. Moreover, by construction the support of $h$ has size at most $|V|/2$. Hence,
\[
\MC*[G] \le \MC[S] \le 8\sqrt{2\lambda} \le 16 \sqrt{\onedimgap{G}}.
\qedhere
\]
\end{Proof}

%\purple{I'm guessing you didn't know about the qedhere command? It's very helpful! :)}

\section{Almost Mixing}
\label{sec:am}

\subsection{Set-Up and Main Result}

The previous section was devoted to estimating mixing-type statistics for the \FMMC problem:
	we controlled the maximal spectral gap $\SG[\PP]$
	amongst all transition matrices $\PP$ on a given graph $G = (V, E)$ which are reversible \wrt the uniform distribution
	in terms of the vertex conductance~of~$G$.
%Typically, we consider $\pi = \unif[V]$.
The purpose of the current section is to relax the condition that the invariant distribution of $\PP$, which we denote $\pi_\PP$, is \emph{exactly} $\unif[V]$:
	we allow $\pi_\PP$ to be $\eps$-far from uniform in \TV.
We show that this can allow a \emph{significant} speed up in the mixing time versus requiring the invariant distribution to be exactly $\pi$:
	we explicitly construct a Markov chain with spectral gap order at least $\eps / (\diam G)^2$.

%Throughout, we associate a Markov chain with its transition matrix and assume that this Markov chain, or transition matrix, is ergodic: it is aperiodic and has a unique invariant distribution.

Recall that we write $\PD{V}$ for the set of positive probability distributions on a set $V$ and
\[
	\PD{\pi, \eps}
=
	\brb{ \pi' \in \PD{V} \midb \mint{x \in V} \pi'(x) / \pi(x) \ge 1 - \eps }
\Qfor
	\pi \in \PD{V}
\Qand
	\eps \in [0,1].
\]

%Let $G = (V, E)$ be a graph.
Let $\uu : E \to \mbr_+$, let $A, B, S \subseteq V$, let $x \in V$ and let $E' \subseteq E$.
We use the following notation.
\begin{gather*}
	\uu(x)
\cq
	\sumt{y \in V}
	\uu\rbr{\bra{x,y}}
\Qand
	\uu(S)
\cq
	\sumt{x \in S}
	\uu(x);
\\
	E(A, B)
\cq
	\brb{ \bra{x, y} \in E \mid x \in A, \: y \in B }
\Qand
	\uu(E')
\cq
	\sumt{e \in E'}
	\uu(e).
\end{gather*}
Define the transition matrix $\PP_\uu \in [0,1]^{V \times V}$ by
\begin{alignat*}{2}
	\PP_\uu(x,y)
&\cq
	\uu\rbr{\bra{x,y}} / \uu(x)&
&\Qfor
	x,y \in V \text{ with } x \ne y;
\\
	\PP_\uu(x,x)
&\cq
	1 - \sumt{y \in V \setminus \bra{x}} \PP_\uu(x,y)&
&\Qfor
	x \in V.
\end{alignat*}
Abbreviate the spectral gap as $\SG[\uu] \cq \SG[\PP_\uu]$.
Define the probability measure $\pi_\uu : V \to [0,1]$
%and flow $Q_\uu : E \to [0,1]$
by
\[
	\pi_\uu(x)
\cq
	\uu(x) / \uu(V)
\Qfor
	x \in V.
%\Qand
%	Q_\uu(e)
%\cq
%	\uu(e) / \uu(V)
%\Qfor
%	e \in E.
\]
$\PP_\uu$ is the transition matrix of the \RW on the weighted graph and $\pi_\uu$ is its invariant distribution.
It is the unique invariant distribution if the edge set $\bra{ e \in E \mid \uu(e) \ne 0 }$ is connected.
%; $Q_\uu$ is sometimes known as its \textit{flow}.
	%
%\end{notationt}

The following theorem is a refinement of \cref{res:intro:main:am}.

\begin{thm}[`Almost Mixing']
\label{res:am:res:main}
Let $G = (V, E)$ be a graph and let $\pi \in \PD{V}$.
There exists an edge weighting $\ww_1 : E \to \mbr_+$,
depending only on $G$ and $\pi$,
with unit total weight, i.e.\ $\ww_1(V) = 1$, and the following property.
Let $\eps \in (0,1)$.
Let $\PP$ be a transition matrix on $G$ which is reversible \wrt $\pi$; it need not be irreducible.
Let $\ww_0 : E \to \mbr_+$ be the unique edge weighting of $G$ with
	$\PP_\ww = \PP$ and $\ww_0(V) = 1$.
Define the transition matrix $\QQ$ via the superposition weighting $\ww \cq \ww_0 + \eps \ww_1$:
	\[
		\QQ \cq \PP_\ww
	\Qwhere
		\uu(e)
	\cq
		\ww_0(e) + \eps \ww_1(e)
	\Qfor
		e \in E.
	\]
Then
\(
	\mint{x \in V}
	\pi_\ww(x) / \pi(x)
\ge
	1 - \eps.
\)
Let $\QQ' \cq \tfrac12 (I + \QQ)$ denote its lazification.
Then
\[
	\SG[\QQ']
\ge
	\tfrac1{48} \eps (\diam G)^{-2}
\Qand
	\MIX[\QQ'][\xi]
\le
	24 \eps^{-1} (\diam G)^2 \log\rbr{ \xi^{-2} \pimin^{-1} }
\Qforall
	\xi \in (0,1).
\]
Immediate consequences of these properties are
\(
	\tv{ \pi_\ww - \pi }
\le
	\eps,
\)
\[
	\QQ(e) \ge (1 - \eps) \PP(e)
\Qforall
	e \in E
\Qand
	\MIX[\QQ']
=
	\MIX[\QQ'][\tfrac14]
\le
	24 \eps^{-1} (\diam G)^2 \log\rbr{ 16 \pimin^{-1} }.
\]
\end{thm}

%\begin{thm}%[Almost Mixing in Diameter Squared]
%\label{res:am:res:main}
%	%
%Let $G = (V, E)$ be a graph and $\pi \in \PD{V}$.
%Let $\eps \in (0,1)$.
%There exists a lazy, reversible, ergodic Markov chain $\PP$ on $G$
%such that
%%	with the following properties:
%\begin{gather*}
%	\mint{x \in V}
%	\pi_\PP(x) / \pi(x)
%\ge
%	1 - \eps,
%\quad
%	\SG[\PP]
%\ge
%	\tfrac1{48} \eps (\diam G)^{-2}
%\quad
%	\text{and}
%\\
%	\tau_\PP(\xi)
%\le
%	24 \eps^{-1} (\diam G)^2 \log\rbr{ \xi^{-2} \pimin^{-1} }
%\Qforall
%	\xi \in (0,1).
%\end{gather*}
%Two immediate consequences of these properties are
%\[
%	\tv{ \pi_\PP - \pi }
%\le
%	\eps
%\Qand
%	\MIX[\PP]
%=
%	\MIX[\PP][\tfrac14]
%\le
%	24 \eps^{-1} (\diam G)^2 \log\rbr{ 16 \pimin^{-1} }.
%\]
%	%
%\end{thm}

We remark briefly on the `independence' of the perturbation by $\eps \ww_1$.

\begin{rmkt}[Independence of Perturbation]
The weighting $\ww = \ww_0 + \eps \ww_1$ can be seen as a perturbation of $\ww_0$ by $\eps \ww_1$, since we are most interested in the case where $\eps$ is very small---indeed, we want the new equilibrium distribution to be very close to $\pi$.
We emphasise that the perturbation weighting $\ww_1$ does not depend on the base weighting $\ww_0$; rather, $\ww_1$ is a function only $G$ and $\pi$.
\end{rmkt}

We fix a graph $G = (V, E)$ and a probability measure $\pi \in \PD{V}$ throughout this section.
We do not always repeat these in statements below.
Also, we write $n \cq \abs V$.

We start by proving a slightly weaker statement.
Assume that $\ww_0$ corresponds to unit-self~loops:
\[
	\ww_0\rbr{\bra{x,y}} = 0
\Qforall
	x, y \in V
\Qwith
	x \ne y
\Qand
	\ww_0\rbr{\bra{v}} = \pi(v)
\Qforall
	v \in V.
\]
The corresponding transition matrix $\PP_{\ww_0}$ is diagonal and thus reversible \wrt any measure.
We then extend the argument to handle arbitrary initial weightings $\ww_0$ in the \cref{sec:am:perturb}.

\subsection{Outline and Proof Given Later Results}
\label{sec:am:outline}

We start by giving a very brief outline with cross-references to the results proved in the following subsections.
We then flesh out this outline, giving a more detailed description.
%Finally, we describe precisely how to deduce the main theorem from results cited in the brief outline.

\begin{Proof}[Outline of Proof: Very Brief]
\qedtriangle
The proof has four key steps.
\begin{enumerate}
	\item 
	We construct a weighted spanning tree;
	see \cref{def:am:pf:tree}.
	
	\item 
	We control the difference between the invariant distribution of the \RW on this weighted tree and the target distribution $\pi$;
	see \cref{res:am:pf:inv-dist}.
	
	\item 
	We estimate the conductance of this weighted tree;
	see \cref{res:am:pf:conductance}.
	
	\item 
	We relate its spectral gap and conductance using canonical paths;
	see \cref{res:am:pf:gap}.
\qedhere
\end{enumerate}
\end{Proof}

\begin{Proof}[Outline of Proof: More Detailed]
We now flesh out the above details somewhat for $\pi = \unif[V]$.
\begin{itemize}
	\item 
	Let $T = (V, F)$ be a \BFS spanning tree of $G$, rooted at $v_\star$.
	We choose a weighting $\ww_\star : F \to (0, \infty)$ such that the weights increase towards the root $v_\star$ in such a way that $\ww_\star$ has edge conductance order $1$.
	We then rescale the $\ww_\star$ to get $\widetilde \ww_\star$ with total weight $\widetilde \ww_\star(V) = \eps n$.
	
	\item 
	Define $\ww : E \to \mbr_+$ by adding unit-weight self-loops to $\widetilde \ww_\star$.
	The total weight $\ww(V) = (1 + \eps) n$ and $\ww(x) \ge 1$ for all $x$.
	Thus the invariant distribution $\pi'$ satisfies
	\(
		\min_{x \in V}
		\pi'(x)
	\ge
	%	(1 + \eps)^{-1} / n
	%\ge
		(1 - \eps) / n.
	\)
	
	\item 
	It remains to analyse the edge conductance of $\ww$, which is intimately related to the original total weight $\ww_\star(V)$.
	We can choose the weights such that $\ww_\star(V) \asymp n \diam G$.
	We then apply a Cheeger-type inequality to deduce a spectral gap lower bound of order $\eps / (\diam G)^2$.
\end{itemize}

We now describe how to choose the weighting $\ww_\star$.
Let $T_x \subseteq T$ denote the subtree rooted at $x$ and consisting of all descendants of $x$.
We choose the weight
\(
	\ww_\star\rbr{\bra{x,\prnt(x)}}
\cq
	\abs{T_x},
\)
where $\prnt(x)$ is the (unique) parent of $x$, for $x \ne o$.
This way, the conductance of a subtree $T_x$ \emph{in the weighted tree $(T, \ww_\star)$} is precisely $1$.
We emphasise that this is in the weighted tree $(T, \ww_\star)$.
We need to rescale $\ww_\star$ and combine it with the unit-weight self-loops to get an approximately uniform weighting.

It turns out that
\(
	\ww_\star(F)
\asymp
	n \diam T
\asymp
	n \diam G.
\)
This then gives rise to a final conductance $\EC* \asymp \eps / \diam G$.
The standard Cheeger inequality then gives $\SG \gtrsim \eps^2 / (\diam G)^2$.
We improve this to $\SG \gtrsim \EC^* / (\diam T)$ by applying the canonical paths method, using the fact that $T$ is a tree.

The proof for general $\pi$ is very similar.
One gives the self-loop at $x$ weight $\pi(x)$ and defines
\(
	\ww_\star\rbr{ \bra{ x, \prnt(x) } }
\cq
	\pi(T_x).
\)
This is the natural extension.
The same arguments follow through.
%Further, one \emph{does not} lose any factors such as $\pimax/\pimin$.
	%
\end{Proof}

\subsection{Preliminaries}
\label{sec:am:prelim}

We introduce some preliminary material which is used throughout the proof, as well as in \cref{sec:cts}.
%and \cref{sec:tinhom}.
The majority of it will be familiar to a reader well-versed in \RWs and mixing time analysis.

%Let $T = (V, F)$ be a graph and $\uu : F \to \mbr_+$ be a weighting.
%We use the following notation.

%\begin{notationt}
%	%
%Let $A, B, S \subseteq V$, let $x \in V$ and let $F' \subseteq F$.
%\begin{gather*}
%	\uu(x)
%\cq
%	\sumt{y \in V}
%	\uu\rbr{\bra{x,y}}
%\Qand
%	\uu(S)
%\cq
%	\sumt{x \in S}
%	\uu(x);
%\\
%	F(A, B)
%\cq
%	\brb{ \bra{x, y} \in F \mid x \in A, \: y \in B }
%\Qand
%	\uu(F')
%\cq
%	\sumt{e \in F'}
%	\uu(e).
%\end{gather*}
%Define the transition matrix $\PP_\uu \in [0,1]^{V \times V}$ by
%\begin{alignat*}{2}
%	\PP_\uu(x,y)
%&\cq
%	\uu\rbr{\bra{x,y}} / \uu(x)&
%&\Qfor
%	x,y \in V \text{ with } x \ne y;
%\\
%	\PP_\uu(x,x)
%&\cq
%	1 - \sumt{y \in V \setminus \bra{x}} \PP_\uu(x,y)&
%&\Qfor
%	x \in V.
%\end{alignat*}
%Define the probability measure $\pi_\uu : V \to [0,1]$ and flow $Q_\uu : F \to [0,1]$ by
%\[
%	\pi_\uu(x)
%\cq
%	\uu(x) / \uu(V)
%\Qfor
%	x \in V
%\Qand
%	Q_\uu(e)
%\cq
%	\uu(e) / \uu(V)
%\Qfor
%	e \in F.
%\]
%$\PP_\uu$ is the transition matrix of the \RW on the weighted graph and $\pi_\uu$ is its invariant distribution; $Q_\uu$ is sometimes known as its \textit{flow}.
%Let $\SG[\uu] \cq \SG[\PP_\uu]$ denote its spectral gap.
%	%
%\end{notationt}

First, we generalise the notation of \textit{edge conductance} of the graph, herein abbreviated \textit{conductance}.
We introduced this in \cref{def:intro:main:vc:edge} for \RWs on unweighted graphs, i.e.\ unit weights on all edges.
Reversible \RWs correspond to a general weighting $\uu : E \to \mbr_+$,
as in the above notation.
The (edge) conductance of a reversible \RW is the (edge) conductance of that weighted graph.%

\begin{defn}[Edge Conductance]
\label{def:am:prelim:conductance}
Let $G = (V, F)$ be a graph and $\uu : F \to \mbr_+$ be a weighting.
The \textit{conductance} $\ECRW[\uu][S]$ of a set $S \subseteq V$ with $\pi_\uu(S) > 0$ \wrt $\uu$ is defined to be
\[
	\ECRW[\uu][S]
\cq
	\frac
		{\pi_\uu\rbb{ F(S, S^c) }}
		{\pi_\uu(S)}
=
	\frac
		{\uu\rbb{ F(S, S^c) }}
		{\uu(S)}.
\]
The \textit{conductance} $\ECRW*[\uu]$ of $\uu$ is defined to be
\[
	\ECRW*[\uu]
\cq
	\MIN{S \subseteq V : 0 < \pi_\uu(S) \le 1/2}
	\ECRW[\uu][S].
\]

The \textit{adjusted conductance} $\aECRW[\uu][S]$ of a set $S \subseteq V$ with $0 < \pi_\uu(S) < 1$ is defined similarly:
\[
	\aECRW[\uu][S]
\cq
	\frac
		{\pi_\uu\rbb{ F(S, S^c) }}
		{\pi_\uu(S) \pi_\uu(S^c)}
\Qand
	\aECRW*[\uu]
\cq
	\MIN{S \subseteq V : 0 < \pi_\uu(S) < 1} \:
	\aECRW[\uu][S].
\]
\end{defn}

\begin{subtheorem}{thm}
\label{rmk:am:prelim:conductance}

\begin{rmkt}[Conductance: Original and Adjusted Relations]
\label{rmk:am:prelim:conductance:adjusted}
The definition of the adjusted conductance does not need the restriction $\pi_\uu(S) \le \tfrac12$ as $\aECRW[\uu][S]$ is invariant under complementation:%
\[
	\aECRW[\uu][S] = \aECRW[\uu][S^c]
\Qforall
	S \subseteq V \text{ with } 0 < \pi_\uu(V) < 1.
\]
The following inequalities between $\ECRW[\uu]$ and $\aECRW[\uu]$ are immediate:
\begin{alignat*}{2}
	&\mathrel{\phantom{\le}} \ECRW[\uu][S] \le \aECRW[\uu][S]&
&\Qforall
	S \subseteq V \text{ with } 0 < \pi_\uu(S) < 1;
\\
	\tfrac12 \aECRW[\uu][S] &\le \ECRW[\uu][S]&
&\Qforall
	S \subseteq V \text{ with } 0 < \pi_\uu(S) \le \tfrac12.
\\
	\tfrac12 \aECRW*[\uu]
&\le
	\ECRW*[\uu]
\le
	\aECRW*[\uu]&
&\quad
	\text{whatever the graph}.
%\nonumber
\tag*{\qedhere}
\end{alignat*}
\end{rmkt}

\begin{rmkt}[Conductance: Connectivity Assumption]
\label{rmk:am:prelim:conductance:connected}
We may assume that $S$ induces a connected subset $T[S]$ when analysing $\ECRW*[\uu]$.
Indeed, if $S = A \mathrel{\dot \cup} B$ with $F(A, B) = \emptyset$, then
\[
	\ECRW[\uu][S]
=
	\frac
		{\uu\rbb{ F(A, A^c) } + \uu\rbb{ F(B, B^c) }}
		{\uu(A) + \uu(B)}
\ge
	\min\brb{ \ECRW[\uu][A], \: \ECRW[\uu][B] },
%\label{eq:am:prelim:conductance:conn}
%\nt
\]
using the fact that
\(
	F(S, S^c)
=
	F(A, A^c) \mathrel{\dot \cup} F(B, B^c)
\)
and that
\[
	\frac{a + b}{a' + b'}
\ge
	\min\brbb{ \frac{a}{a'}, \: \frac{b}{b'} }
\Qforall
	a,a',b,b' > 0.
\qedhere
\]
\end{rmkt}

\end{subtheorem}

Next, we introduce the \textit{canonical paths method} and use it to relate the spectral gap to the conductance in trees.
A proof of \cref{res:am:prelims:canonical:gen} can be found in \cite[Theorem~5]{S:canonical-paths}.

\begin{defn}[Paths]
\label{def:am:prelim:paths}
Let $G = (V, F)$ be a graph.
$\Gamma : \bra{0, ..., L} \to V$ is an \textit{$F$-path} from $x$ to $y$
%, with $x,y \in V$,
if
\[
	\Gamma(0) = x,
\quad
	\Gamma(L) = y
\Qand
	\bra{ \Gamma(\ell-1), \: \Gamma(\ell) } \in F
\Qforall
	\ell \in [1, L].
\]
The \textit{length} of a path $\Gamma : \bra{0, ..., L} \to V$ is defined to be
\(
	\abs \Gamma
\cq
	L.
\)
\end{defn}

\begin{prop}[Canonical Paths: General]
\label{res:am:prelims:canonical:gen}
	Let $G = (V, E)$ be a connected graph and $\uu : E \to \mbr_+$ be a weighting.
	Let $\Gamma_{x,y}$ be an arbitrary $F$-path from $x$ to $y$ for $x,y \in V$.
	The spectral gap $\SG[\uu]$ of the \RW on $(G, \uu)$ satisfies
	\[
		\SG[\uu]
	\ge
		\MIN{e \in E}
		\brb{
			\rbb{ \uu(e) / \uu(V) }
		\big/
			\sumt{x,y \in V}
			\one{e \in \Gamma_{x,y}}
			\pi_\uu(x) \pi_\uu(y) \abs{\Gamma_{x,y}}
		}.
	\]
\end{prop}

\begin{cor}[Canonical Paths: Trees]
\label{res:am:prelims:canonical:trees}
	Let $T = (V, F)$ be a connected tree and $\uu : F \to (0, \infty)$ be a weighting.
	The spectral gap $\SG[\uu]$ of the \RW on $(T, \uu)$ satisfies
	\[
		\SG[\uu]
	\ge
		\aECRW*[\uu] / \diam T.
	\]
\end{cor}

\begin{Proof}
Let $\Gamma_{x,y}$ be the shortest path between $x$ and $y$ for all $x,y \in V$.
Then $\abs{\Gamma_{x,y}} \le \diam T$ for all $x,y \in V$.
Removing the edge $e = \bra{e_-, e_+} \in F$ disconnects the graph, leaving two components, with $e_- \in V$ in one component and $e_+ \in V$ in the other.
Denote the component containing $e_\pm$ by $T^e_\pm$.
The canonical paths method (\cref{res:am:prelims:canonical:gen}) then implies that
\[
	\SG[\uu]
\ge
	\frac1{\diam T}
\cdot
	\MIN{e \in F}
	\frac
		{\uu(e) / \uu(V)}
		{\pi_\uu(T^e_+) \pi_\uu(T^e_-)}
=
	\frac
		{\aECRW[\uu][T^e_+]}
		{\diam T}
\ge
	\frac
		{\aECRW*[\uu]}
		{\diam T}.
\]
This uses the fact that
\(
	T^e_+ \mathrel{\dot \cup} T^e_-
=
	V
\)
and
\(
	F(T_+^e, T_-^e)
=
	e
\)
for all $e \in F$.
\end{Proof}

\begin{rmkt}
	A more general statement is proved by \textcite[Theorem~1]{M:conductance-spectralgap} (article in French).
	He does not require the graph $T$ to be a tree at the cost of replacing $\diam T$ in the denominator of the bound by the longest path in $T$;
	eg, if $T$ has a Hamiltonian path, then the denominator becomes $n-1$.
	He gives two proofs, one of which uses a canonical paths style argument. 
\end{rmkt}

Finally, we introduce some notation for trees and prove a counting lemma.
This lemma, innocuous as it may appear, is fundamental to multiple calculations.
The notation and definitions above were for any graph $T = (V, F)$.
Assume that $T$ is a tree for the rest of this preliminary section.

\begin{defn}[Tree Notation]
\label{def:am:prelim:tree-not}
Let $T = (V, F)$ be a tree rooted at $o$.
\begin{itemize}
	\item 
	Let $\anc(z)$ denote the unique shortest path from $z$ to the root $o$, including both $z$ and $o$.
	
	\item 
	Let $V_y \cq \bra{ z \in V \mid y \in \anc(z) }$ and $T_y \cq T[V_y]$ denote the subtree rooted at $x$.
	
	\item 
	Let $\prnt(x)$ denote the parent of $x \ne o$,
	ie the unique neighbour $y$ of $x$ satisfying $y \in \anc(x)$.
\end{itemize}
\end{defn}

\begin{lem}[Counting Weighted Subtrees]
\label{res:am:prelim:counting-trees}
	For
		all measures $\mu$ on $V$
	and
		all $x \in V$,
	we have
	\[
		\sumt{y \in T_x \setminus \bra{x}}
		\mu(T_y)
	\le
		\mu(T_x) \diam T.
	\]
\end{lem}

\begin{Proof}
	$z \in T_y$ if and only if $y \in \anc(z)$;
	if $y \in T_x \cap \anc(z)$, then $z \in T_x$.
Thus
\[
	\sumt{y \in T_x \setminus \bra{x}}
	\mu(T_y)
&
=
	\sumt{y \in T_x \setminus \bra{x}}
	\sumt{z \in T_y}
	\mu(z)
\\&
=
	\sumt{y \in V}
	\sumt{z \in V}
	\mu(z) \one{y \in T_x \setminus \bra{x}} \one{z \in T_y}
\\&
=
	\sumt{z \in V}
	\mu(z)
	\sumt{y \in V}
	\one{y \in \anc(z) \cap T_x \setminus \bra{x}} \one{z \in T_x}
\\&
\le
	\sumt{z \in T_x}
	\mu(z)
	\depth T_x
=
	\mu(T_x) \depth T_x
\le
	\mu(T_x) \diam T.
\qedhere
\]
\end{Proof}

\subsection{Construction via a Weighted Spanning Tree}
\label{sec:am:pf}

First, we define a weighted spanning tree $(T, \ww)$.

\begin{defn}[Weighted Spanning Tree]
\label{def:am:pf:tree}
	Let $o \in V$ and let $T = (V, F)$ be a \BFS tree rooted at $o$.
	Supplement $F$ with a self-loop at each vertex of $V$.
	Define the following weightings $\ww_{0/1} : F \to \mbr_+$:
	\begin{alignat*}{2}
		\ww_0\rbr{\bra{x,y}}
	&\cq
		\pi(x)\one{x = y}&
	&\Qfor
		x,y \in V;
	\\
		\ww_1\rbr{\bra{x, \prnt(x)}}
	&\cq
		\pi(T_x)&
	&\Qfor
		x \in V \setminus \bra{o}.
	\end{alignat*}
	Define the weighting $\ww : F \to \mbr_+$ via a linear combination:
	\[
		\ww(\cdot)
	\cq
		\ww_0 + \eta \ww_1
	\Qwhere
		\eta \cq \tfrac12 \eps / \diam T.
	\]
	
	$T$ is a \BFS tree, so $\DIAM \cq \diam T$ satisfies $\DIAM \le 2 \diam G$.
\end{defn}

%\begin{itemize}
%	\item 
%	$\ww_0$ is given by weight-$1$ self-loops everywhere:
%	\[
%		\ww_0\rbr{\bra{x, y}}
%	\cq
%		\one{x = y}
%	\Qfor
%		x,y \in V.
%	\]
%	This gives the underlying approximate uniformity to $\ww = \ww_0 + \eta \ww_1$
%	
%	\item 
%	$\ww_1$ is given by the size of subtrees:
%	\[
%		\ww_1\rbr{\bra{x, \prnt(x)}}
%	\cq
%		\abs{T_x}
%	\Qfor
%		x \in V \setminus \bra{o}.
%	\]
%	Typically, vertices $x$ satisfy $\eta \ww_1(x) \ll \eps \ww_0(x)$; e.g.\, $\eta \ww_1(x) = \eta \ll 1 = \ww_0(x)$ if $x \in T$ is a leaf.
%	It is used to amplify the weights through bottlenecks, thus removing them.
%	Some, but few, vertices $x$ have very large $\ww_1(x)$; e.g.\, $\ww_1(o) = n - 1$.
%\end{itemize}

The distribution $\pi_\ww$ induced by $\ww$ is close to $\pi$ in the following sense.

\begin{lem}[$\pi_\ww$ Close to $\pi$]
\label{res:am:pf:inv-dist}
	The weighted tree $(T, \ww)$ and its induced distribution $\pi_\ww$ satisfy
	\[
		\ww(V)
	\le
		1 + \eps
	\Qand
		\MIN{x \in V}
		\pi_\ww(x) / \pi(x)
	\ge
		1 / (1 + \eps)
	\ge
		1 - \eps.
	\]
\end{lem}

\begin{Proof}
First,
the subtree counting lemma (\cref{res:am:prelim:counting-trees}) implies that
\[
	\ww_1(V)
=
	2 \ww_1(F)
=
	2
	\sumt{e \in F}
	\ww_1(e)
&
=
	2
	\sumt{x \in V \setminus \bra{o}}
	\ww_1\rbr{\bra{x, \prnt(x)}}
\\&
=
	2
	\sumt{x \in T_0 \setminus \bra{o}}
	\pi(T_x)
\le
	2 \pi(T_o) \DIAM
=
	2 \DIAM.
\]
Trivially,
\(
	\ww_0(V)
=
	\sumt{x \in V}
	\pi(x)
=
	1.
\)
Thus
\[
	\ww(V)
=
	\ww_0(V) + \eta \ww_1(V)
\le
	1 + (\tfrac12 \eps / \DIAM) \cdot (2 \DIAM)
=
	1 + \eps.
\]

Second,
\(
	\ww(x)
\ge
	\ww_0(x)
=
	\pi(x)
\)
and thus
\(
	\pi_\ww(x)
=
	\ww(x) / \ww(V)
\ge
	\pi(x) / (1 + \eps)
\)
for all $x \in V$.
\end{Proof}

Next, we control the conductance of this weighted spanning tree.

\begin{prop}[Conductance]
\label{res:am:pf:conductance}
	Let $(T, \ww)$ be the weighted spanning tree from \cref{def:am:pf:tree}.
	The conductance $\ECRW*[\ww]$ of the \RW on $(T, \ww)$ satisfies
	\[
		\ECRW*[\ww]
	\ge
		\tfrac16 \eps / \DIAM.
	\]
\end{prop}

\begin{Proof}
First, suppose that $o \notin S$.
Choose $x \in S$ with $\dist(x, o)$ minimal.
We may assume that $T[S]$ is connected, by
%\cref{eq:am:prelim:conductance:conn} in
\cref{rmk:am:prelim:conductance:connected}.
These together imply that $S \subseteq V_x$ and
\[
	\brb{ x, \prnt(x) }
=
	F(T_x, T_x^c)
\subseteq
	F(S, S^c).
\]
This implies that
\[
	\ECRW[\ww][S]
=
	\frac
		{\ww\rbb{ F(S, S^c) }}
		{\ww(S)}
\ge
	\frac
		{\ww\rbb{ F(T_x, T_x^c) }}
		{\ww(T_x)}
=
	\ECRW[\ww][T_x].
\]
The definition of $\ww$ gives.
\[
	\ECRW[\ww][T_x]
=
	\frac
		{\ww\rbr{ \bra{x, \prnt(x)} }}
		{\ww(T_x)}
=
	\frac
		{\eta \pi(T_x)}
		{\pi(T_x) + \sumt{y \in T_x} \eta \pi(T_y)}.
\]
The tree counting lemma (\cref{res:am:prelim:counting-trees}), then implies that
\[
	\ECRW[\ww][T_x]
\ge
	\frac
		{\eta \pi(T_x)}
		{\pi(T_x) + \eta \DIAM \pi(T_x)}
=
	\frac{\eta}{1 + \eta \DIAM}
%=
%	\frac{\frac12 \eps / \DIAM}{1 + \frac12 \eps}
=
	\frac{\eps}{(2 + \eps) \DIAM}
\ge
	\frac{\eps}{3 \DIAM},
\]
recalling that $\eta = \tfrac12 \eps / \DIAM$.
We have thus shown that
\[
	\ECRW[\ww][S]
\ge
	\frac{\eps}{3 \DIAM}
\Qforall
	S \subseteq V \text{ with } o \notin S \ne \emptyset.
%\label{eq:am:pf:conductance:o-notin-S}
%\nt
\]
Importantly, this inequality does not require $\pi_\ww(S) \le \tfrac12$.

Next, suppose that $o \in S$ and $\pi_\ww(S) \le \tfrac12$.
The relations
%\cref{eq:am:prelim:conductance:adj}
of \cref{rmk:am:prelim:conductance:adjusted} imply that
\[
	\ECRW[\ww][S]
\ge
	\tfrac12 \aECRW[\ww][S]
=
	\tfrac12 \aECRW[\ww][S^c]
\ge
	\tfrac12 \ECRW[\ww][S^c].
\]
But now $S' \cq S^c$ satisfies $o \notin S'$ and $S' \ne \emptyset$.
Thus the previous case
%\cref{eq:am:pf:conductance:o-notin-S}
implies that
\[
	\ECRW[\ww][S^c]
=
	\ECRW[\ww][S']
\ge
	\tfrac13 \eps / \DIAM.
%\label{eq:am:pf:conductance:o-in-S}
%\nt
\]
Importantly, this case did not require $\pi_\ww(S') \le \tfrac12$.
Combining
%\cref{eq:am:pf:conductance:o-notin-S,eq:am:pf:conductance:o-in-S}
gives
\[
	\ECRW[\ww][S]
\ge
	\tfrac16 \eps / \DIAM
\Qforall
	S \subseteq V \text{ with } 0 < \pi_\ww(S) \le \tfrac12.
\qedhere
\]
\end{Proof}

%We
%	constructed a weighted spanning tree $(T, \ww)$ in \cref{def:am:pf:tree}
%and
%	controlled its conductance $\ECRW*[\ww]$ in \cref{res:am:pf:conductance}.

Finally, we apply the canonical paths method for trees (\cref{res:am:prelims:canonical:trees}) to deduce a bound on the spectral gap for $(T, \ww)$.

\begin{cor}[Spectral Gap]
\label{res:am:pf:gap}
	Let $(T, \ww)$ be the weighted spanning tree from \cref{def:am:pf:tree}.
	The spectral gap $\SG[\ww]$ of the \RW on $(T, \ww)$ satisfies
	\[
		\SG[\ww]
	\ge
		\tfrac16 \eps / \DIAM^2
	\ge
		\tfrac1{24} \eps / (\diam G)^2.
	\]
\end{cor}

\begin{Proof}
This is an immediate consequence of
	the canonical paths method for trees (\cref{res:am:prelims:canonical:trees})
and
	\cref{res:am:pf:conductance},
along with the relations of \cref{rmk:am:prelim:conductance:adjusted}.
Also, $\DIAM \le 2 \diam G$.
\end{Proof}

We have now almost proved the main result. We just need to make sure the chain is lazy and convert the spectral gap result into a mixing time result.

\begin{Proof}[Proof of \cref{res:am:res:main} when $\ww_0 = \pi$]
The Markov chain constructed is a \RW on a weighted \BFS tree. It is defined in \cref{def:am:pf:tree}.
Denote the invariant distribution of the \RW on this weighted tree by $\pi'$.
\cref{res:am:pf:inv-dist} establishes the claim on the invariant distribution.

The spectral gap bound is proved via \cref{res:am:pf:conductance,res:am:pf:gap}.
Precisely, \cref{res:am:pf:gap} defines a reversible chain $\QQ$ satisfying $\SG[\QQ] \ge \tfrac1{24} \eps / (\diam G)^2$.

The mixing time bound will follow from the spectral gap bound via a standard mixing time--spectral gap relation.
%or \cite[Theorem~12.4]{LPW:markov-mixing}
To apply this relation, we first pass from $\QQ$ to its lazy version $\QQ' \cq \tfrac12(I + \QQ')$.
This ensures that the \emph{spectral gap} and \emph{absolute spectral gap} agree.
$\QQ$ and $\QQ'$ have the same invariant distribution and that
\(
	\SG[\QQ']
=
	\tfrac12 \SG[\QQ].
\)
A simple calculation establishes the mixing time claim using the spectral--mixing relation; see \cite[Lemma~4.23]{AF:book} for details of this relation.
\end{Proof}

It remains to handle the case of general $\ww_0$, i.e.\ where $\ww_0$ is any unit edge weighting with $\pi$ as its induced invariant distribution.
This is done in the next subsection.

\subsection{Perturbation to Arbitrary Base Chain}
\label{sec:am:perturb}

The analysis up to this point has shown the \emph{existence} of a fast `almost mixing' chain.
Precisely, we defined a weighted graph by
	constructing an appropriately weighted \BFS tree
and
	supplementing it with $\pi$-weighted self-loops.
We can think of the self-loops as a `base' weighting which is reversible \wrt $\pi$.
We denoted the `base' weighting $\ww_0$ and the `tree' weighting $\ww_1$; recall \cref{def:am:pf:tree}.

We now explain how to extend this to an arbitrary `base' weighting $\ww_0$.
The analysis is extremely similar to that of the self-loops case above:
	we simply take an arbitrary base weighting $\ww_0$ and superimpose on it the same weighted \BFS tree.
Some small adjustments are needed, but not many.

\begingroup

\renewcommand{\qedsymbol}{\ensuremath{\triangle}}

\medskip

%\begin{Proof}[Notation]
%	%
Let $\ww_0 : E \to \mbr_+$ be an arbitrary unit edge weighting of $E$.
Define
\(
	E'
\cq
	\bra{ e \in E \mid \ww_0(e) \ne 0 },
\)
the edge set of the graph induced by $\ww_0$.
%The corresponding weighted \RW is on $(V, E')$.
$\pi$ is the invariant distribution of the $\ww_0$-weighted \RW.
%	%
%\end{Proof}

\begin{Proof}[Construction]
We define the weighted tree $(T, \ww_1)$ exactly as before in \cref{def:am:pf:tree}:
	$T = (V, F)$ is an arbitrary \BFS tree
and
	$\ww_1\rbr{\bra{x, \prnt(x)}} = \pi(T_x)$ for $x \in V \setminus \bra{o}$;
set
\(
	\ww \cq \ww_0 + \eps \ww_1.
\)
The proof of \cref{res:am:pf:inv-dist} is unchanged,
	showing that
	the induced distribution $\pi_\ww$ is close to $\pi$.
\end{Proof}

\begin{Proof}[Canonical Paths and Adjusted Conductance]
We can no longer use the canonical paths method for trees (\cref{res:am:prelims:canonical:trees}) since the weighting does not necessarily give rise to a tree.
The `extra edges'---ie, those corresponding to non-self-loops in $\ww_0$---can only increase the conductance. Intuitively, these cannot harm mixing, but we must establish this carefully.

First, we adjust the proof of the canonical paths method for trees, i.e.\ the deduction of \cref{res:am:prelims:canonical:trees} from \cref{res:am:prelims:canonical:gen}.
We use the same canonical paths $\Gamma$, defined by paths in the \BFS tree $T$.
The bound on the spectral gap $\SG$ does not require a lower bound on $\EC[S]$ for \emph{arbitrary} $S$; rather, it only needs a lower bound on $\EC[T_x]$ for all $x \in V$.

The fact that $E' = \emptyset$ when $\ww_0$ is only self-loops meant that the set of edges emanating from the set $T_x$ was given by
\(
	F(T_x, T_x^c)
=
	\bra{x, \prnt(x)}.
\)
More generally, it is given by
\(
	\bra{x, \prnt(x)} \cup E'(T_x, T_x^c).
\)
But this is always a superset of $\bra{x, \prnt(x)}$, so the edge conductance is always larger than if the edges from $E'$ were ignored.
Motivated by this, define the following adjustment of edge conductance:
\[
	\widehat \Phi_\ww(T_x)
\cq
	\frac
		{\ww\rbr{\bra{x, \prnt(x)}} / \ww(V)}
		{\pi_\ww(T_x) \pi_\ww(T_x^c)}
\Qfor
	x \in V
\Qand
	\widehat \Phi_\ww^\star
\cq
	\MIN{x \in V}
	\widehat \Phi_\ww(T_x).
\]
This is the edge conductance where only the boundary edges in the tree $T = (V, F)$ are considered.
The same proof as for canonical paths for trees then implies that
\[
	\SG[\ww]
\ge
	\widehat \Phi_\ww^\star / \DIAM,
\Quad{recalling that}
	\DIAM
=
	\diam T.
\qedhere
\]
\end{Proof}

\begin{Proof}[Conductance Analysis]
The analysis of the conductance in \cref{res:am:pf:conductance} needs to be adjusted.
We need only analyse the conductance of complete subtrees $T_x$ and must regard the boundary as only $F(T_x, T_x^c) = \bra{x, \prnt(x)}$, not the full boundary $(E' \cup F)(T_x, T_x^c)$.
The proof of \cref{res:am:pf:conductance} applies almost unchanged to control $\widehat \Phi_\ww^\star$:
we obtain
\[
	\widehat \Phi_\ww^\star
\ge
	\tfrac16 \eps / \DIAM.
\]
The only point to be noted is the establishment of the equality $\ww_0(T_x) = \pi(T_x)$.
Previously, this was obvious from the self-loop weightings.
It still holds here, since the invariant distribution induced by $\ww_0$ is $\pi$ and $\ww_0$ has unit total weight, by assumption.
	Thus, in fact, $\ww_0(x) = \pi(x)$ for all $x \in V$.
\end{Proof}

\endgroup

\begin{Proof}[Conclusion]
We combine the two results above, exactly as before, to obtain
\[
	\SG[\ww]
\ge
	\widehat \Phi_\ww^\star / \DIAM
\ge
	\tfrac16 \eps / \DIAM^2
\ge
	\tfrac1{24} \eps / (\diam G)^2.
\]
The conversion of this into a lazy chain and then into a mixing estimate is unchanged.
\end{Proof}

\section{Continuous-Time Markov Chains}
\label{sec:cts}

\subsection{Set-Up, Main Result and Outline}

We have been studying discrete-time Markov chains throughout this paper.
It is natural to ask the same question for continuous-time chains.
Our attention is devoted to continuous-time Markov chains which are reversible \wrt the uniform distribution.
Such chains can always be represented as a \RW on a weighted graph $(G, \ww)$ where $\ww : E \to \mbr_+$ is a weighting on the edge of $G = (V, E)$.

Our main result for continuous-time chains is simple to state:
	we impose a normalisation of $\abs V^{-1} \sum_{e \in E} \ww(e)$,
		ie the average rate at which the \RW leaves a vertex is at most $1$;
	we define a weighting $\ww$ and show an upper bound of order $(\diam G)^2$ on the spectral gap of this \RW.

\begin{thm}[Fast Mixing Continuous-Time Markov Chain]
\label{res:cts:res:main}
	Let $G = (V, E)$ be a graph.
	There exists a weighting $\ww : E \to \mbr_+$ with
	average rate
	\(
		\abs V^{-1} \sumt{x \in V} \ww(x) \le 1
	\)
	and such that
	the Markov chain induced by this weighting has
	spectral gap $\SG[\ww]$ and mixing time $\MIX[\ww][\cdot]$ satisfying
	\[
		\SG[\ww]
	\ge
		\tfrac1{16} (\diam G)^{-2}
	\Qand
		\MIX[\ww][\xi]
	\le
		8 (\diam G)^2 \log\rbr{ \xi^{-2} \abs V }
	\Qforall
		\xi \in (0,1).
	\]
\end{thm}

\begin{Proof}[Proof of \cref{res:cts:res:main}: Outline]
\qedtriangle
The outline is the same as in discrete-time.
%The proof has four key steps.
\begin{enumerate}
	\item 
	We construct a weighted spanning tree;
	see \cref{def:cts:pf:tree}.
	
	\item 
	We control the total weight of the spanning tree;
	see \cref{res:cts:pf:weight}.
	
	\item 
	We estimate the conductance of this weighted tree;
	see \cref{res:cts:pf:conductance}.
	
	\item 
	We relate its spectral gap and conductance using canonical paths;
	see \cref{res:cts:pf:gap}.
\qedhere
\end{enumerate}
\end{Proof}

We fix a graph $G = (V, E)$ and always take $\pi \cq \unif[V]$ to be the uniform distribution on $V$.
We do not always repeat these in statements below.
Also, we write $n \cq \abs V$.

\subsection{Proof via Adjustments to Discrete-Time Case}
\label{sec:cts:pf}

The proof in continuous-time is surprisingly similar to that used in discrete-time.

\begin{itemize}
	\item 
	We construct the same weighted tree $(T, \ww)$, except that we do not include the self-loops;
	contrast \cref{def:am:pf:tree,def:cts:pf:tree}.
	
	\item 
	The invariant distribution of a \RW on a graph with weights on the edges is always uniform.
	Thus we do not need an analogue of \cref{res:am:pf:inv-dist}.
	We require the \emph{total weight} to be at most $n$, instead of requiring the invariant distribution to be close to a given measure.
	
	\item 
	We use the same argument to control the conductance;
	cf \cref{res:am:pf:conductance}.
	The only differences is that now we do not include the self-loop weight in the calculation.
	
	\item 
	The canonical paths argument applies in continuous-time;
	cf \cref{res:am:prelims:canonical:trees,res:am:pf:gap}.
\end{itemize}

First, we define the weighted spanning tree $(T, \ww$); cf \cref{def:am:pf:tree}.

\begin{defn}[Weighted Spanning Tree]
\label{def:cts:pf:tree}
	Let $o \in V$ and let $T = (V, F)$ be a \BFS tree rooted at $o$.
	Define the following weightings $\ww : F \to \mbr_+$ by
	\[
		\ww\rbr{\bra{x, \prnt(x)}}
	\cq
		\tfrac12 \abs{T_x} / \diam T.
	\Qfor
		x \in V \setminus \bra{o}.
	\]
	
	$T$ is a \BFS tree, so $\DIAM \cq \diam T$ satisfies $\DIAM \le 2 \diam G$.
\end{defn}

This is equivalent to the tree-weight in the discrete-time case; see $\ww_1$ in \cref{def:am:pf:tree}.
The particular scaling is chosen so that the total weight of $\ww$ is at most $\DIAM n$, as the next lemma shows.

\begin{lem}[Total Weight of $\ww$]
\label{res:cts:pf:weight}
	We have
	\(
		\ww(V)
	\le
		n.
	\)
\end{lem}

\begin{Proof}
This is an immediate consequence of the subtree counting lemma (\cref{res:am:prelim:counting-trees}):
\[
	\ww(V)
=
	2 \ww(F)
=
	2
	\sumt{x \in V \setminus \bra{o}}
	\ww\rbr{\bra{x, \prnt(x)}}
=
	\sumt{x \in T_o \setminus \bra{o}}
	\abs{T_x} / \DIAM
\le
	\abs{T_o}
=
	n.
\qedhere
\]
\end{Proof}

Next, we control the conductance of this weighted spanning tree; cf \cref{res:am:pf:conductance}.
To do this, we must first give the precise definition of \textit{conductance} in continuous-time.

\begin{defn}[Conductance]
	Let $T = (V, F)$ be a graph and let $\uu : F \to \mbr_+$ be a weighting.
	The \textit{conductance} $\ECRW[\uu][S]$ of a set $S \subseteq V$ with $\pi_\uu(S) > 0$ \wrt $\uu$ is defined to be
	\[
		\ECRW[\uu][S]
	\cq
		\uu\rbb{ F(S, S^c) }
	\big/
		\abs S.
	\]
	The \textit{conductance} $\ECRW*[\uu]$ of $\uu$ is defined to be
	\[
		\ECRW*[\uu]
	\cq
		\MIN{S \subseteq V : 0 < \pi_\uu(S) \le 1/2}
		\ECRW[\uu][S].
	\]
\end{defn}

\begin{prop}[Conductance]
\label{res:cts:pf:conductance}
	Let $(T, \ww)$ be the weighted spanning tree from \cref{def:am:pf:tree}.
	The conductance $\ECRW*[\ww]$ of the \RW on $(T, \ww)$ satisfies
	\[
		\ECRW*[\ww]
	\ge
		\tfrac14 \DIAM^{-1}.
	\]
\end{prop}

\begin{Proof}
The same reductions as used in \cref{res:am:pf:conductance} show that it suffices to show that
\[
	\Phi_\ww(T_x)
\ge
	\tfrac12 \DIAM^{-1}
\Qforall
	x \in V.
\]
But this is immediate from the definition of $\ww$:
\[
	\Phi_\ww(T_x)
=
	\ww\rbr{\bra{x, \prnt(x)}}
/
	\abs{T_x}
=
	\tfrac12 \abs{T_x} \DIAM^{-1}
/
	\abs{T_x}
=
	\tfrac12 \DIAM^{-1}.
\qedhere
\]
\end{Proof}

Finally, we apply the canonical paths method for trees (\cref{res:am:prelims:canonical:trees}) to deduce a bound on the spectral gap for $(T, \ww)$; cf \cref{res:am:pf:gap}.
We must adjust this to apply in continuous-time; see \cref{res:cts:pf:canonical:gen,res:cts:pf:canonical:trees}.

\begin{cor}[Spectral Gap]
\label{res:cts:pf:gap}
	Let $(T, \ww)$ be the weighted spanning tree from \cref{def:am:pf:tree}.
	The spectral gap $\SG[\ww]$ of the \RW on $(T, \ww)$ satisfies
	\[
		\SG{\ww}
	\ge
		\tfrac14 \DIAM^{-2}.
	\]
\end{cor}

\begin{Proof}
This is an immediate consequence of
	the canonical paths method for trees in continuous-time (\cref{res:cts:pf:canonical:trees})
and
	\cref{res:cts:pf:conductance}.
\end{Proof}

It remains to adjust the method of canonical paths to continuous-time.

\begin{prop}[Canonical Paths in Continuous-Time: General]
\label{res:cts:pf:canonical:gen}
	Let $G = (V, E)$ be a graph and $\uu : E \to \mbr_+$ be a weighting.
	Let $\gamma_{x,y}$ be an $F$-path from $x$ to $y$ for all $x,y \in V$.
	The spectral gap $\SG[\uu]$ of the \RW on $(G, \uu)$ satisfies
	\[
		\SG[\uu]
	\ge
		n
		\MIN{e \in E}
		\brb{
			Q_\uu(e)
		\big/
			\sumt{x,y \in V}
			\one{e \in \gamma_{x,y}}
			\abs{\gamma_{x,y}}
		}.
	\]
\end{prop}

\begin{Proof}
The discrete-time case is proved in \cite[Theorem~5]{S:canonical-paths}.
It involves the variational characterisation of the spectral gap in terms of the Dirichlet form. This characterisation holds both in discrete- and continuous-time; see \cite[\S 3.6]{AF:book}.
The proof in \cite{S:canonical-paths} then passes almost unchanged to the continuous-time set-up, recalling that now the invariant distribution is uniform.

Concretely, one can rescale the weights, setting $\widetilde \ww(\cdot) \cq c \ww(\cdot)$ for some value $c$ such that $\max_{x \in V} \widetilde \ww(x) = 1$. This can then be realised by placing mean-$1$ exponential wait times between jumps of a discrete-time chain $P$.
One then applies the canonical paths method to $P$.
The Dirichlet form is \emph{linear} in this scaling meaning that the scaling can be `undone' at the end.
\end{Proof}

\begin{rmkt*}[Cheeger-Inequality in Continuous-Time]
We remark that while a scaling argument as used above does apply for the usual discrete-time Cheeger inequality,
	namely $\SG \ge \tfrac12 (\EC*)^2$,
the bound is \emph{quadratic} in the scaling.
Thus a factor of $c$ is lost.
See \cite[Theorem~4.40]{AF:book}.
In our set-up,
\(
	\max_{x \in V}
	\ww(x)
\)
may be as large as
\(
	(n - 1) / \DIAM.
\)
This would lead to a lower bound of $1/(n \DIAM)$ on the spectral gap, rather than $1/\DIAM^2$ as we were able to achieve using canonical paths.
\end{rmkt*}

Analogous arguments to those used in the special case that $G$ is a tree,
	ie deducing \cref{res:am:prelims:canonical:trees} from \cref{res:am:prelims:canonical:gen},
apply in the continuous-time set-up too.

\begin{cor}[Canonical Paths in Continuous-Time: Trees]
\label{res:cts:pf:canonical:trees}
	Let $T = (V, F)$ be a tree and $\uu : F \to (0, \infty)$ be a weighting.
	The spectral gap $\SG[\uu]$ of the \RW on $(T, \uu)$ satisfies
	\[
		\SG[\uu]
	\ge
		\ECRW*[\uu] / \diam T.
	\]
\end{cor}

We now have all the ingredients required to deduce the main result.

\begin{Proof}[Proof of \cref{res:cts:res:main}]
The Markov chain constructed is a \RW on a weighted tree. It is defined in \cref{def:am:pf:tree}.
\cref{res:cts:pf:weight} bounds the total weight of this tree by $n$, as required.

The spectral gap bound is proved via \cref{res:cts:pf:conductance,res:cts:pf:gap}.
The mixing time bound is then deduced from the spectral gap bound via the (continuous-time) spectral gap--mixing time relation;
see \cite[Lemma~4.23]{AF:book} for details of this relation.
%\footnote{%
%	Calculation for our reference:
%	the proof of \cite[Lemma~4.23]{AF:book} implies that
%	\[
%		d(t) \le \tfrac12 \pimin^{-1/2} e^{-t \SG}.
%	\Quad{Thus}
%		t \ge \tfrac12 \SG^{-1} \log(\tfrac14 \xi^{-2} \pimin^{-1})
%	\implies
%		d(t) \le \xi.
%	\]}
%or \cite[Theorem~20.6]{LPW:markov-mixing}
	%
\end{Proof}

\subsection{Hitting Time of the Root}

The following result is not needed for the proof, but it is a nice little result and its proof is extremely simple,
given a reference regarding hitting times in trees.
Write $\tau_x$ for the hitting time of vertex $x$.

\begin{lem}[Hit Root in Diameter Squared]
\label{res:cts:hit:main}
	Let $G = (V, E)$ be a graph.
	For all $o \in V$,
	there exists a weighting $\ww : E \to \mbr_+$ with
	average rate
	\(
		\tfrac1n
		\sumt{x \in V}
		\ww(x)
	\le
		1
	\)
	and such that
	the Markov chain induced by this weighting has
	worst-case expected hitting time of $o$ satisfying
	\[
		\MAX{x \in V}
		\ext[x]{\tau_o}
	\le
		8 (\diam G)^2.
	\]
\end{lem}

We use the weighted spanning tree $(T, \ww)$ used above, i.e.\ from \cref{def:cts:pf:tree}.
%; this is fundamentally the same as that used in \cref{def:am:pf:tree}, except without the self-loops.
Precisely, the tree $T = (V, F)$ is rooted at some vertex $o \in V$ and
\[
	\ww\rbr{\bra{x, \prnt(x)}}
\cq
	\tfrac12 \abs{T_x} / \diam T
\Qfor
	x \in V \setminus \bra{o}.
\]
%The \RW is thus a birth-and-death process on a tree.
Moments of hitting times of the root in reversible Markov chains on trees were investigated by \textcite{Z:hitting-birth-death-tree}.
The following result is a special case of \cite[Theorem~1.1]{Z:hitting-birth-death-tree}.

\begin{thm}[Hitting Times in Trees; cf {\cite[Theorem~1.1]{Z:hitting-birth-death-tree}}]
\label{res:cts:hit:zhang}
	Let $T = (V, F)$ be a finite tree and $\qq : F \to (0, \infty)$ a weighting.
	Let $o \in V$ and root the tree at $o$; use the notation of \cref{def:am:prelim:tree-not}.
	Let $\tau_o$ denote the hitting time of the root.
	For all $x \in V$,
	we have
	\[
		\ex[x]{\tau_o}
	=
		\sumt{y \in \anc(x) \setminus \bra{o}}
		\abs{T_y} / \qq\rbr{\bra{y, \prnt(y)}}.
	\]
\end{thm}

The hitting time result of \cref{res:cts:hit:main} follows easily from this.

\begin{Proof}[Proof of \cref{res:cts:hit:main}]
Let $(T, \ww)$ be the weighted spanning tree from \cref{def:am:pf:tree}.
We have
\(
	\abs{T_y} / \ww\rbr{\bra{y, \prnt(y)}}
=
	2 \diam T
\)
for all $y$.
Thus applying \cref{res:cts:hit:zhang} to this weighted tree gives
\[
	\ex[x]{\tau_o}
=
	\abs{\anc(x) \setminus \bra{o}}
\cdot
	2 \diam T
\le
	2 (\diam T)^2
\le
	8 (\diam G)^2,
\]
since $\diam T \le 2 \diam G$.
The weighting satisfies
\(
	\ww(V)
\le
	n
\)
by \cref{res:cts:pf:weight}.
\cref{res:cts:hit:main} follows.
\end{Proof}

\begin{rmkt}[Inspiration]
\label{rmk:cts:hit:inspiration}
This result of \textcite{Z:hitting-birth-death-tree} was the inspiration behind our choice of weighted spanning tree in both the discrete- and continuous-time set-ups.
Particularly, the simplicity of the formula when one takes
\(
	\qq\rbr{\bra{x, \prnt(x)}}
\cq
	\abs{T_x}
\)
encouraged us to try this weighting.

We had originally tried to make the distance to the root behave roughly like an unbiased \RW on the integers---this gives the right ``diameter-squared'' bound.
This means balancing the weights on either `side' of a vertex.
This works for some trees, e.g.\ the path, rooted at one end, and the binary tree.
But it does not combine well: attaching a path and binary tree, each of the same depth, at their root gives rise to a ``diameter-cubed'' hitting time of the root.

Although we have phrased the continuous-time proof as an adjustment to the discrete-time one, we actually developed the continuous-time argument first.
Indeed, this is natural because \emph{any} edge weighting gives rise to the uniform distribution in continuous-time, so there is no need to do any superposition with a `base' weighting, such as self-loops.
\end{rmkt}

\section{Time-Inhomogeneous Markov Chains}
\label{sec:tinhom}

\newcommand{\FOR}
	{\textsf{\textbf{for}}\xspace}
\newcommand{\ENDFOR}
	{\textsf{\textbf{endfor}}\xspace}

\newcommand{\hyp}[1]{\{0,1\}^{#1}}

%%%%%%%%%%

The content of this section is somewhat different from the previous ones.
Our desire, as always, is to sample from a distribution $\pi$ on a set $V$ via some Markov process which only uses transitions permitted by the graph $G$.
The difference here is that we use a \textit{time-inhomogeneous} Markov chain.

Markov chains are typically time-homogeneous.
Discrete-time chains are then described by a transition matrix $P$ and an initial law $\mu_0$. The time-$t$ law $\mu_t$ is then given by applying $P$ $t$ times to $\mu$:
\(
	\mu_t = \mu_0 P^t.
\)
Continuous-time chains are described in a somewhat similar manner.
Time-inhomogeneous chains are allowed to use a different transition matrix at each step:
	$P_1$ is used for the first step,
	$P_2$ for the second
and
	so on.
The time-$t$ law $\mu_t$ is then given by
\(
	\mu_t = \mu_0 P_1 \cdots P_t.
\)
The special case where $P_t = P$ for all $t$, for some $P$, reduces to the time-homogeneous case.

Our main result for time-inhomogeneous chains is simple to state:
	given a graph $G = (V, E)$ and $\pi \in \PD{V}$,
	we exhibit a time-inhomogeneous Markov chain which satisfies
	\(
		\mu_{2 \diam G}
	=
		\pi.
	\)

\begin{thm}%[Time-Inhomogeneous]
\label{res:timhom:res:main}
	Let $G = (V, E)$ be a connected graph and let $\pi \in \PD{V}$.
	There exists a time-inhomogeneous Markov chain on $G$ obtaining perfect mixing after $2 \diam G$ steps:
	\(
		\mu_{2 \diam G}
	=
		\pi.
	\)
\end{thm}

We fix a graph $G = (V, E)$ and a probability measure $\pi \in \PD{V}$ throughout this section.
We do not always repeat these in statements below.

\begin{Proof}[Proof of \cref{res:timhom:res:main}]
Choose $o \in V$ and a breadth-first search (\textit{\BFS}) tree $T = (V, F)$ rooted at $o$ arbitrarily.
Given $x \in V$, let $V_x$ denote the set of vertices for which $x$ lies in the unique shortest path to the root $o$; let $T_x \cq T[V_x]$.
Write $\DIAM\rbr{x} \cq \depth T_x$ for the \emph{depth} of $T_x$, i.e.\ the maximal distance $\max_{y \in T_x} \dist(y,x)$ to the \emph{root}.
We construct the time-inhomogeneous chain inductively

Suppose that $X_0 = o$; we cover the general case later.
We claim that there exists a chain $X^o$ on $T_o$ such that $X^o_{\DIAM\rbr{o}} \sim \fnrestrict{\pi}{V_o}$.
This is trivially true if $T_o$ is a singleton, which has depth $0$.
Assume now that $\abs{T_o} \ge 2$.
Define the first transition matrix $P_1$ to keep $X$ at $o$ with probability $\pi(o)$ and otherwise move to $x$ with $\bra{o,x} \in E$ with probability proportional to $\pi{T_x}$:
\[
	P_1(o, o)
\cq
	\pi(o)
\Qand
	P_1(o, x)
\cq
	\pi(T_x)
\Qfor
	x \in V \text{ with } \bra{o, x} \in V.
\]
Let $P_k(o,o) = 1$ for all $k \ge 2$. Thus if $X_1 = o$ then $X_k = o$ for all $k \ge 1$.

There exists a chain $X^x$ on $T_x$ such that $X^x_{\DIAM\rbr{x}} \sim \fnrestrict{\pi}{V_x}$ for each $x \in V$, by the inductive hypothesis.
Importantly, $T_x \cap T_y = \emptyset$ if $\bra{o,x}, \bra{o,y} \in F$ and $x \ne y$.
We can thus define sequences $P_x \cq (P^x_k)_{k=1}^{\DIAM\rbr{x}}$ of transition matrices for each $x \in V$ with $\bra{o, x} \in E$ and each $P^x$ defined on a disjoint set by induction.
Define $P^x_k \cq I$ for $k > \DIAM\rbr{x}$.
These can then be combined into a single sequence $(P_k)_{k=2}^{\DIAM\rbr{o}}$:
	if $X^o_1 = x$, then we use sequence $P^x$.
Then
\[
	X^o_{\DIAM\rbr{o}} = X^x_{\DIAM\rbr{x}} \sim \fnrestrict{\pi}{V_x}
\Quad{conditional on}
	X^o_1 = x \ne o,
\]
noting that $\DIAM\rbr{x} = \DIAM\rbr{o} - 1$.
Thus
\[
	\pr[o]{ X^o_{\DIAM\rbr{o}} = y }
&
=
	\sumt{x \in V : \bra{o,x} \in E}
	\pr[x]{ X^x_{\DIAM\rbr{x}} = y } \pr[o]{ X^o_1 = x }
\\&
=
	\sumt{x \in V : \bra{o,x} \in E}
	\rbb{ \one{y \in V_x} \fnrestrict{\pi}{V_x}(y) } \cdot \pi(V_x)
\\&
=
	\sumt{x \in V : \bra{o,x} \in E}
	\one{y \in V_x} \cdot \rbb{ \pi(y) / \pi(V_x) } \cdot \pi(V_x)
\\&
=
	\pi(y)
	\sumt{x \in V : \bra{o,x} \in E}
	\one{y \in V_x}
=
	1/\abs{T_o}.
\]

This argument is no more than a formalisation of the following informal verbal description.
\begin{itemize}
	\item 
	If the walk is at $x$, then stay at $x$ with probability $\pi(x) / \pi(V_x)$.
	
	\item 
	Otherwise, move to the children of $x$ with probability proportional to their $\pi$-measure.
	
	\item 
	If at some point the walk stays put, then keep it at the state indefinitely.
\end{itemize}

This completes the argument when $X_0 = o$.
It remains to consider the case that $X_0 \ne o$.
Direct all edges towards $o$ and run for $\DIAM\rbr{o}$ steps.
Precisely, set
\[
	P_0(x,y)
\cq
	\one{x \ne o, \: y = \prnt(x)}
+	\one{x = y = o},
\]
where $\prnt(x)$ is the unique neighbour of $x \ne o$ on the unique shortest path from $x$ to $o$.
If $\DIAM\rbr{o}$ steps are made according to this matrix, then $X_{\DIAM\rbr{o}} = o$, regardless of $X_o$.
We then apply the construction from the case $X_0 = o$, all shifted by $\DIAM\rbr{o}$.

In summary, we obtain a time-inhomogeneous Markov chain $X$ with $X_{2 \DIAM\rbr{o}} \sim \pi$. Finally, $\DIAM\rbr{o} = \depth T \le \diam G$.
Keeping $X$ fixed during $( 2 \DIAM\rbr{o}, \: 2 \diam G ]$ gives $X_{2 \diam G} \sim \pi$.
\end{Proof}

\section{Open Problems and Concluding Remarks}
\label{sec:conc}

We have studied fundamental barriers to fast mixing on graphs. 
There are a few questions left open by our work.
First, throughout this paper we have been mainly focussed on mixing to the uniform distribution. 
Whilst this is arguably the most important case, we believe that generalising our results to arbitrary distributions would be an interesting extension to our work.
In particular, while extending our results on continuous-time chains should not require much effort, more thought is needed in generalising our Cheeger-type inequality between the vertex conductance and the optimal spectral gap.
It is possible to generalise the notion of vertex conductance to any arbitrary positive distribution $\pi$.
We believe that it should be possible to carry over the general idea of our proof to this case by replacing matching conductance with a similar conductance measure based on a generalisation of the weighted vertex cover problem.
The main difficulty then lies in generalising the proof of \cref{res:vc:fm:coarea}, since it crucially depends on the greedy algorithm for maximal matching.

Another open problem prompted by our work is to construct a graph sparsifier, 
i.e., an edge-induced sparse subgraph,
that approximately preserves the vertex conductance of the original graph. 
Our work together with a result by \textcite{BSS} implies that it is possible to construct an order $n$-size sparsifier of $G$ with vertex conductance at most order $\VC*[G]^2 \log n$.
It is then natural to ask if we can obtain a better approximation.
 
Finally, can our results spur new algorithmic applications? We mention two.
First, we would like to design a distributed algorithm to compute a fast mixing Markov chain on $G$, where $G$ also represents the topology of the distributed network.
Second, we ask if it is possible to design a local algorithm, in the spirit of \textcite{ST:nibble}, that outputs a subset of nodes with small vertex conductance.

\renewcommand{\bibfont}{\sffamily\small}
\printbibliography[heading = bibintoc, title = {Bibliography}]

\end{document}